\input amstex\documentstyle{amsppt}  
\pagewidth{12.5cm}\pageheight{19cm}\magnification\magstep1
\topmatter
\title A class of perverse sheaves on a partial flag manifold\endtitle
\author G. Lusztig\endauthor
\address{Department of Mathematics, M.I.T., Cambridge, MA 02139}\endaddress
\thanks{Supported in part by the National Science Foundation.}\endthanks
\endtopmatter   
\document

\define\doa{\dot a}
\define\dda{\ddot{a}}

\define\dw{\dot w}

\define\ds{\dot s}

\define\dZ{\dot Z}
\define\bcl{\bar{\cl}}

\define\uB{\un B}
\define\uT{\un T}
\define\uM{\un M}

\define\utB{\un{\ti B}}

\define\dts{\dot{\ts}}
\define\mpb{\medpagebreak}
\define\Up{\Upsilon}

\define\be{\bar e}

\define\bx{\bar x}

\define\ha{\hat a}
\define\hs{\hat s}

\define\dsv{\dashv}

\define\po{\text{\rm pos}}

\define\si{\sim}

\define\sqc{\sqcup}
\define\ovs{\overset}
\define\qua{\quad}

\define\tcl{\ti\cl}

\define\tiz{\ti\z}

\define\tbcl{\bar{\tcl}}

\define\bY{\bar Y}
\define\bX{\bar X}
\define\bZ{\bar Z}

\define\lb{\linebreak}

\define\eSb{\endSb}

\define\op{\oplus}

\define\part{\partial}
\define\em{\emptyset}
\define\imp{\implies}

\define\n{\notin}
\define\iy{\infty}
\define\m{\mapsto}
\define\do{\dots}

\define\bsl{\backslash}

\define\lra{\leftrightarrow}

\define\sm{\smallmatrix}
\define\esm{\endsmallmatrix}
\define\sub{\subset}    
\define\bxt{\boxtimes}
\define\T{\times}
\define\ti{\tilde}
\define\nl{\newline}
\redefine\i{^{-1}}

\define\un{\underline}

\define\ot{\otimes}
\define\bbq{\bar{\QQ}_l}

\define\Hom{\text{\rm Hom}}
\define\End{\text{\rm End}}
\define\Aut{\text{\rm Aut}}

\define\supp{\text{\rm supp}}

\define\bst{\bigstar}

\define\a{\alpha}
\redefine\b{\beta}
\redefine\c{\chi}
\define\g{\gamma}
\redefine\d{\delta}
\define\e{\epsilon}

\define\io{\iota}
\redefine\o{\omega}
\define\p{\pi}
\define\ph{\phi}

\define\r{\rho}

\redefine\t{\tau}
\define\th{\theta}
\define\k{\kappa}
\redefine\l{\lambda}
\define\z{\zeta}
\define\x{\xi}

\define\vt{\vartheta}

\redefine\G{\Gamma}
\redefine\D{\Delta}

\redefine\L{\Lambda}
\define\Ph{\Phi}

\redefine\aa{\bold a}

\define\boc{\bold c}

\define\kk{\bold k}

\redefine\ss{\bold s}
\redefine\tt{\bold t}
\define\uu{\bold u}
\define\vv{\bold v}
\define\ww{\bold w}

\define\FF{\bold F}

\define\II{\bold I}

\define\NN{\bold N}

\define\QQ{\bold Q}

\define\SS{\bold S}

\define\ZZ{\bold Z}

\define\ca{\Cal A}
\define\cb{\Cal B}

\define\cd{\Cal D}
\define\ce{\Cal E}
\define\cf{\Cal F}
\define\cg{\Cal G}

\define\ci{\Cal I}
\define\cj{\Cal J}

\define\cl{\Cal L}
\define\cm{\Cal M}

\define\cp{\Cal P}

\define\cs{\Cal S}
\define\ct{\Cal T}

\define\cz{\Cal Z}
\define\cx{\Cal X}
\define\cy{\Cal Y}

\define\fd{\frak d}

\define\fA{\frak A}

\define\fD{\frak D}

\define\fF{\frak F}

\define\fP{\frak P}

\define\fT{\frak T}

\define\fY{\frak Y}

\define\tb{\ti b}

\define\td{\ti d}

\define\tf{\ti f}
\define\tg{\ti g}

\define\tj{\ti j}

\define\tl{\ti l}
\define\tm{\ti m}

\define\tio{\ti\o}

\define\ts{\ti s}
\define\tit{\ti t}

\define\tv{\ti v}
\define\tw{\ti w}

\define\ty{\ti y}

\define\tB{\ti B}
\define\tC{\ti C}

\define\tM{\ti M}

\define\tS{\ti S}

\define\tX{\ti X}
\define\tY{\ti Y}
\define\tZ{\ti Z}

\define\sh{\sharp}

\define\sps{\supset}

\define\bul{\bullet}
\define\uZ{\un Z}
\define\che{\check}
\define\cha{\che{\a}}

\define\ucl{\un\cl}

\define\tss{\ti{\ss}}
\define\tSS{\ti{\SS}}
\define\ciss{\ci_\ss}
\define\cits{\ci_{\tss}}
\define\Upss{\Up^\ss}
\define\Uptss{\Up^{\tss}}
\define\tcj{\ti{\cj}}
\define\tcp{\ti{\cp}}
\define\tcf{\ti{\cf}}
\define\tcb{\ti{\cb}}
\define\tcy{\ti{\cy}}
\define\y{\ti r}
\define\ttir{\ti t_{\y}}

\define\sscj{\ss_\cj}
\define\tsstcj{\tss_{\tcj}}
\define\ccl{\check{\cl}}
\define\uucl{\un{\ucl}}
\define\bUp{\bar{\Up}}
\define\tip{\ti\p}
\define\chR{\check R}
\define\BE{Be}
\define\BBD{BBD}
\define\BO{Bo}
\define\DL{DL}
\define\ORA{L1}
\define\CS{L2}
\define\GR{L3}
\define\HB{L4}
\define\PAR{L5}
\define\CSD{L6}

\head Introduction\endhead
Let $f:X@>>>Y$ be a morphism of algebraic varieties over an algebraically closed field 
$\kk$. The theory of perverse sheaves \cite{\BBD} associates to $f$ a collection of 
invariants, namely the collection of simple perverse sheaves on $Y$ which appear as
subquotients of the $l$-adic perverse sheaf $\op_{j\in\ZZ}{}^pH^j(f_!\bbq)$ on $Y$. More 
precisely, if $f$ is equivariant for given actions of a finite group $\G$ on $X$ and $Y$ 
then we denote by $C_\G(f)$ the (finite) collection of simple $\G$-equivariant perverse 
sheaves on $Y$ which appear as subquotients of the $l$-adic $\G$-equivariant  perverse 
sheaf $\op_{j\in\ZZ}{}^pH^j(f_!\bbq)$ on $Y$. (The perverse sheaves in $C_\G(f)$ are not
necessarily simple if the $\G$-equivariant structure is disregarded.) 

In this paper we try to understand the collection $C_\G(f)$ in the case where:

$\kk$ is an algebraic closure of a finite field $\FF_q$, 

$G$ is a connected reductive algebraic group defined over $\FF_q$ with Frobenius map 
$F:G@>>>G$, 

$\G$ is the group of $\FF_q$-rational points of $G$, 

$X$ is the variety of all pairs $(B',g)$ where $B'$ is a Borel subgroup of $G$ and $g\in G$
is such that $g\i F(g)$ is in the unipotent radical of $B'$,

$Y$ is the variety of parabolic subgroups of $G$ of a fixed type,

$f$ associates to $(B',g)\in X$ the unique parabolic subgroup in $Y$ that contains $B'$. 
\nl
(The action of $\G$ on $X$ is $g_0:(B',g)\m(g_0B'g_0\i,gg_0\i)$; the action of $\G$ on $Y$
is by conjugation.)

Note that the variety $X$ can be viewed as a family of varieties of the type considered in
\cite{\DL} indexed by the full flag manifold.

In this paper we define a finite collection $\Bbb S(Y)$ of simple $\G$-equivariant perverse
sheaves on $Y$ by two methods (see Sections 3 and 4). These methods and the proof of their
equivalence are similar to those used in the theory of parabolic character sheaves 
\cite{\PAR}. The second method (see Section 4) gives a description of these perverse 
sheaves in terms of some explicit local systems on some pieces of a finite partition of 
$Y$. This partition, introduced by the author in 1977, and further studied in \cite{\BE},
reduces in the case where $Y$ is the full flag manifold to the partition introduced in 
\cite{\DL}.

In 7.6 we show that $C_f(Y)\sub\Bbb S(Y)$.

In Section 6 we construct an explicit basis for the space of intertwining operators 
between certain cohomologically induced representations of $\G$, extending an idea of
\cite{\HB}. As a bi-product we obtain a disjointness theorem for the objects of $\Bbb S(Y)$
which in the special case where $Y=\{G\}$ reduces to the disjointness theorem 
\cite{\DL, 6.2, 6.3} (but the present proof is quite different from that of \cite{\DL}).

In Section 7 we study the variety $X$ (see above). In particular we show (using results in
Section 6) that $X$ is connected if $G$ is simply connected.

In Section 8 we give a conjecture (based on results in Section 6 and some combinatorial
results in Section 5) which should explain in an intrinsic way the "Jordan decomposition" 
\cite{\ORA} for irreducible representations of $\G$.

\head Contents\endhead
1. Preliminaries.

2. The variety $\cz^\ss$ and the local system $\bcl$.

3. The class $\Bbb S'(\cp_J)$ of simple objects in $\cm_\G(\cp_J)$.

4. The class $\Bbb S(\cp_J)$ of simple objects in $\cm_\G(\cp_J)$.

5. Some computations in the Weyl group.

6. A basis for a space of intertwining operators.

7. The variety $X$.

8. A conjecture.

\head 1. Preliminaries\endhead
\subhead 1.1\endsubhead
Let $\kk$ be an algebraically closed field. In this paper all algebraic varieties are over 
$\kk$. Let $G$ be a connected reductive algebraic group. Let $\cb$ be the variety of Borel
subgroups of $G$. We fix $B\in\cb$ and a maximal torus $T$ of $B$. Let $N(T)$ be the 
normalizer of $T$ in $G$. Let $W=N(T)/T$. Note that $W$ acts on $T$ by conjugation; we use
this action to identify $W$ with a subgroup of $\Aut(T)$. For any $(B',B'')\in\cb\T\cb$ 
there is a unique $w\in W$ such that $gB'g\i=B,gB''g\i=\tw B\tw\i$ for some $g\in G$ and 
some representative $\tw$ of $w$ in $N(T)$. We then write $\po(B',B'')=w$. Let 
$l:W@>>>\NN$ be the length function: $l(w)=\dim(B'/(B'\cap B''))$ where 
$(B',B'')\in\cb\T\cb$, $\po(B',B'')=w$. Let $\II=\{w\in W;l(w)=1\}$. 

If $H$ is a group acting on a set $X$ we denote by $X^H$ the fixed point set of $H$ on $X$.

\subhead 1.2\endsubhead
Let $\le$ be the standard partial order on $W$ regarded as a Coxeter group with generators
$\II$. If $X$ is a subset of $W$ and $w\in X$ we write $w=\min X$ if $l(w)<l(w')$ for all 
$w'\in W-\{w\}$. For $J\sub\II$ let $W_J$ be the subgroup of $W$ generated by $J$. For 
$J,J'\sub W$ let ${}^JW=\{w\in W;w=\min(W_Jw)\}$, $W^{J'}=\{w\in W;w=\min(wW_{J'})\}$, 
${}^JW^{J'}={}^JW\cap W^{J'}$.

Let $w_\II$ be the unique element of maximal length in $W$. 

\subhead 1.3\endsubhead
Let $\cp$ be the variety of parabolic subgroups of $G$. For any $P\in\cp$ let $U_P$ be the
unipotent radical of $P$. We set $U=U_B$. For $J\sub\II$ let $P_J\in\cp$ be the subgroup of
$G$ generated by $B$ and by representatives in $N(T)$ of the various elements of $J$. Let 
$L_J$ be the unique Levi subgroup of $P$ that contains $T$. For $s\in\II$ we write $P_s$ 
instead of $P_{\{s\}}$. For $J\sub\II$ let $\cp_J$ be the $G$-conjugacy class of parabolic
subgroups of $G$ that contains $P_J$. For $B'\in\cb$ let $P_{B',J}$ be the unique subgroup
in $\cp_J$ that contains $B'$. For $P\in\cp_J$, $Q\in\cp_{J'}$ the element 
$\po(P,Q):=\min\{w\in W;w=\po(B',B'')\text{ for some }B'\sub P,B''\sub Q\}$ is well defined
and $\po(P,Q)\in{}^JW^{J'}$. We set 

(a) $P^Q=(P\cap Q)U_P\in\cp_{J\cap uJ'u\i}$ where $u=\po(P,Q)$.
\nl
For any $g\in G$ we define $k(g)\in N(T)$ by $g\in Uk(g)U$.

\subhead 1.4\endsubhead
Let $R\sub\Hom(T,\kk^*)$ be the set of roots. Let $\chR\sub\Hom(\kk^*,T)$ be the set of
coroots. Let $\cha\lra\a$ be the standard bijection $\chR\lra R$. For $\a\in R$ let $U_\a$
be the one dimensional root subgroup (normalized by $T$) corresponding to $\a$. Let 
$R^+=\{a\in R;U_\a\sub B\}$, $R^-=R-R^+$. For $s\in\II$ let $\a_s$ be the unique root such
that $U_{\a_s}\sub P_s$. We write $\cha_s$ instead of $\check{\a_s}$. The natural action of
$W$ on $T$ induces an action of $W$ on $R$.

\subhead 1.5\endsubhead
For any $s\in\II$ we fix a homomorphism $h_s:SL_2(\kk)@>>>G$ such that

$h_s\left(\sm b&0\\0& b\i\esm\right)=\cha_s(b)$ for all $b\in\kk^*$,

$a\m h_s\left(\sm 1&a\\0& 1\esm\right)$ is an isomorphism $x_s:\kk@>\si>>U_{\a_s}$,

$a\m h_s\left(\sm 1&0\\a& 1\esm\right)$ is an isomorphism $y_s:\kk@>\si>>U_{\a_s\i}$.
\nl
We say that $\{B,T,h_s(s\in\II)\}$ is an {\it \'epinglage} of $G$. Let 
$\ds=h_s\left(\sm 0&1\\-1&0\esm\right)\in N(T)$. For any $w\in W$ we set 
$\dw=\ds_1\ds_2\do\ds_n\in N(T)$ where $s_1,s_2,\do,s_n\in\II$ are chosen so that
$w=s_1s_2\do s_n$, $l(w)=n$. (This is independent of the choice.) In particular, 
$\dot 1=1$.

For any sequence $\ww=(w_1,w_2,\do,w_r)$ in $W$ we set $[\ww]=w_1w_2\do w_r\in W$ and
$[\ww]^\bul=\dw_1\dw_2\do\dw_r\in N(T)$. 

\subhead 1.6. Equivariant structures \endsubhead
If $X$ is an algebraic variety we write $\cd(X)$ for the derived category of bounded
constructible $\bbq$-sheaves on $X$. If $K\in\cd(X)$ let $\fD(K)\in\cd(X)$ be the Verdier
dual of $X$. Let $\cm(X)$ be the subcategory of $\cd(X)$ whose objects are perverse 
sheaves. If $\ce$ is a local system on $X$ we denote by $\check\ce$ the dual local system.
For $K\in\cd(X)$ we write ${}^pH^{\cdot}(K)$ instead of $\op_j{}^pH^j(K)$. If $f:X@>>>Y$ is
a smooth morphism between algebraic varieties with connected fibres of dimension $\d$, we 
set $f^\bst A=f^*A[\d]$ for any $A\in\cd(Y)$. 

Let $m:H\T Y@>>>Y$ be an action of an algebraic group $H$ on an algebraic variety $Y$. For
any connected component $C$ of $H$ define $m_C:C\T Y@>>>Y$ by $(g,y)\m gy$ and 
$\p_C:C\T Y@>>>Y$ by $(g,y)\m y$. Let $K\in\cm(Y)$. An $H$-equivariant structure on $K$ is
a collection of isomorphisms $\ph_C:\p_C^\bst K@>\si>>m_C^\bst K$ (one for each $C$) such 
that for any two connected components $C,C'$ of $H$, 
$$(\tm\bxt1)^\bst(\ph_{CC'}):
(\tm\bxt1)^\bst(\p_{CC'}^\bst K)@>>>(\tm\bxt1)^\bst(m_{CC'}^\bst K)$$ 
is equal to the composition
$$\align&(\tm\bxt1)^\bst(\p_{CC'}^\bst K)@>\ovs(1)\to=>>\tip^\bst(\p_{C'}^\bst K)@>\b>>
\tip^\bst(m_{C'}^\bst K)\\&
@>\ovs(3)\to=>>(1\T m_{C'})^\bst(\p_C^\bst K)@>\g>>
(1\T m_{C'})^\bst(m_C^\bst K)@>\ovs(2)\to=>>(\tm\bxt1)^\bst(m_{CC'}^\bst K)\endalign$$
where

$\tm:C\T C'@>>>CC'$ is $(g,g')\m gg'$, $\tip:C\T C'\T Y@>>>C'\T Y$ is $(g,g',y)\m(g',y)$;

the equality $(1)$ comes from $\p_{C'}\tip=\p_{CC'}(\tm\T1):C\T C'\T Y@>>>Y$;

the equality $(2)$ comes from $m_{CC'}(1\T m_{C'})=m_{CC'}(\tm\T1):C\T C'\T Y@>>>Y$;

the equality $(3)$ comes from $m_{C'}\tip=\p_C(1\T m_{C'}):C\T C'\T Y@>>>Y$;

$\b=\tip^\bst\ph_{C'}$, $\g=(1\T m_{C'})^\bst\ph_C$.
\nl
The notion of $H$-equivariant structure on a local system on $Y$ is defined in the same way
(but replace $()^\bst$ by $()^*$). 

We denote by $\cm_H(Y)$ the (abelian) category of perverse 
sheaves on $Y$ endowed with an $H$-equivariant structure; the morphisms are morphisms of 
perverse sheaves which are compatible with the equivariant structures. Note that any object
of $\cm_H(Y)$ has finite length. If $K\in\cm_H(Y)$ is semisimple as an object of $\cm(Y)$ 
and $H$ is finite then $K$ is semisimple as an object of $\cm_H(Y)$.

Let $f:X@>>>Y$ be a morphism compatible with $H$-actions on $X$ and $Y$. If $K\in\cm_H(Y)$
then ${}^pH^i(f^*K)$ is naturally an object of $\cm_H(X)$. If $K'\in\cm_H(X)$ then 
${}^pH^i(f_!K')$ is naturally an object of $\cm_H(Y)$. Similarly, if $\ce$ is a local 
system on $Y$ with a given $H$-equivariant structure then the local system $f^*\ce$ on $X$
has an induced $H$-equivariant structure.

If $Y_0$ is a locally closed smooth $H$-stable subvariety of pure dimension of $Y$ and 
$\cf$ is an $H$-equivariant local system on $Y_0$ we denote by $\cf^\sh$ the complex 
$IC(\bY_0,\cf)$ (where $\bY_0$ is the closure of $Y_0$ in $Y$) extended by $0$ on 
$Y-\bY_0$. Then $\cf^\sh[\dim Y_0]$ is an object of $\cm_H(Y)$.

If $K$ is a simple object of $\cm_H(Y)$ then there exists $Y_0,\cf$ as above so that the 
connected components of $Y_0$ are permuted transitively by the $H$-action and 
$K=\cf^\sh[\dim Y_0]$. Note that $K$ is not necessarily simple as an object of $\cm(Y)$. If
$K,K'$ are objects of $\cm_H(Y)$ we write $K\dsv_\G K'$ if $K$ is isomorphic to a 
subquotient of $K'$.

If $K\in\cm_H(Y)$ then $\fD(K)$ is naturally an object of $\cm_H(Y)$.

If $H$ is connected then $K$ has at most one $H$-equivariant structure.

Now assume that $H$ is finite. Let $m_g:Y@>>>Y,y\m gy$. An $H$-equivariant structure on $K$
is a collection of isomorphisms $\ph_g:K@>>>m_g^*K$ (one for each $g\in H$) such that for 
any $g,g'$ in $H$, $K@>\ph_{gg'}>>m_{gg'}^*K$ is equal to the composition
$K@>\ph_{g'}>>m_{g'}^*K@>m_{g'}^*\ph_g>>m_{gg'}^*K$. The same definition applies to
$H$-equivariant structures on local systems on $Y$. 

\mpb

Assume now that $A,A'$ are two simple objects of $\cm_H(Y)$. Let $S=\supp(A)$, 
$S'=\supp(A')$. Note that the irreducible components of $S$ (resp. $S'$) are permuted 
transitively by $H$ hence $S$ (resp. $S'$) has pure dimmension $d$ (resp. $d'$). We show:

(a) {\it $H^0_c(Y,A\ot A')^H=\bbq(-d)$ if $A'=\fD(A)$ in $\cm_H(Y)$ and
$H^0_c(Y,A\ot A')^H=0$ if $A'\not\cong\fD(A)$ in $\cm_(Y)$.}
\nl
Note that there exists an open dense smooth $H$-stable subset $S_0$ (resp. $S'_0$)
of $S$ (resp. $S'$) and an $H$-equivariant local system $\ce$ (resp. $\ce'$) on $S_0$
(resp. $S'_0$) such that $A=\ce^\sh[d]$, $\fD(A)=\check\ce^\sh[d]$, $A'=\ce'{}^\sh[d']$.
By an argument as in \cite{\CS, II, 7.4} we see that $H^0_c(Y,A\ot A')=0$ if $S\ne S'$.
In the rest of the proof we assume that $S=S'$. Then we can also assume that
$S_0=S'_0$. Again by the argument in \cite{\CS, II, 7.4} we see that 
$H^0_c(Y,A\ot A')=H^{2d}_c(S_0,\ce\ot\ce')$. By Poincar\'e duality this can be
identified with $H^0(S_0,\check\ce\ot\check\ce')^*(-d)$. 
It remains to describe $H^0(S_0,\check\ce\ot\check\ce')^H$. This is
the vector space of homomorphisms of 
local systems $\ce'@>>>\check\ce$ which are compatible with the $H$-equivariant
structures. This is the same as the space of morphisms from $A'$ to $\fD(A)$ in
$\cm_H(Y)$ which is $\bbq$ if $A'=\fD(A)$ and is $0$ if $A'\not\cong\fD(A)$.

\mpb

Let $E$ be a finite dimensional $\bbq$-vector space and let $r:H@>>>GL(E)$ be a 
homomorphism. We regard $E$ as an $H$-equivariant local system over a point
in an obvious way. If $X$ is an algebraic variety with $H$-action, we denote by 
$\e:X@>>>\text{point}$ the obvious map and we set $E_X=\e^*E$; this is naturally
an $H$-equivariant local system on $X$ (since $\e$ is compatible with the $H$-action
on $X$ and the trivial $H$-action on the point).

\subhead 1.7\endsubhead
For any torus $T'$ let $\cs(T')$ be the category whose objects are the local systems of 
rank $1$ on $T'$ that are equivariant for the transitive $T'$-action $z:t\m z^nt$ on $T'$ 
for some $n\in\ZZ_{>0}$ invertible in $\kk$.

\subhead 1.8\endsubhead
Let $f:T'@>>>T''$ be a morphism of tori and let $\cl\in\cs(T'')$. We show that the 
following two conditions are equivalent:

(i) $f^*\cl\cong\bbq$;

(ii) $\cl$ is equivariant for the $T'$-action $t':t''\m f(t')t''$ on $T''$. 
\nl
We can find $\k\in\Hom(T'',\kk^*)$ and $\ce\in\cs(\kk^*)$ such that $\cl\cong\k^*(\ce)$. If
the result holds for $\k f,\ce$ instead of $f,\cl$ then it also holds
for $f,\cl$. Thus we may assume that $T''=\kk^*$. We can assume that 
$T'\ne\{1\}$. Let $T_0$ be a codimension $1$ subtorus of $T'$ contained in $\ker f$. Then
$f$ induces a homomorphism $f':T'/T_0@>>>T''$. If the result holds for $f',\cl$ instead of
$f,\cl$ then it also holds for $f,\cl$. Thus we may assume that $T'=T''=\kk^*$. In this 
case the result is immediate.

\subhead 1.9\endsubhead
In this subsection we assume that $\kk$ is an algebraic closure of a finite field $\FF_q$.
Let $F':T@>>>T$ be the Frobenius map for some $\FF_q$-rational structure on the torus $T$.
Let $T^{F'}=\{t\in T;F'(t)=t\}$. The following three sets coincide:

(i) the set of $\cl\in\cs(T)$ (up to isomorphism) such that $\cl$ is $T$-equivariant for 
the $T$-action $t_0:t\m t_0tF'(t_0)\i$ on $T$; 

(ii) the set of $\cl\in\cs(T)$ (up to isomorphism) such that $F'{}^*\cl\cong\cl$;

(iii) the set of $\cl\in\cs(T)$ (up to isomorphism) such that $\cl$ is a direct summand of
$L_!\bbq$ where $L:T@>>>T$ is $t\m tF'(t\i)$;
\nl
moreover they are in a natural bijection with 

(iv) the set of characters $\Hom(T^{F'},\bbq^*)$.
\nl
For $\cl\in\cs(T)$ we have $L^*\cl\cong\cl\ot F'{}^*(\ccl)$. Hence $\cl$ satisfies (ii) if
and only if it satisfies $L^*\cl\cong\bbq$. The last condition is clearly equivalent to the
condition in (iii) and by 1.8 it is also equivalent to the condition in (i).

If $L$ is as in (iii), $T^{F'}$ acts in an obvious way on $L_!\bbq$ and we have
$L_!\bbq=\op_{\c\in\Hom(T^{F'},\bbq^*)}L_!^\c\bbq$ where $L_!^\c\bbq$ is the subsheaf on 
which $T^{F'}$ according to $\c$. It is clear that $\c\m L_!^\c\bbq$ defines a bijection
between the sets (iv) and (iii). The inverse of this bijection can be described as
follows. Let $\cl$ be as in (i). By restriction of the equivariant $T$-structure we obtain
an equivariant $T^{F'}$-structure on $\cl$. Since $T^{F'}$ acts trivially on $T$, it acts
naturally on the stalk of $\cl$ at $1$; this action is via a character 
$\c_\cl:T'{}^F@>>>\bbq^*$. Now $\cl\m\c_\cl$ is the inverse of the bijection above.

\subhead 1.10\endsubhead
Let $\cl\in\cs(T)$. Let $R_\cl=\{\a\in R;\cha^*\cl\cong\bbq\}$. Then $R_\cl$ is a root 
system and $R^+_\cl=R_\cl\cap R^+$ is a set of positive roots for $R_\cl$. Let $\Pi_\cl$ be
the unique set of simple roots of $R_\cl$ such that $\Pi_\cl\sub R^+_\cl$. Let $W_\cl$ be 
the subgroup of $W$ generated by the reflections with respect to the roots in $R_\cl$. Let
$\II_\cl$ be the set of reflections with respect to the roots in $\Pi_\cl$. Then
$(W_\cl,\II_\cl)$ is a Coxeter group. Let $\chR_\cl=\{\cha;\a\in R_\cl\}$.

\head 2. The variety $\cz^\ss$ and the local system $\bcl$\endhead
\subhead 2.1\endsubhead
In this and the next subsection we assume that $\II$ consists of a single element $s$. Let
$a=\cha_s:\kk^*@>>>T$. To any $\cl\in\cs(T)$ such that 

(a) $a^*\cl\cong\bbq$
\nl
we will associate a local system $\ucl$ of rank $1$ on $G$.

{\it Case 1. $a$ is an imbedding.} We have a diagram 
$T@>c>>T/a(\kk^*)@>d>>G/G_{der}@<e<<G$
where $c,e$ are the obvious maps and $d$ is induced by the inclusion $T\sub G$; note that 
$d$ is an isomorphism. Now $\kk^*$ acts on $T$ by $x:t\m a(x)t$ and on $T/a(\kk^*)$ 
trivially; $c$ is compatible with the $\kk^*$-actions. From (a) we see that $\cl$ is 
$\kk^*$-equivariant. Since $\kk^*$ acts freely on $T$ there is a well defined local system
$\cl_1$ of rank $1$ on $T/a(\kk^*)$ such that $\cl=c^*\cl_1$. We set 
$\ucl=e^*(d\i)^*\cl_1$.

{\it Case 2. $a$ is not an imbedding.} Then the centre $\cz$ of $G$ is connected, the 
obvious homomorphism $G_{der}\T\cz@>>>G$ is an isomorphism and we can identify
$G_{der}=PGL_2(\kk)$ compatibly with the standard \'epinglages. Thus we can identify
$G=PGL_2(\kk)\T\cz$. Let $G'=GL_2(\kk)\T\cz$ and let $\p:G'@>>>G$ be the obvious
homomorphism. Let $K=\ker\p$, a one dimensional torus. Let $T'=\p\i(T)$. Let $\p_0:T'@>>>T$
be the restriction of $\p$. Let $a':\kk^*@>>>T'$ be the coroot of $G'$ such that 
$\p_0a'=a$. Let $\cl'=\p_0^*\cl$. Note that $a'{}^*\cl'\cong\bbq$. Applying the
construction in Case 1 to $G',T',\cl'$ instead of $G,T,\cl$ we obtain a local system
$\ucl'$ on $G'$. Now $K$ acts on $T',G',G'/G'_{der},T'/a'(\kk^*)$ by translation, the
analogues of $c,d,e$ for $G'$ are compatible with the $K$-action and $\cl'$ is 
$K$-equivariant. Hence $\ucl'$ is $K$-equivariant. Since $K$ acts freely on $G'$ there is a
well defined local system $\ucl$ of rank $1$ on $G$ such that $\p^*\ucl=\ucl'$. 

\subhead 2.2\endsubhead
Define $f:G-B@>>>T$ by $f(y)=k(y)\ds\i$ and $f^1:B@>>>T$ by $f^1(y)=k(y)$. We show:

(a) {\it we have canonically $\ucl|_{G-B}=f^*\cl$ and $\ucl|_B=f^{1*}\cl$.}
\nl
In the setup of 2.1, assume first that we are in case 1. Let $j:G-B@>>>G$, $h:B@>>>G$ be
the inclusions. We must show that 
$$j^*e^*(d\i)^*\cl_1=f^*c^*\cl_1,h^*e^*(d\i)^*\cl_1=f^{1*}c^*\cl_1.$$ 
It is enough to show that $d\i ej=cf$, $d\i eh=cf^1$ or that $ej=dcf$ (resp. 
$eh=dcf^1$). Both maps take $ut\ds u'$ (resp. $ut$), where $u,u'\in U,t\in T$, to the 
image of $t$ in $G/G_{der}$. (We use that $u,u',\ds\in G_{der}$.)

Next we assume that we are in case 2. Define $\ds'\in G'$ in terms of the unique 
\'epinglage of $G'$ compatible under $\p$ with that of $G$ in the same way that $\ds\in G$
is defined in terms of the \'epinglage of $G$. Let $B'=\p\i(B)$. Define $f':G'-B'@>>>T'$, 
$f'{}^1:B'@>>>T'$ in terms of $G',B',T',\ds'$ in the same way that $f,f^1$ are defined in 
terms of $G,B,T,\ds$. Let $\p_s:G'-B'@>>>G-B$, $\p_1:B'@>>>B$ be the restrictions of $\p$.
It is enough to show that
$$\p_s^*(\ucl|_{G-B})=\p_s^*f^*\cl,\qua\p_1^*(\ucl|_B)=\p_1^*f^{1*}\cl.$$
We have 
$$\p_s^*(\ucl|_{G-B})=\ucl'|_{G'-B'}, \p_s^*f^*\cl=f'{}^*\cl', \p_1^*(\ucl|_B)=\ucl'|_{B'},
\p_1^*f^{1*}\cl=f'{}^{1*}\cl'.$$
Hence it is enough to show that $\ucl'|_{G'-B'}=f'{}^*\cl'$, $\ucl'|_{B'}=f'{}^{1*}\cl'$.
But these are known from Case 1 applied to $G',\cl'$ instead of $G,\cl$. This proves (a).

\subhead 2.3\endsubhead
We return to the general case. Let $s\in\II$. Let $P=P_s$. Let $\p_P:P@>>>P/U_P$ be the 
obvious map. Note that $P/U_P$ inherits an \'epinglage from $G$ and that $T$, identified 
with its image under $\p_P$ is a maximal torus of $P/U_P$. To any $\cl\in\cs(T)$ such that
$\cha_s^*\cl\cong\bbq$ we associate a local system of rank $1$ on $P$, namely the inverse 
image of the local system $\ucl$ on $P/U_P$ (see 2.1) under $\p_P$; this local system on 
$P$ is denoted again by $\ucl$.

Define $f_s:P-B@>>>T$ by $f_s(y)=k(y)\ds\i$ and $f_s^1:B@>>>T$ by $f_s^1(y)=k(y)$. From 
2.2(a) we deduce by taking inverse image under $\p_P$:

(a) {\it we have canonically $\ucl|_{P-B}=f_s^*\cl$ and $\ucl|_B=f_s^{1*}\cl$,}
\nl
(as local systems over subsets of $P$.)

\subhead 2.4\endsubhead
Let $\ss=(s_1,s_2,\do,s_r)$ be a sequence in $\II$. Let $\cl\in\cs(T)$. Let 
$$\ciss=\{i\in[1,r];s_1s_2\do s_i\do s_2s_1\in W_\cl\}.$$
Let 
$$\cy=\{(y_i)\in G^{[1,r]};y_i\in P_{s_i}(i\in\ciss),y_i\in P_{s_i}-B(i\in[1,r]-\ciss)\}.$$
For $i\in[1,r]$ we define $f_{s_i}:P_{s_i}-B@>>>T$ by $f_{s_i}(y)=k(y)\ds_i\i$ and
$f_{s_i}^1:B@>>>T$ by $f_{s_i}^1(y)=k(y)$. We have obvious projections 
$p_i:\cy@>>>P_{s_i}(i\in\ciss)$, $p_i:\cy@>>>P_{s_i}-B(i\in[1,r]-\ciss)$. Let
$$\uucl=\ot_{i\in[1,r]}\cf_i$$
with
$$\cf_i=p_i^*\un{s_{i-1}^*\do s_2^*s_1^*\cl}\text{ for }i\in\ciss,
\cf_i=p_i^*f_{s_i}^*s_{i-1}^*\do s_2^*s_1^*\cl\text{ for }i\in[1,r]-\ciss.$$
Here $\cf_i,\uucl$ are local systems on $\cy$. Note that if $i\in\ciss$ then the local
system $\un{s_{i-1}^*\do s_2^*s_1^*\cl}$ on $P_{s_i}$ is well defined (see 2.1) since 
$\cha_{s_i}^*(s_{i-1}^*\do s_2^*s_1^*\cl)\cong\bbq$.

For any $\cj\sub\ciss$, let 
$$\cy^\cj=\{(y_i)\in G^{[1,r]};y_i\in P_{s_i}-B(i\in[1,r]-\cj),y_i\in B(i\in\cj)\}$$
and let $\sscj=(s'_1,s'_2,\do,s'_r)$ where $s'_i=s_i$ if $i\in[1,r]-\cj$, $s'_i=1$ if 
$i\in\cj$. Define $f^\cj:\cy^\cj@>>>T$ by 
$(y_i)\m k(y_1)k(y_2)\do k(y_r)[\sscj]^\bul{}\i$. We show:

(a) {\it We have canonically $\uucl|_{\cy^\cj}=(f^\cj)^*\cl$.}
\nl
We have obvious projections $p'_i:\cy^\cj@>>>P_{s_i}-B(i\in[1,r]-\cj)$, 
$p'_i:\cy^\cj@>>>B(i\in\cj)$. Using 2.3(a) we have canonically 
$\uucl|_{\cy^\cj}=\ot_{i\in[1,r]}\cf'_i$ where
$$\cf'_i=p'_i{}^*f_{s_i}^*s_{i-1}^*\do s_2^*s_1^*\cl\text{ for }i\in[1,r]-\cj,
\cf'_i=p'_i{}^*f_{s_i}^{1*}s_{i-1}^*\do s_2^*s_1^*\cl\text{ for }i\in\cj.$$
We define $\tf:\cy^\cj@>>>T^{[1,r]}$ by $\tf=(\tf_i)$ where for $i\in[1,r]$,
$\tf_i:\cy^\cj@>>>T$ is given by
$$\tf_i=s_1s_2\do s_{i-1}f_{s_i}p'_i\text{ for }i\in[1,r]-\cj,
\tf_i=s_1s_2\do s_{i-1}f_{s_i}^1p'_i\text{ for }i\in\cj.$$
Then $\cf'_i=\tf_i^*\cl$ for $i\in[1,r]$ and 
$$\ot_{i\in[1,r]}\cf'_i=\tf^*(\cl\bxt\cl\bxt\do\bxt\cl).$$
For $y\in\cy^\cj)$ we have $f^\cj(y)=\tf_1(y)\tf_2(y)\do\tf_r(y)=m\tf(y)$ where 
$m:T^{[1,r]}@>>>T$ is multiplication. Hence $(f^\cj)^*\cl=\tf^*m^*\cl$. It is then enough 
to show that $m^*\cl=\cl\bxt\cl\bxt\do\bxt\cl$; this is a known property of any local
system in $\cs(T)$. This proves (a).

Let 
$$Y=\{(y_i)\in G^{[1,r]};y_i\in P_{s_i}(i\in[1,r])\}.$$
Note that $\cy$ is an open dense subset of $Y$. Hence $IC(Y,\uucl)$ is well defined. We 
show:
$$IC(Y,\uucl)|_{Y-\cy}=0.\tag b$$
For any $j\in[1,r]-\ciss$ let
$$\D_j=\{\{(y_i)\in G^{[1,r]};y_i\in P_{s_i}(i\in[1,r]-\{j\}),y_j\in B\}.$$
Clearly, $\{\D_j,j\in[1,r]-\ciss\}$ are smooth divisors with normal crossings in the smooth
variety $Y$. Using \cite{\CS, I, 1.6} we see that it suffices to prove the following 
statement.

(c) {\it For $j\in[1,r]-\ciss$, the monodromy of $\uucl$ around the divisor $\D_j$ is
non-trivial.}
\nl
We define a cross-section $\xi:\kk@>>>Y$ to $\D_j$ in $Y$ by
$$\xi(a)=(\ds_1,\do,\ds_{j-1},y_{s_j}(-a),\ds_{j+1},\do,\ds_r).$$
We have $\xi(0)\in\D_j$, 
$\xi(a)\in\cy$ for $a\in\kk^*$. Let $\xi':\kk^*@>>>\cy$ be the restriction of $\xi$. It is
enough to show that $\xi'{}^*\uucl\not\cong\bbq$ or, with notation in 2.4, that 
$$\xi'{}^*p_j^*f_{s_j}^*s_{j-1}^*\do s_2^*s_1^*\cl\not\cong\bbq$$
or that
$$(s_1s_2\do s_{j-1}f_{s_j}p_j\xi')^*\cl\not\cong\bbq,
(s_1s_2\do s_{j-1}\cha_{s_j})^*\cl\not\cong\bbq.$$
(We have 
$f_{s_j}p_j\xi'(a)=k(y_{s_j}(-a))\ds_j\i=\cha_{s_j}(a)$.) This follows from the fact that 
$j\notin\ciss$. This proves (c) and hence (b).

(A similar result with a similar proof appears in \cite{\CSD, VI, 28.10(b)}.)

\mpb

{\it In the remainder of this paper we assume that $\kk$ is an algebraic closure of 
$\FF_q$, a finite field with $q$ elements and that we are given a fixed $\FF_q$-rational 
structure on $G$ such that $B$ and $T$ are defined over $\FF_q$.} Let $F:G@>>>G$ be the 
corresponding Frobenius map. We set 
$$\G=\{g\in G;F(g)=g\}.$$
Now $F:G@>>>G$ induces an isomorphism $F:T@>>>T$. For $w\in W$ we write $\cl\in\cs(T)^{wF}$
instead of "$\cl\in\cs(T)$ and $(wF)^*\cl\cong\cl$".

Define $F_0:T@>>>T$ by $t\m t^q$. For any $t\in T$ we have $F(t)=F_0(\boc(t))=\boc(F_0(t))$
where $\boc:T@>>>T$ is a well defined automorphism. Then $w\m\boc w\boc\i$ is an 
automorphism of $W$ denoted also by $w\m\boc(w)$. This restricts to a bijection 
$\II@>\si>>\II$. Let
$$\uB:=\{(b_i)\in B^{[0,r]};k(b_r\i F(b_0))=1\}.$$
We show:

(d) {\it If $\cl\in\cs(T)^{[\ss]F}$ then the local system $\uucl$ is equivariant for the 
$\uB$-action 
$$(b_0,b_1,\do,b_r):(y_1,y_2,\do y_r)\m(b_0y_1b_1\i,b_1y_2b_2\i,\do,b_{r-1}y_rb_r\i)$$
on $\cy$.}
\nl
By (a), the restriction of $\uucl$ to the $\uB$-stable open dense subset $\cy^\em$ of $\cy$
is $(f^\em)^*\cl$. Since $\cy$ is smooth, it is enough to show that the local system 
$(f^\em)^*\cl$ on $\cy^\em$ is $\uB$-equivariant. Now $\uB$ acts on $T$ by
$(b_0,b_1,\do,b_r):t\m k(b_0)\i t([\ss]F(b_0))$ and $f^\em:\cy^\em@>>>T$ is compatible with
the $\uB$-actions. Hence it is enough to show that $\cl$ is $\uB$-equivariant. An 
equivalent statement is that $\cl$ is $T$-equivariant for the $T$-action 
$t_0:t\m t_0\i t([\ss]F(t_0))$. This follows from our assumption on $\cl$; (d) is proved.

\subhead 2.5\endsubhead
Let $\ww=(w_1,w_2,\do,w_r)$ be a sequence in $W$. Let
$$\align&Z^\ww=\{(B_i)\in\cb^{[0,r]};\po(B_{i-1},B_i)=w_i(i\in[1,r]),B_r=F(B_0)\},\\&
\dZ=\{(g_iU)\in(G/U)^{[0,r]};k(g_{i-1}\i g_i)=\dw_i(i\in[1,r]),g_r\i F(g_0)\in U\}.
\endalign$$ 
Let 
$$\fT=\{(t_i)\in T^{[0,r]};t_i=w_i\i(t_{i-1}) (i\in[1,r]),t_r=F(t_0)\},$$
a finite subgroup of $T^{[0,r]}$ which may be identified via $(t_i)\m t_0$ with $T^{F'}$
where $F':T@>>>T$ is $t\m[\ww]F(t)$. The free $\fT$-action $(t_i):(g_iU)\m(g_it_i\i U)$ 
on $\dZ$ makes $\dZ$ into a principal $\fT$-bundle over $Z^\ww$ via the map 
$f:\dZ@>>>Z^\ww$, $(g_iU)\m(g_iBg_i\i)$. Now $f_!\bbq$ is a local system on $Z^\ww$ with a
free action of $\fT=T^{F'}$ on each stalk. We have 
$f_!\bbq=\op_{\c\in\Hom(T^{F'},\bbq^*)}f_!^\c\bbq$ where $f_!^\c\bbq$ is the subsheaf of 
$f_!\bbq$ on which $T^{F'}$ acts according to $\c$. 

Now $\G$ acts on $Z^\ww$ by $g:(B_i)\m(gB_ig\i)$, and on $\dZ$ by $g:(g_iU)\m(gg_iU)$. This
last action commutes with the $\fT$-action. Hence $f_!\bbq$ has a natural $\G$-equivariant
structure and each $f_!^\c\bbq$ inherits a $\G$-equivariant structure from $f_!\bbq$.

We now give an alternative construction of the local systems $f_!^\c\bbq$. Let
$$\uZ=\{(g_iU)\in(G/U)^{[0,r]};g_{i-1}\i g_i\in B\dw_iB(i\in[1,r]),g_r\i F(g_0)\in U\}.$$
Define $\g:\uZ@>>>Z^\ww$ by $(g_iU)\m(g_iBg_i\i)$. Define $\p_\ww:\uZ@>>>T$ by 
$$(g_iU)\m k(g_0\i g_1)k(g_1\i g_2)\do k(g_{r-1}\i g_r)[\ww]^\bul{}\i.$$
The torus $\uT:=\{(t_i)\in T^{[0,r]};t_r=F(t_0)\}$ acts on $\uZ$ by 
$(t_i):(g_iU)\m(g_it_i\i U)$ and on $T$ by $(t_i):t\m t_0t([\ww]F(t_0\i))$. These actions 
are compatible with $\p_\ww$. Let $\cl\in\cs(T)^{[\ww]F}$. By
1.8, $\cl$ is equivariant for the $\uT$-action on $T$. Hence $\p_\ww^*\cl$ is equivariant 
for the (free) $\uT$-action on $\uZ$. Hence $\p_\ww^*\cl=\g^*\cl_\ww$ for a well defined 
local system $\cl_\ww$ on $Z^\ww$. 

Now $\G$ acts on $\uZ$ by $g:(g_iU)\m(gg_iU)$ and on $T$ trivially. Also, $\cl$ has a 
natural $\G$-equivariant structure in which $\G$ acts trivially on each stalk of $\cl$. 
Since $\p_\ww$ is compatible with the $\G$-actions it follows that $\p_\ww^*\cl$ has a 
natural $\G$-equivariant structure. Since $\g$ is compatible with the $\G$-actions it 
follows that $\cl_\ww$ has a natural $\G$-equivariant structure. 

Now assume that $\cl$ and $\c\in\Hom(T^{F'},\bbq^*)$ correspond to each other as in 1.9. 
Thus we assume that $\cl=L^\c_!\bbq$ where $L:T@>>>T$ is as in 1.9. We 
show 
$$\cl_\ww=f_!^\c\bbq.\tag a$$
Since $\g$ is smooth with connected fibres it is enough to show that 
$\p_\ww^*\cl=\g^*f_!^\c\bbq$. Let 
$$\align
\fP=&\{(g_iU,\t_i)\in(G/U\T T)^{[0,r]};k(\t_{i-1}g_{i-1}\i g_i\t_i\i)=\dw_i(i\in[1,r]),\\&
\t_r=F(\t_0),g_r\i F(g_0)\in U\}.\endalign$$
Define $f':\fP@>>>\uZ$ by $(g_iU,\t_i)\m(g_iU)$ and $\g':\fP@>>>\dZ$ by 
$(g_iU,\t_i)\m(g_i\t_i\i U)$. Define $\p':\fP@>>>T$ by $(g_iU,\t_i)\m\t_0$. Now $\fT$ acts
on $\fP$ by $(t_i):(g_iU,\t_i)\m(g_iU,t_i\t_i)$, making $f'$ into a principal $\fT$-bundle.
We have a cartesian diagram of principal $\fT=T^{F'}$ bundles:
$$\CD
\dZ@<\g'<<\fP@>\p'>>T\\
@VfVV   @Vf'VV   @VLVV\\
Z^\ww@<\g<<\uZ@>\p_\ww>>T
\endCD$$ 
It follows that $\g^*(f_!\bbq)=f'_!\bbq=\p_\ww^*(L_!\bbq)$, and taking $\c$-eigenspaces:
$\g^*(f_!^\c\bbq)=f'_!{}^\c\bbq=\p_\ww^*(L_!^\c\bbq)$. Thus, 
$\g^*(f_!^\c\bbq)=\p^*_\ww\cl$, as required.

From the definitions we see that (a) is compatible with the $\G$-equivariant structures.

When $\ww=(w)$ is a one term sequence with $w\in W$ we can identify $Z^\ww$ with
$$\cb_w=\{B'\in\cb;\po(B',F(B'))=w\}$$ 
via $B'\lra(B',F(B')$. Note that $\cb_w$ is stable under conjugation by $\G$.
For $\cl\in\cs(T)^{wF}$, the local 
system $\cl_{(w)}$ on $Z^\ww$ can be then identified with a local system $\cl_w$ on 
$\cb_w$. The subvarieties $\cb_w$ of $\cb$ and the local systems $\cl_w$ were introduced 
in \cite{\DL}. 

\subhead 2.6\endsubhead
In the remainder of this section we fix a sequence $\ss=(s_1,s_2,\do,s_r)$ in $\II$ and
$\cl\in\cs(T)^{[\ss]F}$. Let $\ciss$ be as in 2.4. Let 
$$\bZ^\ss=\{(B_i)\in\cb^{[0,r]};\po(B_{i-1},B_i)\in\{1,s_i\}(i\in[1,r]),B_r=F(B_0)\},$$
$$\align&\cz^\ss=\{(B_i)\in\cb^{[0,r]};\po(B_{i-1},B_i)\in\{1,s_i\}(i\in\ciss)\\&
\po(B_{i-1},B_i)=s_i(i\in[1,r]-\ciss),B_r=F(B_0)\}.\endalign$$
The variety $\bZ^\ss$ was introduced in \cite{\DL} (in the case where $l([\ss])=r$).
For $\cj\sub\ciss$ the variety $Z^{\sscj}$ (as in 2.5) can be also described as 
$$\align
Z^{\sscj}&=\{(B_i)\in\cb^{[0,r]};B_{i-1}=B_i(i\in\cj), \po(B_{i-1},B_i)=s_i(i\in[1,r]-\cj),\\&
B_r=F(B_0)\}\sub\cz^\ss.\endalign$$
Consider the commutative diagram
$$\CD     
\bZ@<\d_0<<\bZ_0@<\d_1<<\bZ_1@>\d_2>>\bZ_2\\
@A\e AA @A\e_0AA @A\e_1AA @A\e_2AA\\
Z@<d_0<<Z_0@<d_1<<Z_1@>d_2>>Z_2\\
@AeAA @Ae_0AA @Ae_1AA @Ae_2AA\\
Z^\cj@<d'_0<<Z_0^\cj@<d'_1<<Z_1^\cj@>d'_2>>Z_2^\cj\\
@.@Vf_0^\cj VV @Vf_1^\cj VV @Vf_2^\cj VV \\  
{}@.T@<=<<T@>=>>T\endCD$$
where the following notation is used.

$\bZ=\bZ^\ss,Z=\cz^\ss,Z^\cj=Z^{\sscj}$.

$\bZ_0$ is the set of all $(g_0U,g_1U,\do,g_rU)\in(G/U)^{[0,r]}$ such that 
$g_{i-1}\i g_i\in P_{s_i}$ for $i\in[1,r]$ and $g_r\i F(g_0)\in U$.

$Z_0$ is the set of all $(g_0U,g_1U,\do,g_rU)\in(G/U)^{[0,r]}$ such that 
$g_{i-1}\i g_i\in P_{s_i}$ for $i\in\ciss$, $g_{i-1}\i g_i\in P_{s_i}-B$ for 
$i\in[1,r]-\ciss$ and $g_r\i F(g_0)\in U$.

$Z_0^\cj$ is the set of all $(g_0U,g_1U,\do,g_rU)\in(G/U)^{[0,r]}$ such that 
$g_{i-1}\i g_i\in P_{s_i}-B$ for $i\in[1,r]-\cj$, $g_{i-1}\i g_i\in B$ for $i\in\cj$ and
$g_r\i F(g_0)\in U$.

$d_0,\d_0$ are given by $(g_0U,g_1U,\do,g_rU)\m(g_0Bg_0\i,g_1Bg_1\i,\do,g_rBg_r\i)$.

$f^\cj_0$ is given by
$(g_0U,g_1U,\do,g_rU)\m k(g_0\i g_1)\do k(g_{r-1}\i g_r)[\sscj]^\bul{}\i$.

$\bZ_1$ is the set of all $(y_0,y_1,\do,y_r)\in G^{[0,r]}$ such that 
$y_i\in P_{s_i}(i\in[1,r])$, $y_0\i F(y_0)\in y_1y_2\do y_rU$.

$Z_1$ is the set of all $(y_0,y_1,\do,y_r)\in G^{[0,r]}$ such that 
$y_i\in P_{s_i}(i\in\ciss)$, $y_i\in P_{s_i}-B(i\in[1,r]-\ciss)$, 
$y_0\i F(y_0)\in y_1y_2\do y_rU$.

$Z_1^\cj$ is the set of all $(y_0,y_1,\do,y_r)\in G^{[0,r]}$ such that 
$y_i\in P_{s_i}-B(i\in[1,r]-\cj)$, $y_i\in B(i\in\cj)$, $y_0\i F(y_0)\in y_1y_2\do y_rU$.

$d_1,\d_1$ are given by $(y_0,y_1,\do,y_r)\m(y_0U,y_0y_1U,\do,y_0y_1\do y_rU)$.

$f^\cj_1$ is $(y_0,y_1,\do,y_r)\m k(y_1)k(y_2)\do k(y_r)[\sscj]^\bul{}\i$. 

$\bZ_2=Y, Z_2=\cy, Z_2^\cj=\cy^\cj$. (See 2.4.)

$d_2,\d_2$ are given by $(y_0,y_1,\do,y_r)\m(y_1,\do,y_r)$.

$f^\cj_2$ is $(y_1,\do,y_r)\m k(y_1)k(y_2)\do k(y_r)[\sscj]^\bul{}\i$. 
\nl
The maps $e,e_i,\e_i(i\in[0,2])$ are the obvious imbeddings. For $i\in[0,2]$ the map $d'_i$
is the restriction of $d_i$. From the definitions we have:

(a) {\it In our commutative diagram, all squares that do not involve $T$ are cartesian.}

\subhead 2.7\endsubhead
$\G$ acts:

 on $Z$ and $\bZ$ by $g:(B_0,B_1,\do,B_r)\m(gB_0g\i,gB_1g\i,\do,gB_rg\i)$;

on $Z_0$ and $\bZ_0$ by $g:(g_0U,g_1U,\do,g_rU)\m(gg_0U,gg_1U,\do,gg_rU)$;

on $Z_1$ and $\bZ_1$ by $g:(y_0,y_1,\do,y_r)\m(gy_0,y_1,y_2,\do,y_r)$;

on $Z_2$, $\bZ_2$ trivially.

The subsets $Z^\cj,Z^\cj_i$ of $Z,Z_i(i\in[0,2])$ are stable under the $\G$-action. The 
maps $d_i,d'_i,\d_i,f^\cj_i$ are compatible with the $\G$-actions.

$\uB$ (see 2.4) acts: 

on $Z$ trivially;

on $Z_0$ by $(b_0,b_1,\do,b_r):(g_0U,g_1U,\do,g_rU)\m(g_0b_0\i U,g_1b_1\i U,\do,g_rb_r\i U)$;

on $Z_1$ by 

$(b_0,b_1,\do,b_r):(y_0,y_1,\do,y_r)\m
(y_0b_0\i,b_0y_1b_1\i,b_1y_2b_2\i,\do,b_{r-1}y_rb_r\i)$;

on $Z_2$ by 
$(b_0,b_1,\do,b_r):(y_1,\do,y_r)\m (b_0y_1b_1\i,b_1y_2b_2\i,\do,b_{r-1}y_rb_r\i)$.

The maps $d_i$ are compatible with the $\uB$-actions.

\subhead 2.8\endsubhead
Now 

(a) $d_0,\d_0,d'_0$ are principal $\uB/U^{[0,r]}$-bundles.

(b) $d_1,\d_1,d'_1$ are principal $U^{[0,r]}$-bundles. (The action of $U^{[0,r]}$ on
$\bZ_1,Z_1,Z^\cj_1$ is by restriction of the $\uB$-action.) 

(c) Each of $d_2,\d_2,d'_2$ is a composition of a principal $\G$-bundle with a principal 
$U$-bundle.

\subhead 2.9\endsubhead
We show for $i\in[0,2]$ that 

(a) {\it $Z_i$ is smooth of pure dimension say $\fd_i$ and it is open dense in $\bZ_i$.}
\nl
Let $P_i$ be the property expressed by (a). It is obvious that $P_2$ holds. Using $P_2$, 
2.6(a) and 2.8(c) we see that $P_1$ holds. Using $P_1$, 2.6(a) and 2.8(b) we see that $P_0$
holds. Thus (a) holds. Using $P_0$, 2.6(a) and 2.8(a) we see that 

(b) {\it $Z$ is smooth of pure dimension say $\fd$ and it is open dense in $\bZ$.}
\nl
We show:

(c) $\fd=r$.
\nl
From the definitions we see that $\fd_2=r(\dim B+1)$. From the arguments above we see
successively that $\fd_1=r(\dim B+1)+\dim U$, $\fd_0=\fd_1-(r+1)\dim U$
$\fd=\fd_0-r\dim T$; (c) follows.

We show:

(d) {\it The natural $\G$-action on the set of connected components of $Z_i$ ($i\in[0,2]$)
or of $Z$ is transitive.}
\nl
For $Z_2$ this is clear since $Z_2$ is connected. This also implies the result for $Z_1$
(see 2.8(c)). Using 2.8(b),(a) we deduce that the result also holds for $Z_0$ and for $Z$.

\subhead 2.10\endsubhead
Let $\cj\sub\ciss$. For $i\in[0,2]$ we set $\cl^\cj_i=(f^\cj_i)^*\cl$, a local system of
rank $1$ on $Z_i^\cj$. Since $\cl$ has a natural $\G$-equivariant structure (with $\G$ 
acting trivially on each stalk) and $f^\cj_i$ is compatible with the $\G$-actions we see 
that $\cl_i^\cj$ has a natural $\G$-equivariant structure. From the definitions we have
isomorphisms compatible with the $\G$-equivariant structures as follows:
$$d'_1{}^*\cl_1^\cj\cong\cl_0{}^\cj; d'_2{}^*\cl_2^\cj\cong\cl_1^\cj.\tag a$$
From the definitions we see that

(b) {\it $\G$ acts trivially on any stalk of $\cl_2^\cj$.}
\nl
Let $\cl^\cj$ be the local system on $Z^\cj=Z^{\sscj}$ denoted in 2.5 by $\cl_\ww$ 
where $\ww=\sscj$. This is well defined since for $\cj\sub\ciss$ we have
$\cl\in\cs(T)^{[\sscj]F}$. (We use that $\cl\in\cs(T)^{[\ss]F}$ and 
$(s_1s_2\do s_j\do s_2s_1)^*\cl\cong\cl$ for any $j\in\cj$.) As in 2.5, $\cl^\cj$ has a 
natural $\G$-equivariant structure. From the definitions we have 
$d'_0{}^*\cl^\cj\cong\cl_0^\cj$ compatibly with the $\G$-equivariant structures.

\subhead 2.11\endsubhead
For $i\in[1,2]$ we define a local system $\bcl_i$ on $Z_i$ by $\bcl_2=\uucl$, 
$\bcl_1=d_2^*\bcl_2$ where $\uucl$ is as in 2.4. From 2.4(b) and the results in 2.7 we see
that $\bcl_i$ is $\uB$-equivariant. Since $d_0d_1:Z_1@>>>Z$ is a principal $\uB$-bundle we
see that there is a well defined local system $\bcl$ on $Z$ such that 
$(d_0d_1)^*\bcl=\bcl_1$. Let $\bcl_0=d_0^*\bcl$. Then $\bcl_1=d_1^*\bcl_0$. We regard 
$\bcl_2$ as a $\G$-equivariant local system on $Z_2$ with $\G$ acting trivially on each 
stalk. Since each $d_i$ is compatible with the $\G$-actions we see that $\bcl_i$ 
($i\in[0,2]$) and $\bcl$ have natural $\G$-equivariant structures which are compatible with
$d_i^*$. 

\subhead 2.12\endsubhead
We show:

(a) {\it For any $\cj\sub\ciss$ we have $\bcl|_{Z^\cj}\cong\cl^\cj$ compatibly with the
$\G$-equivariant structures.}

(b) {\it For any $i\in[0,2]$ and $\cj\sub\ciss$ we have 
$\bcl_i|_{Z^\cj_i}\cong\cl_i^\cj$ compatibly with the $\G$-equivariant structures.}
\nl
Note that (b) holds for $i=2$ by 2.4(a) (the compatibility with the $\G$-equivariant 
structures is automatic since $\G$ acts trivially on each stalk of the local systems 
involved). From this we get (using 2.10, 2.11) that (b) holds for $i=1$, then
for $i=0$, and then that (a) holds.

\subhead 2.13\endsubhead
We show:

(a) {\it We have $IC(\bZ,\bcl)|_{\bZ-Z}=0$.}

(b) {\it For $i\in[0,2]$ we have $IC(\bZ_i,\bcl_i)|_{\bZ_i-Z_i}=0$.}
\nl
Note that the $IC$ complexes in (a),(b) are well defined by 2.9(a),(b). Now (b) holds for 
$i=2$ by 2.4(b). From this we get (using 2.8, 2.11) that (b) holds for $i=1$, 
then for $i=0$, and then that (a) holds.

\subhead 2.14\endsubhead
Assume that $r\ge2,h\in[2,r]\cap\ciss,s_{h-1}=s_h$. We set \lb
$\ss':=(s_1,\do,s_{h-1},s_{h+1},\do,s_r)$. Then $\cl\in\cs(T)^{[\ss']F}$ so that 
$Z':=\cz^{\ss'}$ is defined as in 2.6. We have a commutative diagram
$$\CD
Z@<d_0<<Z_0@<d_1<<Z_1@>d_2>>Z_2\\
@V\b VV     @V\b_0VV      @V\b_1VV      @V\b_2VV\\
Z'@<d'_0<<Z'_0@<d'_1<<Z'_1@>d'_2>>Z'_2                    \endCD$$
where 

the upper row is as in 2.6, 

the lower row is defined analogously in terms of $\ss',\cl$ instead of $\ss,\cl$,

$\b$ is $(B_0,B_1,\do,B_r)\m(B_0,B_1,\do,B_{h-2},B_h,\do,B_r)$,

$\b_0$ is $(g_0U,g_1U,\do,g_rU)\m(g_0U,g_1U,\do,g_{h-2}U,g_hU,g_rU)$,

$\b_1$ is $(y_0,y_1,\do,y_r)\m(y_0,y_1,\do,y_{h-2},y_{h-1}y_h,y_{h+1},\do,y_r)$,

$\b_2$ is $(y_1,\do,y_r)\m(y_1,\do,y_{h-2},y_{h-1}y_h,y_{h+1},\do,y_r)$.
\nl
Let $\bcl,\bcl_i(i\in[0,2])$ be the local systems on $Z,Z_i$ defined in 2.11; let 
$\bcl',\bcl'_i$ be the analogous local systems on $Z',Z'_i$. We show:
$$\bcl\cong\b^*(\bcl').\tag a$$
It is enough to show that $\d_1^*\d_0^*\bcl\cong\d_1^*\d_0^*\b^*(\bcl')$ or equivalently 
that $\bcl_1\cong\b_1^*\bcl'_1$. Hence it is enough to show that 
$\d_2^*\bcl_2\cong\b_1^*\d'_2{}^*\bcl'_2$ or equivalently that
$\d_2^*\bcl_2\cong\d_2^*\b_2{}^*\bcl'_2$. It is enough to show that 
$\bcl_2\cong\b_2{}^*\bcl'_2$. From the definition of $\uucl$ in 2.4 and with the notation 
in 2.4 we see that it is enough to show that $m^*{}'\ucl\cong{}'\ucl\bxt{}'\ucl$ (local systems
on $P_{s_h}$) where ${}'\cl=s_{h-1}^*\do s_2^*s_1^*\cl\cong s_h^*s_{h-1}^*\do s_2^*s_1^*\cl$ 
and $m:P_{s_h}\T P_{s_h}@>>>P_{s_h}$ is multiplication. It is enough to show that, in the
setup of 2.1 we have $m'{}^*\ucl\cong\ucl\bxt\ucl$ (local systems on $G$) where
$m':G\T G@>>>G$ is multiplication. This follows from the definitions in 2.1 using the
isomorphism $m_1^*\cl\cong\cl\bxt\cl$ (local systems on $T$) where $m_1:T\T T@>>>T$ is 
multiplication. 

\subhead 2.15\endsubhead
Assume that $r\ge2,h\in[2,r],h\n\ciss,s_{h-1}=s_h$. Let $Z^1$ be the open subset of $Z$ 
defined by the condition $\po(B_{h-2},B_h)=s_i$. Define $\b:Z^1@>>>\cb^r$ by
$(B_0,B_1,\do,B_r)\m(B_0,B_1,\do,B_{h-2},B_h,\do,B_r)$. We show:
$$\b_!(\bcl|_{Z^1})=0.\tag a$$
Let $p=(B_0,B_1,\do,B_{h-2},B_h,\do,B_r)\in\cb^r$ be such that $\Ph:=\b\i(p)\ne\em$. Then 
$$\Ph
=\{(B_0,B_1,\do,B_{h-2},\tB,B_h,\do,B_r);\tB\in\cb,\po(B_{h-2},\tB)=\po(\tB,B_h)=s_h\}.$$
It is enough to show that $H^*_c(\Ph,\bcl)=0$. Let 
$$\align\Ph'=&\{(B_0,B_1,\do,B_{h-2},\tB,B_h,\do,B_r);
\tB\in\cb,\po(B_{h-2},\tB)=s_h,\\&\po(\tB,B_h)\in\{1,s_h\}\}.\endalign$$
Then $\Ph'$ is an affine line which is a cross section in $\bZ$ to the divisor $\D_h$ (see
2.4) and $\Ph'\cap\D_h$ is the point $p'=(B_0,B_1,\do,B_{h-2},B_h,B_h,\do,B_r)$. Moreover,
$\Ph=\Ph'-\{p'\}$. The vanishing $H^*_c(\Ph,\bcl)=0$ follows from 2.4(c).

Now let $Z^2=Z-Z^1$, that is the closed subset of $Z$ defined by the condition
$B_{h-2}=B_h$. Let $\ss':=(s_1,\do,s_{h-2},s_{h+1},\do,s_r)$. We have 
$\cl\in\cs(T)^{[\ss']F}$ so that $Z':=\cz^{\ss'}$ is defined as in 2.6. Define 
$\b':Z^2@>>>Z'$ by $(B_0,B_1,\do,B_r)\m(B_0,B_1,\do,B_{h-2},B_{h+1},\do,B_r)$, an affine 
line bundle. Let $\bcl$ be the local systems on $Z$ defined in 2.11; let 
$\bcl'$ be the analogous local system on $Z'$. From the definitions we have:
$$\bcl|_{Z^2}=\b'{}^*(\bcl').\tag b$$

\head 3. The class $\Bbb S'(\cp_J)$ of simple objects in $\cm_\G(\cp_J)$\endhead
\subhead 3.1\endsubhead
Let $J\sub\II$. We view $\cp_J$ as a variety with $\G$-action (conjugation). Hence 
$\cm_\G(\cp_J)$ is well defined.

\subhead 3.2\endsubhead
Let $\cl\in\cs(T)$. If $\ww$ is as in 2.5 and $\cl\in\cs(T)^{[\ww]F}$ then the local system
$\cl_\ww$ on $Z^\ww$ has a natural $\G$-equivariant structure (see 2.5). The map
$\Pi^\ww:Z^\ww@>>>\cp_J$, $(B_0,B_1,\do,B_r)\m P_{B_0,J}$, commutes with the $\G$-actions.
Hence for any $j\in\ZZ$, ${}^pH^j(\Pi^\ww_!\cl_\ww)$ is an object of $\cm_\G(\cp_J)$.

If $\ss$ is as in 2.6 and $\cl\in\cs(T)^{[\ss]F}$ then the local system $\bcl$ on $\cz^\ss$
has a natural $\G$-equivariant structure (see 2.6, 2.11). Hence $\bcl^\sh=IC(\bZ^\ss,\bcl)$
(see 2.6, 2.11, 2.13) has a natural $\G$-equivariant structure. Define 
$\Upss:\cz^\ss@>>>\cp_J$ and $\bUp^\ss:\bZ^\ss@>>>\cp_J$ by 
$(B_0,B_1,\do,B_r)\m P_{B_0,J}$. These maps commute with the $\G$-actions. Hence for any 
$j\in\ZZ$,
$${}^pH^j(\Upss_!\bcl)={}^pH^j(\bUp^\ss_!\bcl^\sh)\tag a$$
is an object of $\cm_\G(\cp_J)$. (The equality in (a) follows from by 2.13.)

In 3.3-3.7 we will show that the following conditions for a simple object $K$ in 
$\cm_\G(\cp_J)$ are equivalent:

(i) $K\dsv_\G{}^pH^{\cdot}(\Pi^\ww_!\cl_\ww)$ for some one term sequence 
$\ww$ in $W$ such that $\cl\in\cs(T)^{[\ww]F}$.

(ii) $K\dsv_\G{}^pH^{\cdot}(\Pi^\ww_!\cl_\ww)$ for some sequence $\ww$
in $W$ such that $\cl\in\cs(T)^{[\ww]F}$.

(iii) $K\dsv_\G{}^pH^{\cdot}(\Pi^\ww_!\cl_\ww)$ for some sequence $\ww$
in $\II\cup\{1\}$ such that $\cl\in\cs(T)^{[\ww]F}$.

(iv) $K\dsv_\G{}^pH^{\cdot}(\Pi^\ss_!\cl_\ss)$ for some sequence $\ss$
in $\II$ such that $\cl\in\cs(T)^{[\ss]F}$.

(v) $K\dsv_\G{}^pH^{\cdot}(\Upss_!\bcl)$ for some sequence $\ss$ in 
$\II$ such that $\cl\in\cs(T)^{[\ss]F}$.

(vi) $K\dsv_\G{}^pH^{\cdot}(\bUp^\ss_!\bcl^\sh)$ for some sequence $\ss$ in $\II$ such that
$\cl\in\cs(T)^{[\ss]F}$.

\subhead 3.3\endsubhead
Let $\ww=(w_1,\do,w_r)$ be a sequence in $W$ such that $\cl\in\cs(T)^{[\ww]F}$. Assume 
that for some $i\in[1,r]$, $w'_i,w''_i\in W$ satisfy $w_i=w'_iw''_i$ and 
$l(w_i)=l(w'_i)+l(w''_i)$. Let $\ww'=(w_1,\do,w_{i-1},w'_i,w''_i,w_{i+1},\do,w_r)$. Define
an isomorphism $Z^{\ww'}@>\si>>Z^\ww$ by 
$$(B_0,B_1,\do,B_{r+1})\m(B_0,B_1,\do,B_{i-1},B_{i+1},\do,B_{r+1}).$$
This isomorphism is compatible with the $\G$-actions, with the maps $\Pi^{\ww'},\Pi^\ww$
and with the local systems $\cl_{\ww'},\cl_\ww$. Hence for any $j$ we have
$${}^pH^j(\Pi^\ww_!\cl_\ww)={}^pH^j(\Pi^{\ww'}_!\cl_{\ww'})\tag a$$
(as objects of $\cm_\G(\cp_J)$).
Applying (a) repeatedly we see that conditions 3.2(ii), 3.2(iii), 3.2(iv) are equivalent.

\subhead 3.4\endsubhead
We prove the equivalence of conditions 3.2(iii), 3.2(v).

Let $\ss=(s_1,\do,s_r)$ be a sequence in $\II$ such that $\cl\in\cs(T)^{[\ss]F}$. Define 
a sequence ${}^0\cz\supset{}^1\cz\supset\do$ of closed subsets of $\cz^\ss$ by 
${}^i\cz=\cup_{\cj\sub\ciss;|\cj|\ge i}Z^{\sscj}$ (notation of 2.6). Let
$f^i:{}^i\cz@>>>\cz^\ss$, $f'{}^i:{}^i\cz-{}^{i+1}\cz@>>>\cz^\ss$ be the inclusions. The
natural distinguished triangle
$$(\Upss_!f'{}^i_!f'{}^{i*}\bcl,\Upss_!f^i_!f^{i*}\bcl,\Upss_!f^{i+1}_!(f^{i+1*}\bcl))$$
gives rise for any $i\ge 0$ to a long exact sequence in $\cm_\G(\cp_J)$:
$$\align&\do@>>>{}^pH^{j-1}(\Upss_!f^{i+1}_!f^{i+1*}\bcl)@>>>
\op_{\cj\sub\ciss;|\cj|=i}{}^pH^j(\Pi^{\sscj}_!\cl_{\sscj})
@>>>{}^pH^j(\Upss_!f^i_!f^{i*}\bcl)\\&@>>>{}^pH^j(\Upss_!f^{i+1}_!f^{i+1*}\bcl)@>>>
\op_{\cj\sub\ciss;|\cj|=i}{}^pH^{j+1}(\Pi^{\sscj}_!\cl_{\sscj})@>>>\do.\tag a\endalign$$
Here we have used the equality
$$\Upss_!f'{}^i_!f'{}^{i*}\bcl=\op_{\cj\sub\ciss;|\cj|=i}\Pi^{\sscj}_!\cl_{\sscj}.$$
(See 2.12(a).) Note that $\Upss_!f^0_!f^{0*}\bcl=\Upss_!\bcl$ and 
$\Upss_!f^i_!f^{i*}\bcl=0$ for $i$ large.

If $K$ does not satisfy 3.2(iii) then from (b) we see that for any $i\ge0$
we have $K\dsv_\G{}^pH^{\cdot}(\Upss_!f^i_!f^{i*}\bcl)$ if and only if 
$K\dsv_\G{}^pH^{\cdot}(\Upss_!f^{i+1}_!f^{i+1*}\bcl)$. Since 
$K\not\dsv_\G{}^pH^{\cdot}(\Upss_!f^i_!f^{i*}\bcl)$ with large $i$ it follows that 
$K\not\dsv_\G{}^pH^{\cdot}(\Upss_!f^0_!f^{0*}\bcl)$ that is 
$K\not\dsv_\G{}^pH^{\cdot}(\Upss_!\bcl)$. Thus, $K$ does not satisfy 3.2(v).

Assume now that $K$ satisfies 3.2(iii). We may assume that 
$K\dsv_\G{}^pH^{\cdot}(\Pi^\ss_!\cl_\ss)$ where $\ss$ as in 3.2(iii) has a 
minimum possible number of terms in $\II$. By the equivalence of 3.2(ii), 3.2(iii) we see 
that we may assume that all terms of $\ss$ are in $\II$ and that 
$K\not\dsv_\G{}^pH^{\cdot}(\Pi^{\sscj}_!\cl_{\sscj})$
for and any $\cj$ such that $|\cj|>0$. Then from (b) we see that for any $i>0$,
$K\dsv_\G{}^pH^{\cdot}(\Upss_!f^i_!f^{i*}\bcl)$ if and only if 
$K\dsv_\G{}^pH^{\cdot}(\Upss_!f^{i+1}_!f^{i+1*}\bcl)$. Since 
$K\not\dsv_\G{}^pH^{\cdot}(\Upss_!f^i_!f^{i*}\bcl)$ with large $i$ it follows that 
$K\not\dsv_\G{}^pH^{\cdot}(\Upss_!f^1_!f^{1*}\bcl)$. Using again (b) (with $i=0$) we see 
that $K\dsv_\G{}^pH^{\cdot}(\Upss_!f^0_!f^{0*}\bcl)$ hence 
$K\dsv_\G{}^pH^{\cdot}(\Upss_!\bcl)$. Thus, $K$ satisfies 3.2(v). The equivalence of 
3.2(iii), 3.2(v) is proved.

\subhead 3.5\endsubhead
Let $\ss=(s_1,s_2,\do,s_r)$ be a sequence in $\II$ such that $\cl\in\cs(T)^{[\ss]F}$.
Assume that $r\ge2,h\in[2,r]\cap\ciss$, $s_{h-1}=s_h$. Let $\ss',\bcl',\b$ be as in 2.14. 
We have $\Upss=\Up^{\ss'}\b$ and using 2.14(a) we have
$\Upss_!(\bcl)=\Up^{\ss'}_!\b_!\b^*\bcl'$. Since $\b$ is a projective line bundle we have
an exact sequence in $\cm_\G(\cp_J)$:
$$\do@>>>{}^pH^{j-2}(\Up^{\ss'}_!\bcl')(-1)@>>>{}^pH^j(\Upss_!\bcl)@>>>
{}^pH^j(\Up^{\ss'}_!\bcl')@>>>\do\tag a$$

\subhead 3.6\endsubhead
Let $\ss=(s_1,s_2,\do,s_r)$ be a sequence in $\II$ such that $\cl\in\cs(T)^{[\ss]F}$.
Assume that $r\ge2,h\in[2,r],h\n\ciss$, $s_{h-1}=s_h$. Let $Z^1,Z^2,\b,\b',\ss',\bcl'$ 
be as in 2.15; let 
$f_1:Z^1@>>>\cz^\ss,f_2:Z^2@>>>\cz^\ss$ be the inclusions. We have a distinguished triangle
$$(\Upss_!f_{1!}f_1^*\bcl,\Upss_!\bcl,\Upss_!f_{2!}f_2^*\bcl).$$
We have $\Upss_!f_{1!}=e_!\b_!$ where $e:\cb^r@>>>\cp_J$ is 
$(B_0,B_1,\do,B_{h-2},B_h,\do,B_r)\m P_{B_0,J}$. Using 2.15(a) we have 
$\Upss_!f_{1!}f_1^*\bcl=e_!\b_!(\bcl|_{Z^1})=0$. Hence the distinguished triangle above 
yields $\Upss_!\bcl=\Upss_!f_{2!}f_2^*\bcl$. We have $\Upss f_2=\Up^{\ss'}\b'$.
Using 2.15(b) we have 
$$\Upss_!f_{2!}f_2^*\bcl=\Up^{\ss'}_!\b'_!(\bcl|_{Z^2})=\Up^{\ss'}_!\b'_!\b'{}^*\bcl'=
\Up^{\ss'}_!(\bcl'\ot\b'_!\b'{}^*\bbq).$$
We see that
$$\Upss_!\bcl=\Up^{\ss'}_!\bcl'[-2](-1).\tag a$$
Hence for any $j$ we have
$${}^pH^j(\Upss_!\bcl)={}^pH^{j-2}(\Up^{\ss'}_!\bcl')(-1)\tag b$$
in $\cm_\G(\cp_J)$.

\subhead 3.7\endsubhead
Assume that 3.2(iv) holds. We show that 3.2(i) holds.

We may assume that $K\dsv_\G{}^pH^{\cdot}(\Pi^\ss_!\cl_\ss)$ for some sequence 
$\ss=(s_1,s_2,\do,s_r)$ in $\II$ such that $\cl\in\cs(T)^{[\ss]F}$ and that
$r$ is minimum possible. From the  proof in 3.4 we see that 
$K\dsv_\G{}^pH^{\cdot}(\Upss_!\bcl)$ and that $r$ is also minimal for this property.

Assume first that $l(s_1s_2\do s_r)<r$. We can find $h\in[2,r]$ such that 
$$l(s_hs_{h+1}\do s_r)=r-h+1,l(s_{h-1}s_h\do s_r)<r-h+2.$$
 We can find $s'_h,s'_{h+1},\do,s'_r$ in $\II$ such that 
$$s'_hs'_{h+1}\do s'_r=s_hs_{h+1}\do s_r=y$$
and $s'_h=s_{h-1}$. Let 
$$\uu'=(s_1,s_2,\do,s_{h-1},s'_h,s'_{h+1},\do,s'_r), \uu''=(s_1,s_2,\do,s_{h-1},y).$$
From 3.3(a) we see that
$\Pi^\ss_!\cl_\ss=\Pi^{\uu'}_!\cl_{\uu'}=\Pi^{\uu''}_!\cl_{\uu''}$. Hence we may assume
that $s_h=s_{h-1}$. 

If $h\in\ciss$ then using 3.5(a) we see that $K\dsv_\G{}^pH^{\cdot}(\Up^{\ss'}_!\bcl')$
(notation of 3.5); since $\ss'$ has $r-1$ terms this is a contradiction. If $h\n\ciss$ then
using 3.6 we see that $K\dsv_\G{}^pH^{\cdot}(\Up^{\ss'}_!\bcl')$ (notation of 3.6); since 
$\ss'$ has $r-2$ terms this is a contradiction.

We see that $l(s_1s_2\do s_r)=r$. Using 3.3(a) repeatedly we see that
${}^pH^j(\Pi^\ss_!\cl_\ss)={}^pH^j(\Pi^\ww_!\cl_\ww)$ where $\ww=(w_1)$,
$w_1=s_1s_2\do s_r$. Thus, 3.2(i) holds.

Since the implication 3.2(i)$\imp$3.2(ii) is obvious and the equivalence of 3.2(v), 3.2(vi)
follows from 3.2(a) we see that the equivalence of 3.2(i)-3.2(vi) is established.

\mpb

For an object $A$ of $\cm_\G(\cp_J)$ we write $A\in\Bbb S'(\cp_J)$ instead of "$A$ 
satisfies the equivalent conditions of 3.2(i)-3.2(vi) for some $\cl\in\cs(T)$".

\subhead 3.8\endsubhead
The results in this and the next subsection are not used in the subsequeny sections.

Let $\ss=(s_1,s_2,\do,s_r)$ be a sequence in $\II$ such that $\cl\in\cs(T)^{[\ss]F}$,
$s_1\in J$.
Let $\ss'=(s'_1,s'_2,\do,s'_r)$ where $s'_i=s_{i+1}$ for $i\in[1,r-1]$ and
$s'_r=\boc(s_1)$ where $\boc:W@>>>W$ is as in 2.4. Let $\cl'=s_1^*\cl$. We have 
$\cl'\in\cs(T)^{[\ss']F}$. Let $\bcl$ be the local system on $\cz^\ss$ defined in 2.11 and
let $\bcl'$ be the analogous local system on $\cz^\ss$ defined in terms of $\cl'$. We show:
$$\Upss_!\bcl\cong\Up^{\ss'}_!\bcl'.\tag a$$
Define $\ci'_{\ss'}$ in terms of $\ss',\cl'$ in the same way as $\ciss$ was defined in 2.4
in terms of $\ss,\cl$. If $i\in[2,r]$ we have $i\in\ciss$ if and only if 
$i-1\in\ci'_{\ss'}$. Moreover, we have $1\in\ciss$ if and only if $r\in\ci_{\ss'}$. It 
follows that $f:\cz^\ss@>>>\cz^{\ss'}$, $(B_0,B_1,\do,B_r)\m(B_1,B_2,\do,B_r,F(B_1))$, is 
well defined. From the definitions we see that $f^*\bcl'\cong\bcl$. Hence 
$\Upss_!\bcl=\Upss_!f^*\bcl'$. It remains to show that $\Up^\ss f=\Up^{\ss'}$. The first 
(resp. second) map takes $(B_0,B_1,\do,B_r)$ to $P_{B_1,J}$ (resp. $P_{B_0,J}$). It is 
enough to show that $P_{B_1,J}=P_{B_0,J}$. This follows from the fact that 
$\po(B_0,B_1)\in J$.

\subhead 3.9\endsubhead
Let $\ss=(s_1,s_2,\do,s_r)$ be a sequence in $\II$ such that $\cl\in\cs(T)^{[\ss]F}$. Let 
$s\in\II$ be such that $s\n W_\cl$. Let $\uu=(s,s_1,s_2,\do,s_r,\boc(s))$, $\cl'=s^*\cl$. 
Let $\vv=(s,s,s_1,s_2,\do,s_r)$. We have $\cl'\in\cs(T)^{[\uu]F}$, $\cl\in\cs(T)^{[\vv]F}$.
Let $\bcl$ be as in 2.11; let $\bcl',\bcl''$ be the analogous local systems on $\cz^\uu$,
$\cz^\vv$ defined in terms of $\cl',\cl$. We show
$$\Upss_!\bcl[-2](-1)=\Up^\uu_!\bcl'.\tag a$$
From 3.8 we have  $\Up^\uu_!\bcl'=\Up^\vv_!\bcl''$. From 3.6(a) we have
$\Up^\vv_!\bcl''=\Upss_!\bcl[-2](-1)$ and (a) follows. We see that
$${}^pH^j(\Up^\uu_!\bcl')={}^pH^{j-2}(\Upss_!\bcl)(-1).\tag b$$

\head 4. The class $\Bbb S(\cp_J)$ of simple objects in $\cm_\G(\cp_J)$\endhead  
\subhead 4.1\endsubhead
In this section we fix $J\sub\II$. 

In 1977 the author generalized the partition $(\cb_w)_{w\in W}$ of $\cb$ (see 2.5) by 
defining a partition of $\cp_J$ into 
finitely many pieces stable under conjugation by $\G$, as follows. To any $P\in\cp_J$ we 
associate a sequence $P^0\supset P^1\supset P^2\supset\do$ in $\cp$ by 
$$P^0=P,\qua P^n=(P^{n-1})^{F(P^{n-1})}\text{ for }n\ge1,$$
a sequence $J_0\supset J_1\supset J_2\supset\do$ of subsets of $\II$ by $P^n\in\cp_{J_n}$ 
and a sequence $w_0,w_1,w_2,\do$ in $W$ by
$$w_n=\po(P^n,F(P^n)).$$
We have 

(a) $J_0=J$,

(b) $J_n=J_{n-1}\cap w_{n-1}\boc(J_{n-1})w_{n-1}\i$ for $n\ge1$.
\nl
(see 1.3(a)),

(c) $w_n\in{}^{J_n}W^{\boc(J_n)}$ for $n\ge0$.
\nl   
Clearly, for $n\ge|\II|$ we have $P^n=P^{n+1}=\do$ hence 

(d) $w_n=w_{n+1}=\do$ and $J_n=J_{n+1}=\do$.
\nl
We set $P^\iy=P^n$ for $n\ge|\II|$, $w_\iy=w_n$ for $n\ge|\II|$, $J_\iy=J_n$ for 
$n\ge|\II|$. 

For any $\tt=(J_n,w_n)_{n\ge0}$ where $J_0\sps J_1\sps J_2\sps\do$ are subsets of $J$
satisfying (a) and $w_0,w_1,w_2,\do$ are elements of $W$ satisfying (b),(c) let $\cp_J^\tt$
be the set of all $P\in\cp_J$ which give rise to $\tt$ by the procedure above. Let 
$\ct'(J,\boc)$ be the set of all sequences $\tt$ as above such that $\cp_J^\tt\ne\em$. From
(d) we see that $\ct'(J,\boc)$ is a finite set. From (a),(b) we see that for 
$(J_n,w_n)_{n\ge0}\in\ct'(J,\boc)$, the $J_n$ are uniquely determined by the $w_n$. The 
locally closed subvarieties $\cp_J^\tt$, $\tt\in\ct'(J,\boc)$, form the desired partition 
of $\cp_J$. (See \cite{\PAR, I, 1.3, 1.4} for some examples in the classical groups.) 

\subhead 4.2\endsubhead
In this subsection we review some results in the R.B\'edard's Ph.D. Thesis (M.I.T. 1983),
see also \cite{\BE}.

(a) {\it $\ct'(J,\boc)$ is precisely the set of all $(J_n,w_n)_{n\ge0}$ with $J_n\sub\II$, 
$w_n\in W$ such that 4.1(a),(b),(c) hold and $w_n\in W_{J_n}w_{n-1}W_{\boc(J_{n-1})}$ for 
$n\ge1$.}
\nl
(With notation in \cite{\PAR, I, 2.2}, we have $\ct'(J,\boc)=\ct(\boc(J),\boc\i)$.) 

(b) {\it The assignment $(J_n,w_n)_{n\ge0}\m w_\iy$ defines a bijection
$\ct'(J,\boc)@>\si>>{}^JW$.}

(c) {\it Let $z\in{}^JW^{\boc(J)}$, $J_1=J\cap z\boc(J)z\i$. Let 
$V=\{P\in\cp_J;\po(P,F(P))=z$, $V'=\{Q\in\cp_{J_1};\po(Q,F(Q))\in zW_{\boc(J)}\}$. Then
$f:V@>>>V'$, $P\m P^1:=P^{F(P)}$ is an isomorphism.}
\nl
Define $V'@>>>\cp_J$ by $Q\m P$ where $P$ is the unique parabolic in $\cp_J$ such that 
$Q\sub P$. We have automatically $P\in V$ hence $Q\m P$ is a map $f':V'@>>>V$. Clearly 
$f'f=1$. We show $ff'=1$. It is enough to show that, if $Q,P$ are as above, then 
$P^{F(P)}=Q$. We have $\po(Q,F(Q))=zu$ where $u\in W_{\boc(J)}$. We can find 
$B_0,B_1\in\cb$ 
such that $B_0\sub Q$, $B_1\sub F(Q)$, $\po(B_0,B_1)=zu$. Since $l(zu)=l(z)+l(u)$ we can 
find $B_2\in\cb$ such that $\po(B_0,B_2)=z,\po(B_2,B_1)=u$. Since $u\in W_{\boc(J)}$ and 
$B_1\sub F(P)$ we have $B_2\sub F(P)$. Since $B_0\sub P,B_2\sub F(P)$, $\po(B_0,B_2)=z$, we
have $B_0\sub P^{F(P)}$. Since $Q,P^{F(P)}$ are in $\cp_{J_1}$ and both contain $B_0$ we 
have $Q=P^{F(P)}$. This proves (c).

Let $\tt=(J_n,w_n)_{n\ge0}\in\ct'(J,\boc)$. For $m\ge0$ we set 
$\tt_m=(J'_n,w'_n)_{n\ge0}$ where $J'_n=J_{n+m},w'_n=w_{n+m}$. We have 
$\tt_m\in\ct'(J_m,\boc)$. We set $\tt_\iy=(J'_n,w'_n)_{n\ge0}$ where 
$J'_n=J_\iy,w'_n=w_\iy$. We have $\tt_\iy\in\ct'(J_\iy,\boc)$. Clearly, $P\m P^1$ is a map 
$$\vt:\cp^\tt_J@>>>\cp^{\tt_1}_{J_1}.$$
 
(d) {\it The map $P\m P^1$ is an isomorphism $\cp^\tt_J@>\si>>\cp^{\tt_1}_{J_1}$. The map 
$P\m P^\iy$ is an isomorphism $\cp^\tt_J@>\si>>\cp^{\tt_\iy}_{J_\iy}$.}
\nl
The first assertion of (d) follows from (c). The second assertion follows using the first 
assertion repeatedly.

(e) {\it Let $\tt=(J_n,w_n)_{n\ge0}\in\ct'(J,\boc)$ be such that $J_n=J$ and $w_n=w$ for 
all $n\ge0$ where $w\in W$. We have $\boc(J)=w\i Jw$, 
$w\in{}^JW^{\boc(J)}$. If $P\in\cp_J^\tt$ then
$P^n=P$ for $n\ge0$ and $\po(P,F(P))=w$. From $P=P^{F(P)}$ we see that $P,F(P)$ have a 
common Levi. We have $\cp_J^\tt=\{P\in\cp_J;\po(P,F(P))=w\}$.}

(f) {\it Let $(J_n,w_n)_{n\ge0}\in\ct'(J,\boc)$. For $n\ge0$ we have 
$w_n=\min(W_Jw_\iy W_{\boc(J_n)})$.}

\subhead 4.3\endsubhead
In the setup of 4.2(e) we show:
$$\dw\i L_J\dw=L_{\boc(J)}=F(L_J).\tag a$$
From $\po(P_J,\dw P_{\boc(J)}\dw\i)=w$ 
we see that $P_J,\dw P_{\boc(J)}\dw\i$ have a common Levi
subgroup containing $T$ which must be $L_J$ and also $\dw L_{\boc(J)}\dw\i$.

Let 
$$\tcp_J^\tt=\{gU_{P_J}\in G/U_{P_J};g\i F(g)\in U_{P_J}\dw U_{P_{\boc(J)}}\}.$$
Define $F':L_J@>>>L_J$ by $g\m\dw F(g)\dw\i$. This is the Frobenius map for an 
$\FF_q$-rational structure on $L_J$. We set
$$L_J^{F'}=\{l\in L_J;F'(l)=l\}.$$
The finite group $L_J^{F'}$ acts freely on $\tcp_J^\tt$ by $l:gU_{P_J}\m gl\i U_{P_J}$ and the map $f:\tcp_J^\tt@>>>\cp_J^\tt$,
$gU_{P_J}\m gP_Jg\i$ is constant on the orbits of this action. We show:

(b) {\it $f$ is a principal $L_J^{F'}$-bundle.}
\nl
We only show this at the level of sets. If $P\in\cp_J;\po(P,F(P))=w$, we have $P=gP_Jg\i$ 
where $g\in G$ satisfies 
$$\align w=&\po(gP_Jg\i,F(g)F(P_J)F(g\i))=\po(P_J,g\i F(g)P_{\boc(J)}F(g\i)g)\\&
=\po(P_J,\dw P_{\boc(J)}\dw\i).\endalign$$
Hence there exists $y\in P_J$ such that 
$g\i F(g)P_{\boc(J)}F(g\i)g=y\dw P_{\boc(J)}\dw\i y\i$ 
hence $g\i F(g)\in P_J\dw P_{\boc(J)}$ that is $g\i F(g)\in l'U_{P_J}\dw U_{P_{\boc(J)}}$ 
for some $l'\in L_J$. (We use (a).) By Lang's theorem for $F'$ we can find $l\in L_J$ such
that $l\i F'(l)=l'$. Then $gl\i U_{P_J}\in\tcp_J^\tt$. We see that $f$ is surjective.

Assume that $gU_{P_J},g'U_{P_J}$ in $\tcp_J^\tt$ have the same image under $f$ that is 
$gP_Jg\i=g'P_Jg'{}\i$. Then $g'=gp\i$ where $p\in P_J$. We may assume that $g'=gl\i$, 
$l\in L_J$. We have $g\i F(g)\in P_J\dw P_{\boc(J)}$ and  
$(gl\i)\i F(gl\i)\in U_{P_J}\dw U_{P_{\boc(J)}}$ that is 
$$g\i F(g)=U_{P_J}l\i\dw F(l)U_{P_{\boc(J)}}\text{ and }
U_{P_J}l\i F'(l)\dw U_{P_{\boc(J)}}=U_{P_J}\dw U_{P_{\boc(J)}}.$$
Using \cite{\GR, 3.2} we deduce $l\i F'(l)\dw=\dw$ hence $l\in L_J^{F'}$. 

Let 
$${}'\tcp_J^\tt=\{g(U_{P_J}\cap F\i(\dw\i U_{P_J}\dw))\in 
G/(U_{P_J}\cap F\i(\dw\i U_{P_J}\dw));g\i F(g)\in U_{P_J}\dw\}.$$
We show:

(c) {\it The map $g(U_{P_J}\cap F\i(\dw\i U_{P_J}\dw))\m gU_{P_J}$ is an isomorphism
$\g:{}'\tcp_J^\tt@>\si>>\tcp_J^\tt$.}
\nl
We only show this at the level of sets. Let $gU_{P_J}\in\tcp_J^\tt$. We have
$g\i F(g)=u\dw F(u')$ for some $u\in U_{P_J},u'\in U_{P_J}$. Then 
$(gu'{}\i)\i F(gu'{}\i) =u'u\dw$ so that 
$\g(gu'{}\i(U_{P_J}\cap F\i(\dw\i U_{P_J}\dw)))=gU_{P_J}$.
We see that $\g$ is surjective. The injectivity is immediate.

Next we show:

(d) {\it $\cp_J^\tt$ is a smooth variety of pure dimension equal to 
\lb $\dim U_{P_J}-\dim(U_{P_J}\cap F\i(\dw\i U_{P_J}\dw))$ and its connected components are
permuted transitively by the $\G$-action on $\cp_J^\tt$.}
\nl
By (b),(c) it is enough to show that ${}'\tcp_J^\tt$ is smooth, of pure dimension equal to
$\dim U_{P_J}-\dim(U_{P_J}\cap F\i(\dw\i U_{P_J}\dw))$ and its connected components are
permuted transitively by the $\G$-action 

$g_0:g(U_{P_J}\cap F\i(\dw\i U_{P_J}\dw))\m g_0g(U_{P_J}\cap F\i(\dw\i U_{P_J}\dw))$ 
\nl
on ${}'\tcp_J^\tt$. This follows from the fact that $\{g\in G;g\i F(g)\in U_{P_J}\dw\}$ is
smooth of dimension $\dim U_{P_J}$ and $U_{P_J}\dw$ is connected.

\subhead 4.4\endsubhead
We now consider a general $\tt=(J_n,w_n)_{n\ge0}\in\ct'(J,\boc)$. Let 
$\cp^\tt_J@>\si>>\cp^{\tt_\iy}_{J_\iy}$ be the isomorphism in 4.2(d). By 4.3 for $\tt_\iy$
instead of $\tt$, $\tcp^{\tt_\iy}_{J_\iy}$ is defined. Let
$$\tcp_J^\tt:=\tcp^{\tt_\iy}_{J_\iy}=\{gU_{P_{J_\iy}}\in G/U_{P_{J_\iy}};
g\i F(g)\in U_{P_{J_\iy}}\dw U_{P_{\boc(J_\iy)}}\}.$$
Define $F':L_{J_\iy}@>>>L_{J_\iy}$ by $l\m\dw_\iy F(l)\dw_\iy\i$. (We have 
$\dw_\iy\i L_{J_\iy}\dw_\iy=L_{\boc(J_\iy)}=F(L_{J_\iy})$, see 4.3(a).)
Let $\L=L_{J_\iy}^{F'}$. The finite group $\G\T\L$ acts on $\tcp_J^\tt$ by 
$$(g_0,l):gU_{P_{J_\iy}}\m g_0gl\i U_{P_{J_\iy}}$$
and on $\cp_J^\tt$ by $(g_0,l):P\m g_0Pg_0\i$. From 4.3(b) we see that 

(a) {\it The map $f:\tcp_J^\tt@>>>\cp_J^\tt$, $gU_{P_{J_\iy}}\m gP_Jg\i$ (which is 
compatible with the $\G\T\L$-actions) is a principal $\L$-bundle.}
\nl
From 4.3(d) we see that:

(b) {\it $\cp_J^\tt$ is a smooth variety of pure dimension equal to \lb 
$\dim U_{P_{J_\iy}}-\dim(U_{P_{J_\iy}}\cap F\i(\dw_\iy\i U_{P_{J_\iy}}\dw_\iy))$ and its 
connected components are permuted transitively by the $\G$-action on $\cp_J^\tt$.}
\nl
Let $M$ be a finite dimensional $\L$-irreducible module over $\bbq$. We view 
$M$ as a $\G\T\L$-module with $\G$ acting trivially and we form the 
$\G\T\L$-equivariant local system $M_{\tcp_J^\tt}$ on $\tcp_J^\tt$ as in 1.6.
Using (a) we see that there is a well defined $\G\T\L$-equivariant local 
system $\uM$ on $\cp_J^\tt$ with trivial action of $\L$ such that 
$f^*\uM=M_{\tcp_J^\tt}$ as $\G\T\L$-equivariant local systems. We will regard
$\uM$ as a $\G$-equivariant local system. 

Let $d=\dim\cp_J^\tt$ (see (b)). We show:

(c) {\it $\uM[d]$ (an object of $\cm_\G(\cp_J^\tt)$ by (b)) is simple.}
\nl
Let $r:\tC@>>>C$ be a finite principal covering with finite group $H$. Assume that $C$ is
connected. There is an obvious functor $E\m E'$ from $H$-modules of finite dimension over 
$\bbq$ to local systems on $C$ which are direct summands of $r_!\bbq$. If $E$ is 
irreducible then $E'$ is irreducible as a local system.

We apply this statement in the case where $C$ is a connected component of $\cp_J^\tt$, 
$\tC=f\i(C)$ ($f$ as in (a)), $r$ is the restriction of $f$, $H=\L$ and $E=M$. 
Note that $E'=\uM|_C$. We see that the local system $\uM|_C$ is irreducible. It remains to
use the transitivity statement in (b).

For an object $A$ of $\cm_\G(\cp_J^\tt)$ we write $A\in\Bbb S(\cp_J^\tt)$ instead of
"$A$ is isomorphic to $\uM[d]$ for some $M$ as above."

From (c) we see that:

(d) {\it $\uM^\sh[d]$ is a simple object of $\cm_\G(\cp_J)$.}

For an object $A$ of $\cm_\G(\cp_J)$ we write $A\in\Bbb S(\cp_J)$ instead of "$A$ is 
isomorphic to $\uM^\sh[d]\in\cm_\G(\cp_J)$ for some $\tt\in\ct'(J,\boc)$ and some $M$ as 
above."

\subhead 4.5\endsubhead
For $w\in W$ we identify $\cb_w$ with $Z^{(w)}$ as in 2.5. We show:

(a) {\it Let $a_1,a_2,b\in W$. Let $\cl\in\cs(T)^{a_1a_2F}$.
Let $V$ be a locally closed $\G$-stable subvariety of $\cb_b$. Let 
$X_{a_1,a_2}=\{(B_0,B_1,B_2)\in Z^{(a_1,a_2)};B_1\in V\}$, (see 2.5). Define 
$\k:X_{a_1,a_2}@>>>V$ by $(B_0,B_1,B_2)\m B_1$. Let $\ce=\cl_{(a_1,a_2)}|_{X_{a_1,a_2}}$.
Let $V'$ be an algebraic variety with a $\G$-action and let $m:V@>>>V'$ be a morphism 
compatible with the $\G$-actions. Let $A'$ be a simple object of $\cm_\G(V')$ such that 
$A'\dsv_\G{}^pH^{\cdot}((m\k)_!\ce)$. Then there exists $e\in W$ such that 
$(bF)^*(e^*\cl)\cong e^*\cl$ and $A'\dsv_\G{}^pH^{\cdot}(m_!(e^*\cl)_b|_V)$.}
\nl
We argue by induction on $l(a_1)$. If $l(a_1)=0$ then $a_2=b$, $\k$ is an isomorphism and 
the result is obvious (with $e=1$). Assume now that $l(a_1)>0$. We can find $s\in\II$ such
that $l(a_1)>l(sa_1)$. Let $\ce_1=(s^*\cl)_{(sa_1,a_2\boc(s)}|_{X_{sa_1,a_2\boc(s)}}$.

Assume first that $l(a_2\boc(s))=l(a_2)+1$. We have an isomorphism 
$\io:X_{a_1,a_2}@>>>X_{sa_1,a_2\boc(s)}$, $(B_0,B_1,B_2)\m(B'_0,B_1,F(B'_0))$ where 
$B'_0\in\cb$ is defined by 
$$\po(B_0,B'_0)=s,\po(B'_0,B_1)=sa_1.\tag b$$
Define $\k':X_{sa_1,a_2\boc(s)}@>>>V$ by $(B_0,B_1,B_2)\m B_1$. We have $\k=\k'\io$,
$\io^*\ce_1=\ce$ hence $\io_!\ce=\ce_1$. Thus
$(m\k)_!\ce=(m\k')_!\ce_1$ and $A'\dsv_\G{}^pH^{\cdot}((m\k')_!\ce_1)$. By the 
induction hypothesis there exists $e'\in W$ such that 
$(bF)^*(e'{}^*s^*\cl)\cong e'{}^*s^*\cl$ and 
$A'\dsv_\G{}^pH^{\cdot}(m_!(e'{}^*s^*\cl)_b|_V)$. The result follows with $e=se'$.

Assume next that  $l(a_2\boc(s))=l(a_2)-1$. We have a partition
$X_{a_1,a_2}=X'\cup X''$ where $X'$ (resp. $X''$)
is the open (resp. closed) subset of $X_{a_1,a_2}$ defined by $\po(B_1,F(B'_0))=a_2$ (resp.
$\po(B_1,F(B'_0))=a_2\boc(s)$). Let $j'=\k|_{X'}$, $j''=\k|_{X''}$. By 
general principles we have either

(c) $A'\dsv_\G{}^pH^{\cdot}((mj')_!(\ce|_{X'}))$ or

(d) $A'\dsv_\G{}^pH^{\cdot}((mj'')_!(\ce|_{X''}))$.
\nl
Assume that (d) holds. We have $j''=\k''\io''$ where $\k'':X_{sa_1,a_2\boc(s)}@>>>V$ is 
given
by $(B_0,B_1,B_2)\m B_1$ and $\io'':X''@>>>X_{sa_1,a_2\boc(s)}$ is 
$(B_0,B_1,B_2)\m(B'_0,B_1,F(B'_0))$ with $B'_0$ as in (b). We have 
$\ce|_{X''}=\io''{}^*\ce_1$. Now $\io''$ is an affine line bundle hence 
$\io''_!(\ce|_{X''})=\ce_1[-2](-1)$. Hence $A'\dsv_\G{}^pH^{\cdot}((\m\k'')_!\ce_1)$. 
By the induction hypothesis there exists $e'\in W$ such that 
$(bF)^*(e'{}^*s^*\cl)\cong e'{}^*s^*\cl$ and 
$A'\dsv_\G{}^pH^{\cdot}(m_!(e'{}^*s^*\cl)_b|_V)$.
The result follows with $e=se'$.

Assume now that (c) holds. We have $j'=\k'\io'$ where $\k':X_{sa_1,a_2}@>>>V$ is 
$(B_0,B_1,B_2)\m B_1$ and $\io':X'@>>>X_{sa_1,a_2}$ is 
$(B_0,B_1,B_2)\m(B'_0,B_1,F(B'_0))$ with $B'_0$ as in (b). Note that $\io'$ makes 
$X'$ into the complement of a section of an affine line bundle over 
$X_{sa_1,a_2}$. If $s\n W_\cl$ then by an argument as in the proof of 2.15 we see that
$\io'_!(\ce|_{X'})=0$ contradicting (c). Thus we may assume that $s\in W_\cl$. 
Then $\ce_2=\cl_{(sa_1,a_2)}|_{X_{sa_1,a_2}}$ is defined and
$\ce|_{X'}=\io'{}^*\ce_2$. Hence we have a distinguished triangle  
$(\io'_!\ce_{X'},\ce_2,\ce_2[-2](-1))$ hence a distinguished triangle 
$(m_!j'_!(\ce|_{X'},m_!\k'_!\ce_2,m_!\k'_!\ce_2[-2](-1))$. It follows that
$A'\dsv_\G{}^pH^{\cdot}(m_!\k'_!\ce_2)$. 
By the induction hypothesis there exists $e\in W$ such that 
$(bF)^*(e^*\cl)\cong e^*\cl$ and $A'\dsv_\G{}^pH^{\cdot}(m_!(e^*\cl)_b|_V)$. This 
completes the proof of (a).

\subhead 4.6\endsubhead
Let $\tt=(J_n,w_n)_{n\ge0}\in\ct'(J,\boc)$. Let 
$\cb_\tt=\{B'\in\cb;P_{B',J}\in\cp_J^\tt\}$. For $a\in W$ let 
$\cb_{\tt,a}=\cb_\tt\cap\cb_a$. Define $\x_{\tt,a}:\cb_{\tt,a}@>>>\cp_J^\tt$ by 
$B'\m P_{B',J}$. We show:

(a) {\it Let $\cl\in\cs(T)^{aF}$. Let $A$ be a simple object of $\cm_\G(\cp_J^\tt)$ such 
that $A\dsv_\G{}^pH^{\cdot}(\x_{\tt,a!}(\cl_a|_{\cb_{\tt,a}}))$. 
Then there exist $b,e\in W$ such that $b^*e^*\cl\cong e^*\cl$ and
$A\dsv_\G{}^pH^{\cdot}(\vt^*\x_{\tt_1,b!}((e^*\cl)_b|_{\cb_{\tt_1,b}}))$.}
\nl
Since $\x_{\tt,a!}(\cl_a|_{\cb_{\tt,a}})\ne0$ we have $\cb_{\tt,a}\ne\em$. Thus there 
exists $B'\in\cb$ such that $\po(B',F(B'))=a$, $P_{B',J}\in\cp_J^\tt$. We have 
$\po(P_{B',J},F(P_{B',J}))=w_0$. Since $B'\sub P_{B',J}$, $F(B')\sub F(P_{B',J})$, it 
follows that $a\in W_Jw_0W_{\boc(J)}$ and $w_0=\min(W_JaW_{\boc(J)})$.

Define $\ph:\cb_{\tt,a}@>>>\cb_{\tt_1}$ by $\ph(B')=(P_{B',J})^{F(B')}$. (For 
$B'\in\cb_\tt$ we have $(P_{B',J})^{F(B')}\in\cb_{\tt_1}$ since
$(P_{B',J})^{F(B')}\sub(P_{B',J})^{F(P_{B',J})}\in\cp^{\tt_1}_{J_1}$.) We have a partition 
$\cb_{\tt_1}=\sqc_{b\in W}\cb_{\tt_1,b}$. Setting $\cb_{\tt,a,b}=\ph\i(\cb_{\tt_1,b})$ we
get a partition $\cb_{\tt,a}=\sqc_{b\in W}\cb_{\tt,a,b}$. Let 
$\x_{\tt,a,b}:\cb_{\tt,a,b}@>>>\cp_J^\tt$ be the restriction of $\x_{\tt,a}$. By general 
principles we have $A\dsv_\G{}^pH^{\cdot}(\x_{\tt,a,b!}(\cl_a|_{\cb_{\tt,a,b}}))$ for 
some $b\in W$. Let
$\ph_b:\cb_{\tt,a,b}@>>>\cb_{\tt_1,b}$ be the restriction of $\ph$. We have
$\vt\x_{\tt,a,b}=\x_{\tt_1,b}\ph_b$ (both compositions carry $B'$ to 
$(P_{B',J})^{F(P_{B',J})}$.) Hence $\x_{\tt,a,b}=\vt\i\x_{\tt_1,b}\ph_b$. Thus, 
$A\dsv_\G{}^pH^{\cdot}((\vt\i)_!\x_{\tt_1,b!}\ph_{b!}(\cl_a|_{\cb_{\tt,a,b}}))$ and
$$\vt_!A\dsv_\G{}^pH^{\cdot}(\x_{\tt_1,b!}\ph_{b!}(\cl_a|_{\cb_{\tt,a,b}})).\tag b$$
\nl
We can write uniquely $a=a_1a_2$ where $a_1\in W_J,a_2\in{}^JW$. We show that for any
$\tB\in\cb_{\tt_1,b}$ we have
$$\ph_b\i(\tB)=\{B'\in\cb;\po(B',\tB)=a_1,\po(\tB,F(B'))=a_2\}.\tag c$$
Assume first that $B'\in\ph_b\i(\tB)$. We know that $\po(B',F(B'))=a$. We have 
$\po(B',(P_{B',J})^{F(B')})\in W_J$ since $B',(P_{B',J})^{F(B')}$ are two Borel subgroups
of $P_{B',J}$. From the definitions we have $\po((P_{B',J})^{F(B')},F(B'))\in{}^JW$. We 
have automatically $\po(B',(P_{B',J})^{F(B')})=a_1$, $\po((P_{B',J})^{F(B')},F(B'))=a_2$ 
that is $\po(B',\tB)=a_1$, $\po(\tB,F(B'))=a_2$.  

Conversely, assume that $B'$ belongs to the right hand side of (c). We have 
$l(a_1a_2)=l(a_1)+l(a_2)$ hence $\po(B',F(B'))=a_1a_2=a$. Since the properties 
$\po(B',\tB)\in W_J$, $\po(\tB,F(B'))\in{}^JW$ characterize $\tB$ and 
$(P_{B',J})^{F(B')}$ has the same properties, it follows that $P^{F(B')}=\tB$ where 
$P=P_{B',J}$. Since $B'\sub P,F(B')\sub F(P)$ we have 
$\po(P,F(P))=\min(W_J\po(B',F(B'))W_{\boc(J)})=\min(W_JaW_{\boc(J)})=w_0$.
It follows that $P^{F(P)}\in\cp_{J\cap w_0\boc(J)w_0\i}=J_1$. Clearly,
$\tB=P^{F(B')}\sub P^{F(P)}$. Hence $P^{F(P)}=P_{\tB,J_1}\in\cp_{J_1}^{\tt_1}$. It follows 
that $P\in\cp_J^\tt$. Thus, $B'\in\cb_{\tt,a}$ and $\ph(B')=\tB$. Since 
$\tB\in\cb_{\tt_1,b}$ we see that $B'\in\ph_b\i(\tB)$. This proves (c).

From (c) we see that $(B_0,B_1,B_2)\m B_0$ is an isomorphism $X_{a_1,a_2}@>>>\cb_{\tt,a,b}$
where $X_{a_1,a_2}$ is defined as in 4.5(a) in terms of $V=\cb_{\tt_1,b}$. Under this
isomorphism, $\ph_b$ corresponds to $X_{a_1,a_2}@>>>\cb_{\tt_1,b}$, $(B_0,B_1,B_2)\m B_1$.
Applying 4.5(a) with $V'=\cp_{J_1}^{\tt_1}$, $m=\x_{\tt_1,b}$, $A'=\vt_!A$ we see that 
there exists $e\in W$ such that $(bF)^*(e^*\cl)\cong e^*\cl$ and
$\vt_!A\dsv_\G{}^pH^{\cdot}(\x_{\tt_1,b!}((e^*\cl)_b|_{\cb_{\tt_1,b}}))$. (The 
assumption of 4.5(a) is verified by (b).) Since $\vt$ is an isomorphism, it follows that 
$A\dsv_\G{}^pH^{\cdot}(\vt^*\x_{\tt_1,b!}((e^*\cl)_b|_{\cb_{\tt_1,b}}))$. This 
proves (a).

\subhead 4.7\endsubhead
Let $\tt=(J_n,w_n)_{n\ge0}\in\ct'(J,\boc)$. Let $d=\dim\cp_J^\tt$. Let $a\in W$. We show:

(a) {\it Let $\cl\in\cs(T)^{aF}$. Let $A$ be a simple object of
$\cm_\G(\cp_J^\tt)$ such that $A\dsv_\G{}^pH^{\cdot}(\x_{\tt,a!}(\cl_a|_{\cb_{\tt,a}}))$. 
Then $A\in\Bbb S(\cp_J^\tt)$.}
\nl
More generally we show that (a) holds when $J,\tt$ are replaced by $J_n,\tt_n$, $n\ge0$. 
First we show:

(b) if the result holds for $n=1$ then it holds for $n=0$.
\nl
Let $A$ be as in (a). By 4.6(a) there exist $b,e\in W$ such that 
$(bF)^*(e^*\cl)\cong e^*\cl$ and
$A\dsv_\G{}^pH^{\cdot}(\vt^*\x_{\tt_1,b!}((e^*\cl)_b|_{\cb_{\tt_1,b}}))$. Since 
$\vt$ is an isomorphism, there exists a simple object $A'$ of $\cm_\G(\cp_{J_1}^{\tt_1})$ 
such that $A=\vt^*A'$. From our assumption we have 
$A'\dsv_\G{}^pH^{\cdot}(\x_{\tt_1,b!}((e^*\cl)_b|_{\cb_{\tt_1,b}}))$. Since (b) 
holds for $n=1$ we have $A'\cong\uM[d]$ for some irreducible $L_{J_\iy}^{F'}$-module $M$. 
Hence $A$ is of the same form. Thus (b) holds.

Similarly, if the result holds for some $n\ge1$ then it holds for $n-1$. In this way we see
that it suffices to prove the result for $n$ large. Thus in the remainder of this proof we
assume, as we may, that $J_0=J_1=\do=J$ and $w_0=w_1=\do=w$. We can write uniquely 
$a=a_1a_2$ where $a_1\in W_J,a_2\in{}^JW$. Since $\x_{\tt,a!}(\cl_a|_{\cb_{\tt,a}})\ne0$,
we have $\cb_{\tt,a}\ne\em$. From this we deduce as in the proof of 4.6(a) that 
$a\in W_JwW_{\boc(J)}$. Since $w\boc(J)w\i=J$ we must have $wW_{\boc(J)}=W_Jw$ so that 
$a\in W_Jw$. Since $w\in{}^JW$ (see 4.2(b)) it follows that $w=a_2$. In particular we have
$a_2\in{}^JW^{\boc(J)}$. Hence if $B''\in\cb_a$ then $\po(P_{B'',J},F(P_{B'',J}))=a_2=w$. 
Since in our case $\cp_J^\tt=\{P\in\cp_J;\po(P,F(P))=w\}$ (see 4.2(e)) we see that 
$P_{B'',J}\in\cp_J^\tt$. Thus we have $\cb_a=\cb_{\tt,a}$ and
$A\dsv_\G{}^pH^{\cdot}(\x_{\tt,a!}(\cl_a))$. Let 
$\ti\cb_a=\{gU\in G/U;g\i F(g)\in U\doa U\}$. Recall from 4.3(b) that 
$f_a:\ti\cb_a@>>>\cb_a$, $gU\m gBg\i$ is a finite principal covering with group 
$\fT=\{t\in T;\doa F(t)\doa\i=t\}$. From 2.5 we see that $\cl_a$ is a direct summand of 
$f_{a!}\bbq$. Thus we have $A\dsv_\G{}^pH^{\cdot}(\x_{\tt,a!}f_{a!}\bbq)$.

By 4.3(a), $L_J=\dw L_{\boc(J)}\dw\i$ is a common Levi subgroup of 
$P_J,\dw P_{\boc(J)}\dw\i$.  Let $F':L_J@>>>L_J$ be as in 4.3. Let $\L=L_J^{F'}$.
Let $\cb'$ be the variety of Borel subgroups
of $L_J$. For $\b,\b'$ in $\cb'$ we have $\po(\b U_{P_J},\b'U_{P_J})=y$ for a unique 
$y\in W_J$; we then also write $\po'(\b,\b')=y$. Let 
$Y'=\{\b\in\cb';\po'(\b,F'(\b))=a_1\}$. Now $\L$ acts on $Y'$ by conjugation. Let 
$\tcp_J^\tt$ be as in 4.3. Note that $\L$ acts (freely) on $\tcp_J^\tt\T Y'$ by
$l:(gU_{P_J},\b)\m(gl\i U_{P_J},l\b l\i)$ and we can form the orbit space
$\L\bsl(\tcp_J^\tt\T Y')$.

We define $\psi:\tcp_J^\tt\T Y'@>>>\cb_a$ by $(gU_{P_J},\b)\m g\b U_{P_J}g\i$. We 
show that $\psi$ is well defined, that is 
$\po(g\b U_{P_J}g\i,F(g)F(\b)U_{P_{\boc(J)}}F(g\i))=a$
for $gU_{P_J}\in\tcp_J^\tt$. Since $a=a_1w$, $l(a)=l(a_1)+l(w)$, it is enough to show that
for $(gU_{P_J},\b)\in\tcp_J^\tt\T Y'$ we have

(c) $\po(g\b U_{P_J}g\i,gF'(\b)U_{P_J}g\i)=a_1$,

(d) $\po(gF'(\b)U_{P_J}g\i,F(g)F(\b)U_{P_{\boc(J)}}F(g\i))=w$.
\nl
Now the left hand side of (c) is equal to 
$\po(\b U_{P_J},F'(\b)U_{P_J})=\po'(\b,F'(\b))=a_1$, proving (c). We have 
$g\i F(g)=u\dw u'$ where $u\in U_{P_J},u'\in U_{P_{\boc(J)}}$. 
The left hand side of (d) is 
equal to
$$\align&\po(F'(\b)U_{P_J},g\i F(g)F(\b)U_{P_{\boc(J)}}F(g\i)g)\\&=
\po(F'(\b)U_{P_J},u\dw u'F(\b)U_{P_{\boc(J)}}u'{}\i\dw\i u\i)=
\po(F'(\b)U_{P_J},\dw F(\b)U_{P_{\boc(J)}}\dw\i)\\&=
\po(F'(\b)U_{P_J},F'(\b)U_{\dw P_{\boc(J)}\dw\i})=\po(P_J,\dw P_{\boc(J)}\dw\i).
\endalign$$
(The last equality follows from \cite{\PAR, II, 8.3}.) This equals $w$. 

We see that $\psi$ is well defined. We show:

(e) {\it The map $\bar\psi:\L\bsl(\tcp_J^\tt\T Y')@>>>\cb_a$ induced by $\psi$ is an 
isomorphism.}
\nl
Let $B'\in\cb_a$. We can find uniquely $\tB\in\cb$ such that $\po(B',\tB)=a_1$,\lb
$\po(\tB,F(B'))=w$. Then $P_{B',J}=P_{\tB,J}$ and $\po(P_{\tB,J},P_{F(B'),\boc(J)})=w$. 
Hence $\po(P_{B',J},F(P_{B',J}))=w$ that is, $P_{B',J}\in\cp_J^\tt$. Since 
$\tcp_J^\tt@>>>\cp_J^\tt$ is surjective (see 4.3) we can find $gU_{P_J}\in\tcp_J^\tt$ such
that $gP_Jg\i=P_{B',J}$; note that $gU_{P_J}$ is unique up to the action of $\L$. 
Since $B'\sub P_{B',J}$ we have $g\i B'g\sub P_J$ hence there is a unique $\b\in\cb'$ such
that $g\i B'g=\b U_{P_J}$. As above we see that 
$$\po(gF'(\b)U_{P_J}g\i,F(g)F(\b)U_{P_{\boc(J)}}F(g\i))=w$$
that is 
$\po(gF'(\b)U_{P_J}g\i,F(B'))=w$. Thus $gF'(\b)U_{P_J}g\i$ is a Borel subgroup of 
$P_{B',J}$ whose relative position with $F(B')$ is $w$. But there is only one such Borel
subgroup. Therefore $gF'(\b)U_{P_J}g\i=\tB$. Since $\po(B',\tB)=a_1$ we have
$\po(g\b U_{P_J}g\i,gF'(\b)U_{P_J}g\i)=a_1$ hence $\po(\b U_{P_J},F'(\b)U_{P_J})=a_1$ that
is, $\po'(\b,F'(\b))=a_1$. Thus to $B'\in\cb_a$ we have associated 
$(gU_{P_J},\b)\in\tcp_J^\tt\T Y'$; its $\L$-orbit is well defined. Thus we have a 
well defined map $\cb_a@>>>\L\bsl(\tcp_J^\tt\T Y')$. Clearly, this is the inverse of 
$\bar\psi$. This proves (e).

Now let $\b'=B\cap L_J$, and let $U'$ be the unipotent radical of $\b'$. Let 
$U''=U'\cap F'{}\i(\doa_1\i U'\doa_1)$. Let $\tY'=\{lU''\in L_J/U'';l\i F'(l)\in U'\doa_1\}$.
As in 4.3(b),(c) the map $\r:\tY'@>>>Y'$, $lU''\m l\b'l\i$ is a finite principal covering 
with group $\{t\in T;\doa_1F'(t)\doa_1\i=t\}=\fT$. Now $\L$ acts freely on 
$\tcp_J^\tt\T\tY'$ by
$$l_0:(gU_{P_J},lU'')\m(gl_0\i U_{P_J},l_0lU'')$$
and we can form the orbit space $\L\bsl(\tcp_J^\tt\T\tY')$.

The map $\tcp_J^\tt\T\tY'@>>>\ti\cb_a$, $(gU_{P_J},lU'')\m glU$ induces a 
map $\ti\psi:\L\bsl(\tcp_J^\tt\T\tY')@>>>\ti\cb_a$ which is easily seen to be an
isomorphism. We have a commutative diagram
$$\CD
\L\bsl(\tcp_J^\tt\T\tY')@>\ti\psi>>\ti\cb_a\\
@V\x VV                                   @Vf_aVV     \\
\L\bsl(\tcp_J^\tt\T Y')@>\bar\psi>>\cb_a
\endCD$$
where $\x$ is induced by $(gU_{P_J},lU'')\m(gU_{P_J},l\b'l\i)$ (a principal $\fT$-bundle).
Note that the horizontal maps in this diagram are isomorphisms.

Under the isomorphism (e), the map $\x_{\tt,a}$ becomes the map 
$\x':\L\bsl(\tcp_J^\tt\T Y')@>>>\cp_J^\tt$ induced by $(gU_{P_J},\b)\m gP_Jg\i$. It 
follows that $\x_{\tt,a!}f_{a!}\bbq=(\x'\x)_!\bbq$ and 
$A\dsv_\G{}^pH^{\cdot}((\x'\x)_!\bbq)$. We extend the natural $\G$-actions on 
$\cp_J^\tt$ and on $\L\bsl(\tcp_J^\tt\T\tY')$ to $\G\T\L$-actions with 
$\L$ acting trivially. Then $A,{}^pH^j((\x'\x)_!\bbq)$ are naturally objects of 
$\cm_{\G\T\L}(\cp_J^\tt)$ and $A\dsv_{\G\T\L}{}^pH^{\cdot}((\x'\x)_!\bbq)$.
Let $f:\tcp_J^\tt@>>>\cp_J^\tt$ be as in 4.3. Now $\G\T\L$ acts on $\tcp_J^\tt$ as 
in 4.4 (compatibly with $f$). Moreover $f^*A$, 
$f^*({}^pH^j((\x'\x)_!\bbq))={}^pH^j(f^*(\x'\x)_!\bbq)$ are objects of 
$\cm_{\G\T\L}(\tcp_J^\tt)$ and $f^*A\dsv_{\G\T\L}{}^pH^{\cdot}(f^*(\x'\x)_!\bbq)$.
Let $p_1:\tcp_J^\tt\T\tY'@>>>\tcp_J^\tt$ be the first projection. Now $\G\T\L$ acts
on $\tcp_J^\tt\T\tY'$ by $(g_0,l_0):(gU_J,lU'')\m(g_0gl_0\i U_J,l_0lU'')$ and $p_1$ is 
compatible with the $\G\T\L$-actions.
 We have $f^*(\x'\x)_!\bbq=p_{1!}\bbq$. We see that
$f^*A\dsv_{\G\T\L}{}^pH^{\cdot}(p_{1!}\bbq)$.
Let $p_2:\tcp_J^\tt\T\tY'@>>>\tY'$ be the second projection. Now $\G\T\L$ acts on 
$\tY'$ by $(g_0,l_0):lU''\m l_0lU''$ and $p_2$ is compatible with the
$\G\T\L$-actions. The obvious maps $\tcp_J^\tt@>e>>\text{point}@<e'<<\tY'$ are 
again compatible with
the $\G\T\L$-actions (the action on the point is trivial). We have
$p_{1!}\bbq=e^*e'_!\bbq$ hence $f^*A\dsv_{\G\T\L}{}^pH^{\cdot}(e^*e'_!\bbq)$.
We have a spectral sequence in $\cm_{\G\T\L}(\tcp_J^\tt)$ with 
$E_2={}^pH^{\cdot}(e^*{}^pH^{\cdot}(e'_!\bbq))$ and 
$E_\iy$ is an associated graded of ${}^pH^{\cdot}(e^*e'_!\bbq)$. We have
$f^*A\dsv_{\G\T\L}E_\iy$ hence $f^*A\dsv_{\G\T\L}E_2$. Now 
${}^pH^{\cdot}(e'_!\bbq)$ is just a $\G\T\L$-module $M$ with trivial action of $\G$
and
$e^*{}^pH^{\cdot}(e'_!\bbq)$ is the local system $f^*\uM$ (notation of 4.4). Since 
$\tcp_J^\tt$
is smooth of pure dimension $d$, ${}^pH^j(e^*{}^pH^{\cdot}(e'_!\bbq))$ is $0$ if $j\ne d$ 
and is $f^*\uM[d]$ if $j=d$. Thus $E_2=f^*\uM[d]$. We see that
$f^*A\dsv_{\G\T\L}f^*(\uM)[d]$. It follows that $A\dsv_\G\uM[d]$. This implies (a).

\subhead 4.8\endsubhead
Let $\tt=(J_n,w_n)_{n\ge0}\in\ct'(J,\boc)$. Let $d=\dim\cp_J^\tt$. Let $\cl\in\cs(T)$ and 
$\ss$ be as in 2.6. Let $\bcl$ be the local system on $\cz^\ss$ as in 2.11. Let 
$\Upss:\cz^\ss@>>>\cp_J$, $\bUp^\ss:\bZ^\ss@>>>\cp_J$ be as in 3.2. We show:

(a) {\it Let $A$ be a simple object of $\cm_\G(\cp_J)$ such that 
$A\dsv_\G{}^pH^{\cdot}(\Upss_!\bcl)$. Let $h:\cp_J^\tt@>>>\cp_J$ be the inclusion.
Let $A'$ be a simple object of $\cm_\G(\cp_J^\tt)$ such that 
$A'\dsv_\G{}^pH^{\cdot}(h^*A)$.
Then $A'\cong\uM[d]$ for some irreducible $L_{J_\iy}^{F'}$-module $M$.}
\nl
Let $K=\bUp^\ss_!\bcl^\sh$. Using 2.13(a) we
see that $\Upss_!\bcl=K$. Hence $A\dsv_\G{}^pH^{\cdot}(K)$. We have a spectral sequence in
$\cm_\G(\cp_J^\tt)$ with $E_2={}^pH^{\cdot}(h^*({}^pH^{\cdot}K))$ and $E_\iy$ is an 
associated graded of ${}^pH^{\cdot}(h^*K)$. Since $\bUp^\ss$ is proper, we see from the 
decomposition theorem \cite{\BBD} that $K\cong\op_i{}^pH^i(K)[-i]$ and that each 
${}^pH^i(K)$ is semisimple as an object of $\cm(\cp_J)$ hence also as an object of 
$\cm_\G(\cp_J)$. It also follows that $h^*K\cong\op_ih^*({}^pH^i(K))[-i]$ hence 
${}^pH^j(h^*K)\cong\op_i{}^pH^{j-i}(h^*({}^pH^i(K)))$. This shows that $E_2\cong E_\iy$ as
objects of $\cm(\cp_J^\tt)$. Thus the spectral sequence above is degenerate when regarded
in $\cm(\cp_J^\tt)$. But then it is also degenerate in $\cm_\G(\cp_J^\tt)$. Using
$A'\dsv_\G{}^pH^{\cdot}(h^*A)$ and the fact that $A$ is a direct summand of
${}^pH^{\cdot}(K)$ (in $\cm_\G(\cp_J)$) we see that $A'\dsv_\G E_2$. It follows that 
$A'\dsv_\G E_\iy$ that is $A'\dsv_\G{}^pH^{\cdot}(h^*\Upss_!\bcl)$. As in the proof of the
implication 3.2(v)$\imp$3.2(i) we deduce that 
$A'\dsv_\G{}^pH^{\cdot}(h^*\Pi^{(a)}_!\cl_{(a)})$ ($\Pi^{(a)}$ as in 3.2) for some 
$a\in W$. Now (a) follows from 4.7(a).

\subhead 4.9\endsubhead
Let $A$ be a simple object of $\cm_\G(\cp_J)$ such that 
$A\dsv_\G{}^pH^{\cdot}(\Upss_!\bcl)$ with $\cl,\ss$ as in 2.6. We show:

(a) {\it There exists $\tt=(J_n,w_n)_{n\ge0}\in\ct'(J,\boc)$ and an irreducible 
$L_{J_\iy}^{F'}$-module $M$ such that $A\cong\uM^\sh[d]$ where $\uM^\sh$ is as in 4.4 and 
$d=\dim\cp_J^\tt$.}
\nl
Since $\cp_J=\cup_{\tt\in\ct'(J,\boc)}\cp_J^\tt$, we can find $\tt\in\ct'(J,\boc)$ such that 
$\supp(A)\cap\cp_J^\tt$ is open dense in $\supp(A)$. Then, denoting by
$h:\cp_J^\tt@>>>\cp_J$ the inclusion, we see that $h^*A$ is a simple object of 
$\cm_\G(\cp_J^\tt)$. As in 4.8(a), we have $h^*A\cong\uM[d]$ for some irreducible 
$L_{J_\iy}^{F'}$-module $M$. It follows that $A$ is of the required form.

\subhead 4.10\endsubhead
Let $\tt=(J_n,w_n)_{n\ge0}\in\ct'(J,\boc)$. We set $w=w_\iy\in{}^JW$. We show:

(a) {\it for any $b\in W_{J_\iy}$, $B'\m P_{B',J}$ is a well defined map 
$\cb_{bw}@>>>\cp_J^\tt$;}

(b) {\it for $b=1$, the map $\cb_w@>>>\cp_J^\tt$ in (a) is surjective.}
\nl
We prove (a). Let $B'\in\cb_{bw}$. Let $P=P_{B',J}$. By 4.2(f) we have
$\po(P,F(P))=\min(W_JwW_{\boc(J)})=w_0$. Define $P^n$ in terms of $P$ as in 4.1. We have 
$P^1=P^{F(P)}\in\cp_{J_1}$. As in the proof of 4.6(c) (with $a_1=b,a_2=w$) we see that
$b=\po(B',P^{F(B')})$, $w=\po(P^{F(B')},F(B'))$. Since $b\in W_{J_1}$ we have
$\po(B',P^{F(B')})\in W_{J_1}$ hence $B'\sub P^1$. From the definitions we have
$w_1=\min(W_{J_1}wW_{\boc(J_1)})$ hence
$$\po(P^1,F(P^1))=\min(W_{J_1}\po(B',F(B'))W_{\boc(J_1)})=\min(W_{J_1}bwW_{\boc(J_1)})=w_1.
$$
By the same argument applied to $B',P^1,\tt_1$ instead of $B',P,\tt$ we see that 
$P^2\in\cp_{J_2}$ and $\po(P^2,F(P^2))=w_2$. (We have $w\in{}^{J_1}W$ since $J_1\sub J$.) 
Continuing in this way we see that $P^n\in\cp_{J_n}$ and $\po(P^n,F(P^n))=w_n$ for all 
$n\ge0$. Thus $P\in\cp_J^\tt$. This proves (a).

We prove (b). Let $P\in\cp_J^\tt$. Define $P^\iy$ in terms of $P$ as in 4.1. We have 
$\po(P^\iy,F(P^\iy))=w$. Hence
$$\align pr_2:&\{(B',B'')\in\cb\T\cb;B'\sub P^\iy,B''\sub F(P^\iy),\po(B',B'')=w\}\\&
@>>>\{B''\in\cb;B''\sub F(P^\iy)\}\endalign$$
is a bijection with inverse $B''\m((P^\iy)^{B''},B'')$. The condition that $(B',B'')$ in 
the domain of $pr_2$ satisfies $B''=F(B')$ is that $B''$ is a fixed point of the map
$B''\m F((P^\iy)^{B''})$ of the flag manifold of $F(P^\iy)$ into itself. This map may be 
identified with the map induced by $F':L_J@>>>L_J$ (see 4.3) on the flag manifold of $L_J$
hence it has at least one fixed point. Thus there exist $(B',B'')\in\cb\T\cb$ such that
$B''=F(B')$, $B'\sub P^\iy$, $\po(B',B'')=w$. Then $B'\in\cb_a$ and $B'\sub P$ (since
$P^\iy\sub P$). This proves (b).

\subhead 4.11\endsubhead
Let $\tt=(J_n,w_n)_{n\ge0}\in\ct'(J,\boc)$. We set $w=w_\iy\in{}^JW$. We show:

(a) {\it Let $b\in W_{J_\iy}$. Let $\cl\in\cs(T)^{bwF}$. Let 
$A$ be a simple object of $\cm_\G(\cp_J^{\tt_1})$ such that 
$A\dsv_\G{}^pH^{\cdot}(\x_{\tt_1,w\boc(b)!}(b^*\cl)_{w\boc(b)})$. Then 
$\vt^*A\dsv_\G{}^pH^{\cdot}(\x_{\tt,bw!}\cl_{bw})$.}
\nl
The result makes sense since $\cb_{\tt,bw}=\cb_{bw}$ by 4.10(a) and
$\cb_{\tt_1,w\boc(b)}=\cb_{w\boc(b)}$ (this follows from 4.10(a) applied to 
$\tt_1,w,w\boc(b)w\i$ instead of $\tt,w,b$; note that $w\boc(J_\iy)w\i=J_\iy$ by 4.2(e) 
hence 
$w\boc(b)w\i\in W_{J_\iy}$). This shows also that $l(w\boc(b))=l(w)+l(\boc(b))$. (Since 
$w\in{}^JW$ we have $l((w\boc(b)w\i)w)=l(w\boc(b)w\i)+l(w)$.)
We define $h:\cb_{bw}@>>>\cb_{w\boc(b)}$ by $B'\m B''$ with $B''$ defined by
$\po(B',B'')=b,\po(B'',F(B'))=w$. This is an isomorphism whose inverse 
$\cb_{w\boc(b)}@>>>\cb_{bw}$ is given by $B''\m B'$ with $B'$ defined by
$\po(B',B'')=b,\po(B'',F(B'))=w$. We have a commutative diagram
$$\CD 
\cb_{bw}@>h>>\cb_{w\boc(b)}\\
@V\x_{\tt,bw}VV           @V\x_{\tt_1,w\boc(b)}VV     \\
\cp_J^\tt@>\vt>>\cp_{J_1}^{\tt_1}     \endCD$$
where the horizontal maps are isomorphism. From the definitions we see that 
$h^*((b^*\cl)_{w\boc(b)})=\cl_{bw}$. The result follows.

\subhead 4.12\endsubhead
Let $\tt=(J_n,w_n)_{n\ge0}\in\ct'(J,\boc)$. We set $w=w_\iy\in{}^JW$. Let 
$d=\dim\cp_J^\tt$.
Let $\L=L_{J_\iy}^{F'}$. Let $M$ be a finite dimensional irreducible $\L$-module over 
$\bbq$. Let $\uM[d]\in\cm_\G(\cp_J^\tt)$ be as in 4.4(c). We show:

(a) {\it there exists $b\in W_{J_{\iy}}$ and $\cl\in\cs(T)^{bwF}$
such that $\uM[d]\dsv_\G{}^pH^{\cdot}(\x_{\tt,bw!}\cl_{bw})$;}

(b) {\it there exists $b\in W_{J_{\iy}}$ and $\cl\in\cs(T)^{bwF}$
such that $\uM^\sh[d]\dsv_\G{}^pH^{\cdot}(\Pi^{(bw)}_!\cl_{(bw)})$.}
\nl
We prove (a). More generally we show that for any $n\ge0$:

(c) {\it there exists $b_n\in W_{J_{\iy}}$ and $\cl\in\cs(T)^{b_nwF}$ such that $\uM[d]$ 
(regarded as an object of $\cm_\G(\cp_{J_n}^{\tt_n})$) satisfies 
$\uM[d]\dsv_\G{}^pH^{\cdot}(\x_{\tt_n,b_nw!}\cl_{b_nw})$.}
\nl
If (c) holds for $n=1$ then, by 4.11(a) it holds for $n=0$. Similarly, if (c) holds for 
some $n\ge1$ then it holds for $n-1$. Hence it is enough to prove (c) for large $n$. Thus
we may assume that $J_0=J_1=\do=J$ and $w_0=w_1=\do=w$. By \cite{\DL, 7.7} applied to $\L$,
there exists $a_1\in W_{J_\iy}$ such that $M\dsv_\L {}^pH^{\cdot}(e'_!\bbq)$ (in 
$\cm_\L(\text{ point})$). (The reference to \cite{\DL} could be replaced by a selfcontained
proof, see 7.9.) By definition, $f^*\uM[d]=e^*M[d]$ with $f$ as in 4.3, $e$ as in 4.7. We 
have $e^*M[d]\dsv_{\G\T\L}e^*({}^pH^{\cdot}(e'_!\bbq))[d]$ hence 
$f^*\uM[d]\dsv_{\G\T\L}e^*({}^pH^{\cdot}(e'_!\bbq))[d]$ 
in $\cm_{\G\T\L}(\tcp_J^\tt)$. Since
$\tcp_J^\tt$ is smooth of pure dimension $d$ we have 
$e^*({}^pH^i(e'_!\bbq))[d]={}^pH^{i+d}(e^*e'_!\bbq)$. Since 
$$e^*e'_!\bbq=p_{1!}\bbq=f^*(\x'\x)_!\bbq=f^*\x_{\tt,a!}f_{a!}\bbq$$
(notation of 4.7 with $a=a_1w$) we have 
$f^*\uM[d]\dsv_{\G\T\L}{}^pH^{\cdot}(f^*\x_{\tt,a!}f_{a!}\bbq)$.
Since $f$ is a finite principal covering we have also
$f^*\uM[d]\dsv_{\G\T\L}f^*({}^pH^{\cdot}(\x_{\tt,a!}f_{a!}\bbq))$.
It follows that $\uM[d]\dsv_\G{}^pH^{\cdot}(\x_{\tt,a!}f_{a!}\bbq)$.
Now $f_{a!}\bbq=\op_\cl\cl_a$ where $\cl$ runs over the local systems in $\cs(T)^{aF}$ 
(up to isomorphism). Hence for some such $\cl$ we have\lb
$\uM[d]\dsv_\G{}^pH^{\cdot}(\x_{\tt,a!}\cl_a)$. This proves (a).

We prove (b). Let $b,\cl$ be such that (a) holds. Let $K=\Pi^{(bw)}_!\cl_{(bw)}$,
$K'=\x_{\tt,bw!}\cl_{bw}$. We have $\uM[d]\dsv_\G{}^pH^{\cdot}(K')$. Let 
$\k:\cp_J^\tt@>>>\cp_J$ be the inclusion. We have $K=\k_!K'$. Let $C$ be the closure of 
$\cp_J^\tt$ in $\cp_J$. Let $\k':\cp_J^\tt@>>>C$, $\k'':C@>>>\cp_J$ be the inclusions. Let
$\tM=IC(C,\uM)[d]$. Let $K_1=\k'_!K$. We show:
$$\tM\dsv_\G{}^pH^{\cdot}(K_1).\tag d$$
We have $\k'{}^*K_1=K'$ and $\k'{}^*({}^pH^j(K_1))={}^pH^j(K')$. Let 
$0=X_0\sub X_1\sub\do\sub X_m={}^pH^j(K_1)$ be a composition series of ${}^pH^j(K_1)$ in 
$\cm_\G(C)$. Aplying $\k'{}^*$ to the exact sequence
$0@>>>X_{i-1}@>>>X_i@>>>X_i/X_{i-1}@>>>0$ (where $1\le i\le m$) we get an exact sequence
$0@>>>\k'{}^*(X_{i-1})@>>>\k'{}^*(X_i)@>>>\k'{}^*(X_i/X_{i-1})@>>>0$ in 
$\cm_\G(\cp_\G^\tt)$. Since $X_i/X_{i-1}$ is simple and $\cp_\G^\tt$ is open dense in $C$ 
we see that $\k'{}^*(X_i/X_{i-1})$ is either simple or $0$. It follows that any composition
factor of $\k'{}^*(X_m)$ is isomorphic to $\k'{}^*(X_i/X_{i-1})$ for some $i$. In 
particular, $\uM[d]\cong\k'{}^*(X_i/X_{i-1})$ for some $i$. It follows that
$\tM\cong X_i/X_{i-1}$ for some $i$. Thus (d) holds.

Applying $\k''_!$ to (d) we obtain 
$$\uM^\sh[d]=\k''_!\tM\dsv_\G \k''_!{}^pH^{\cdot}(K_1)={}^pH^{\cdot}(\k''_!K_1)
={}^pH^{\cdot}(K).$$
This proves (b).

\proclaim{Theorem 4.13} Let $K$ be a simple object of $\cm_\G(\cp_J)$. Then 
$K\in\Bbb S'(\cp_J)$ (see 3.7) if and only if $K\in\Bbb S(\cp_J)$.
\endproclaim
If $K\in\Bbb S(\cp_J)$ then by 4.12(b) it satisfies 3.2(i). If $K$ satisfies 3.2(v) then by
4.9(a), $K\in\Bbb S(\cp_J)$. This completes the proof.

\proclaim{Theorem 4.14} Let $\tt=(J_n,w_n)_{n\ge0}\in\ct'(J,\boc)$. Let 
$h:\cp_J^\tt@>>>\cp_J$ be the inclusion. Let $A\in\Bbb S(\cp_J)$. Let $A'$ be a simple 
object of $\cm_\G(\cp_J^\tt)$ such that $A'\dsv_\G{}^pH^{\cdot}(h^*A)$. Then 
$A'\in\Bbb S(\cp_J^\tt)$.
\endproclaim
By assumption, $A$ is as in 4.8(a). The result now follows from 4.8(a).

\head 5. Some computations in the Weyl group\endhead
\subhead 5.1\endsubhead
Let $\cl\in\cs(T)$. Define $\tl:W@>>>\NN$ by $\tl(w)=|(\a\in R^+_\cl;w(\a)\in R^-\}|$.

Let $\ss=(s_1,\do s_r)$ be a sequence in $\II\cup\{1\}$. Let
$$\ciss=\{i\in[1,r];s_i\ne1,s_1s_2\do s_i\do s_2s_1\in W_\cl\}.$$
This agrees with the definition in 2.4.

\proclaim{Lemma 5.2} We have $|\ciss|\ge\tl(s_r\do s_1)$ with equality if 
$l(s_r\do s_1)=l(s_1)+\do+l(s_r)$.
\endproclaim
Let 
$$\align&X=\{\a\in R^+_\cl;(s_r\do s_1)(\a)\in R^-\},\\&X'=\{\a\in R^+_\cl;
\a=s_1s_2\do s_{i-1}(\a_{s_i})\text{ for some }i\in[1,r]\text{ such that }s_i\ne1\}.
\endalign$$
We have $X\sub X'$. We have $|X|=\tl(s_r\do s_1)$ hence $\tl(s_r\do s_1)\le|X'|$. Define 
$f:\ciss@>>>R_\cl$ by $f(i)=s_1s_2\do s_{i-1}(\a_{s_i})$; then $X'=f(\ciss)\cap R^+_\cl$. 
Hence $|X'|\le|f(\ciss)|\le|\ciss|$ and the desired inequality is proved. Assume now that 
$l(s_r\do s_1)=l(s_1)+\do+l(s_r)$. Then $s_1s_2\do s_{i-1}(\a_{s_i}) (i\in[1,r],s_i\ne1)$ 
are distinct in $R^+$. Hence for $i\in\ciss$, $s_1s_2\do s_{i-1}(\a_{s_i})$ are distinct 
elements of $X$. Thus $|\ciss|\le|X|$. It follows that $|\ciss|=|X|$. 

\proclaim{Lemma 5.3}Let $j\in\ciss$. Let $\ss(j)=(s'_1,s'_2,\do,s'_r)$ with $s'_i=s_i$ for
$i\ne j$, $s'_j=1$. We have $\ci_{\ss(j)}=\ciss-\{j\}$.
\endproclaim  
Let $h\in\ciss-\{j\}$. We have $s_1s_2\do s_h\do s_2s_1\in W_\cl$. Hence if $j<h$, we have
$$\align&s_1s_2\do\hs_j\do s_h\do\hs_j\do s_2s_1\\&=
(s_1s_2\do s_j\do s_2s_1)(s_1s_2\do s_h\do s_2s_1)(s_1s_2\do s_j\do s_2s_1)\in W_\cl
\endalign$$
so that $h\in\ci_{\ss(j)}$. ($\hat{}$ denotes an omitted symbol.) If $j>h$ then 
$h\in\ci_{\ss(j)}$ is obvious.

Conversely, assume that $h\in\ci_{\ss(j)}$. Clearly, $h\ne j$. Assume first that $j<h$. We
have
$$\align&s_1s_2\do s_h\do s_2s_1\\&=(s_1s_2\do s_j\do s_2s_1)
(s_1s_2\do\hs_j\do s_h\do\hs_j\do s_2s_1)(s_1s_2\do s_j\do s_2s_1)\in W_\cl.\endalign$$
Hence $h\in\ciss$. If $j>h$ then it is clear that $h\in\ciss$. The lemma is proved.

\proclaim{Lemma 5.4}Let $h$ be the smallest element of $\ciss$. Then 
$s_1s_2\do s_h\do s_2s_1\in\II_\cl$.
\endproclaim
Let $\ss'=(s_1,s_2,\do,s_h,\do,s_2,s_1)$. We have

$s_1s_2\do s_h\do s_2s_1\in W_\cl$,

$s_1s_2\do s_{h-1}\do s_2s_1\n W_\cl$ or $s_{h-1}=1$,

$s_1s_2\do s_{h-2}\do s_2s_1\n W_\cl$ or $s_{h-2}=1$,

$\do$.
\nl
Hence the middle term in $\ss'$ has index in $\ci_{\ss'}$ but all terms preceding it have 
index not in $\ci_{\ss'}$. We show that the term in $\ss'$ immediately following the middle
term has index not in $\ci_{\ss'}$. Otherwise it would be $\ne1$ and 
$$s_1s_2\do s_{h-1}s_hs_{h-1}s_hs_{h-1}\do s_2s_1\n W_\cl.$$
Multiplying on the left and right by $s_1s_2\do s_h\do s_2s_1$ we find 
$$s_1s_2\do s_{h-1}\do s_2s_1\in W_\cl,$$
a contradiction. Similarly we see that all terms in $\ss'$ following the middle term have 
index not in $\ci_{\ss'}$. Thus $|\ci_{\ss'}|=1$. By Lemma 5.2 we have 
$$\tl(s_1s_2\do s_h\do s_2s_1)\le1.$$
Since $s_1s_2\do s_h\do s_2s_1\in W_\cl-\{1\}$, it follows that 
$s_1s_2\do s_h\do s_2s_1\in\II_\cl$. The lemma is proved.

\subhead 5.5\endsubhead
We write the elements of $\ciss$ in ascending order: $i_1<i_2<\do<i_b$. Define a sequence 
$\SS=(S_1,S_2,\do,S_b)$ in $W$ by

$S_1=s_1s_2\do s_{i_1}\do s_2s_1$,

$S_2=s_1s_2\do\hs_{i_1}\do s_{i_2}\do\hs_{i_1}\do s_2s_1$,

$\do$,

$S_b=s_1s_2\do\hs_{i_1}\do\hs_{i_2}\do\hs_{i_{b-1}}\do s_{i_b}
\do\hs_{i_{b-1}}\do\hs_{i_2}\do\hs_{i_1}\do s_2s_1$,

$\o=s_1\do\hs_{i_1}\do\hs_{i_2}\do\hs_{i_b}\do s_r$.

\proclaim{Lemma 5.6}(a) $S_1,S_2,\do,S_b$ belong to $\II_\cl$.

(b) $\tl(\o\i)=0$.
        
(c) $s_1s_2\do s_r=S_1S_2\do S_b\o$. 
\endproclaim
We use induction on $b$. Assume first that $b=0$. By 5.2 we have $\tl(s_r\do s_1)=0$ that 
is $\tl(\o\i)=0$ and the lemma is clear. Assume now that $b\ge1$ and that the lemma holds 
for $b-1$. Let $h=i_1$. By 5.3 we have $\ci_{\ss(h)}=\ciss-\{h\}$. By the induction 
hypothesis, $S_2,S_3,\do,S_b$ belong to $\II_\cl$, 
$s_1\do\hs_{i_1}\do s_r=S_2S_3\do S_b\o$, $\tl(\o\i)=0$. By 5.4 we have $S_1\in\II_\cl$. It
follows that 
$$s_1s_2\do s_r=(s_1\do s_{i_1}\do s_1)(s_1\do\hs_{i_1}\do s_r)=S_1S_2\do S_b\o.$$
The lemma is proved.

\subhead 5.7\endsubhead
Let $W'_\cl=\{w\in W;w^*\cl\cong\cl\}$, a subgroup of $W$. Let 
$W'{}^0_\cl=\{w\in W'_\cl;w(R^+_\cl)=R^+_\cl\}$, a subgroup of $W'_\cl$. Note that
$W'_\cl=W'{}^0_\cl W_\cl$ (semidirect product with $W_\cl$ normal).

In the remainder of this section we fix an automorphism $c$ of finite order of $T$ that 
induces a permutation of $R$, one of $\chR$ and one of $R^+$. 

\proclaim{Lemma 5.8} Assume that $F_0^*(s_1s_2\do s_rc)^*\cl\cong\cl$. Let $c'=\o c$ so 
that $s_1s_2\do s_rc=S_1S_2\do S_bc'$ and $F_0^*(S_1S_2\do S_bc')^*\cl\cong\cl$. Then $c'$
is an automorphism of finite order of $T$ that induces a permutation of $R_\cl$, one of 
$\chR_\cl$ and one of $R^+_\cl$. 
\endproclaim
We first show:
$$c'(R_\cl)=R_\cl.\tag i$$
It is enough to show that for any coroot $\k:\kk^*@>>>T$ such that $\k^*(\cl)\cong\bbq$ we
have $(c'\k)^*\cl\cong\bbq$. From $F_0^*(S_1S_2\do S_bc')^*\cl\cong\cl$ and 
$S_1S_2\do S_r\in W_\cl$ we deduce $F_0^*c'{}^*\cl\cong\cl$. Hence 
$\k^*F_0^*c'{}^*\cl\cong\k^*\cl\cong\bbq$. Now $F_0\k=\k F'_0$ where $F'_0:\kk^*@>>>\kk^*$
is $x\m x^q$. Hence $F'_0{}^*\k^*c'{}^*\cl=\bbq=F'_0{}^*\bbq$ and $\k^*c'{}^*\cl=\bbq$. 
This proves (i).

By 5.6(b) we have $\o\i(R^+_\cl)\sub R^+$. Moreover, $c\i(R^+)=R^+$ hence \lb 
$c\i\o\i(R^+_\cl)\sub R^+$. Using (i) we see that $c'{}\i(R^+_\cl)\sub R^+_\cl$ hence
$c'(R^+_\cl)=R^+_\cl$. The lemma is proved.

\subhead 5.9\endsubhead
Let $\tss=(\ts_1,\ts_2,\do,\ts_{\y})$ be a second sequence in $\II$.
Then the subset $\cits$ is defined in terms of $\tss,\cl$ in the same way as $\ciss$ is 
defined in terms of $\ss,\cl$.

We write the elements of $\cits$ in ascending order: $j_1<j_2<\do<j_{\tb}$. Define 
$\tSS=(\tS_1,\tS_2,\do,\tS_{\tb}),\tio$ in terms of $\tss,\cl$ in the same way that
$\SS,\o$ are defined in 5.5 in terms of $\ss,\cl$. As in 5.6 we have $\tS_i\in\II_\cl$ for
$i\in[1,\tb]$ and $\ts_1\ts_2\do\ts_{\y}=\tS_1\tS_2\do\tS_{\tb}\tio$. We assume that 
$$F_0^*([\ss]c)^*\cl\cong\cl, F_0^*([\tss]c)^*\cl\cong\cl.$$
From 5.8 we see that $\o c$ and $\tio c$ are automorphisms of finite order of $T$ that 
induce a permutation of $R_\cl$, one of $\chR_\cl$ and one of $R^+_\cl$. Moreover, we have
$$F_0^*(S_1S_2\do S_b\o c)^*\cl\cong\cl,F_0^*(\tS_1\tS_2\do\tS_{\tb}\tio c)^*\cl\cong\cl,$$
or equivalently $F_0^*(\o c)^*\cl\cong\cl,F_0^*(\tio c)^*\cl\cong\cl$. Let 
$$\fF=\{f\in W'{}^0_\cl;f\i\o cf=\tio c\}.$$
Now for $f\in\fF$ we have $f\II_\cl f\i=\II_\cl$; we set 
$${}^f\tSS=(f\tS_1f\i,f\tS_2f\i,\do f\tS_{\tb}f\i),$$
a sequence in $\II_\cl$. Since $F_0^*(\o c)^*\cl\cong\cl$, we have
$$F_0^*((f\tS_1f\i)(f\tS_2f\i)\do(f\tS_{\tb}f\i)\o c)^*\cl\cong\cl.$$
We set $\r=r+\y$. Let $\ca(W,c,\cl,\ss,\tss)$ be the set of all sequences 
$(a_0,a_1,\do,a_\r)$ in $W$ such that 

$a_{j-1}\i a_j\in\{1,\ts_j\}$ for $j\in\cits$;

$a_{j-1}\i a_j=\ts_j$ for $j\in[1,\y]-\cits$;

$a_{\y+i}a_{\y+i-1}\i\in\{1,s_i\}$ for $i\in\ciss$;

$a_{\y+i}a_{\y+i-1}\i=s_i$ for $i\in[1,r]-\ciss$;

$a_\r=ca_0c\i$;    

$a_0^*\cl\cong\cl$.
\nl
Replacing here $W,c,\cl,\ss,\tss$ by $W_\cl,\o c,\bbq,\SS,{}^f\tSS$, we obtain for any
$f\in\fF$ a set $\ca(W_\cl,\o c,\bbq,\SS,{}^f\tSS)$. From the definition,
$\ca(W_\cl,\o c,\bbq,\SS,{}^f\tSS)$ consists of all sequences $(A_0,A_1,\do,A_{b+\tb})$ 
in $W_\cl$ such that 

$A_{j-1}\i A_j\in\{1,f\tS_jf\i\}$ for $j\in[1,\tb]$;

$A_{\tb+i}A_{\tb+i-1}\i\in\{1,S_i\}$ for $i\in[1,b]$;

$A_{b+\tb}=(\o c)A_0(\o c)\i$.
\nl
We now state the following result. 

\proclaim{Proposition 5.10} There is a canonical bijection
$$\Psi:\ca(W,c,\cl,\ss,\tss)@>\si>>\sqc_{f\in\fF}\ca(W_\cl,\o c,\bbq,\SS,{}^f\tSS).$$
\endproclaim
Let $(a_0,a_1,\do,a_\r)\in\ca(W,c,\cl,\ss,\tss)$. Consider the product
$$\align a_{\y}&(a_{\y+i_1}\i a_{\y+i_1-1})a_{\y}\i
=(a_{\y}a_{\y+1}\i)(a_{\y+1}\a_{\y+2}\i)\do(a_{\y+i_1-2}a_{\y+i_1-1}\i)\\&
\T(a_{\y+i_1-1}a_{\y+i_1}\i)(a_{\y+i_1-1}a_{\y+i_1-2}\i)\do
(a_{\y+2}a_{\y+1}\i)(a_{\y+1}a_{\y}\i).\endalign$$
The right hand side is equal to $s_1s_2\do s_{i_1-1}xs_{i_1-1}\do s_2s_1$ where $x$ is 
either 
$s_{i_1}$ or $1$. Hence it is equal to $s_1s_2\do s_{i_1-1}s_{i_1}s_{i_1-1}\do s_2s_1$ or 
to $1$. Thus we have
$$a_{\y}(a_{\y+i_1}\i a_{\y+i_1-1})a_{\y}\i\in\{S_1,1\}.$$
Similarly, 
$$\align a_{\y}&(a_{\y+i_1-1}\i a_{\y+i_1})(a_{\y+i_2-1}\i a_{\y+i_2})\do
(a_{\y+i_{e-1}-1}\i a_{\y+i_{e-1}})(a_{\y+i_e}\i a_{\y+i_e-1})\\&\T
(a_{\y+i_{e-1}}\i a_{\y+i_{e-1}-1})
\do(a_{\y+i_2}\i a_{\y+i_2-1})(a_{\y+i_1}\i a_{\y+i_1-1})
a_{\y}\i\in\{S_e,1\}\endalign$$
for $e\in[1,b]$,
$$\align&a_0\i(a_{j_1-1}a_{j_1}\i)(a_{j_2-1}a_{j_2}\i)\do(a_{j_{e-1}-1}a_{j_{e-1}}\i)
(a_{j_e}a_{j_e-1}\i)\\&(a_{j_{e-1}}a_{j_{e-1}-1}\i)\do(a_{j_2}a_{j_2-1}\i)
(a_{j_1}a_{j_1-1}\i)\in\{\tS_e,1\}\endalign$$
for $e\in[1,\tb]$,
$$\o=a_{\y}(a_{\y+i_1-1}\i a_{\y+i_1})(a_{\y+i_2-1}\i a_{\y+i_2})\do
(a_{\y+i_b-1}\i a_{\y+i_b})a_\r\i,$$
$$\tio=a_0\i(a_{j_1-1}a_{j_1}\i)(a_{j_2-1}a_{j_2}\i)\do(a_{j_{\tb-1}}a_{j_{\tb}}\i)a_{\y}.
$$
Setting 
$$\ha_e=a_{j_e}a_{j_e-1}\i\text{ for }e\in[1,\tb],\qua
\ha_{\tb+e}=a_{\y}a_{\y+i_e}\i a_{\y+i_e-1}a_{\y}\i\text{ for }e\in[1,b]\tag a$$
we see that
$$\ha_{\tb+1}\i\ha_{\tb+2}\i\do\ha_{\tb+e-1}\i\ha_{\tb+e}\ha_{\tb+e-1}\do\ha_{\tb+2}
\ha_{\tb+1}\in\{S_e,1\}\text{ for }e\in[1,b],\tag b$$
$$a_0\i\ha_1\i\ha_2\i\do\ha_{e-1}\i\ha_e\ha_{e-1}\do\ha_2\ha_1a_0\in\{\tS_e,1\}\text{ for }
e\in[1,\tb],\tag b$$
$$\o=\ha_{\tb+1}\i\ha_{\tb+2}\i\do\ha_{\tb+b}\i a_{\y}a_\r\i,\tag c$$
$$\tio=a_0\i\ha_1\i\ha_2\i\do\ha_{\tb}\i a_{\y}.\tag d$$
Since $S_e=S_e\i\in W_\cl$, $\tS_e=\tS_e\i\in W_\cl$ it follows by induction on $e$ that
$$\ha_{\tb+e}=\ha_{\tb+e}\i\in W_\cl \text{ for }e\in[1,b],\tag e$$
$$a_0\i\ha_ea_0=a_0\i\ha_e\i a_0\in W_\cl\text{ for }e\in[1,\tb].$$
Since $a_0$ normalizes $W_\cl$ it follows that
$$\ha_e=\ha_e\i\in W_\cl\text{ for }e\in[1,\tb].\tag f$$
Since $a_0^*\cl\cong\cl$ we can write uniquely 
$$a_0=A_0f\text{ with }A_0\in W_\cl,f\in W'{}^0_\cl.\tag g$$ 
We set
$$A_e=\ha_e\do\ha_2\ha_1a_0f\i\qua(e\in[1,\tb]),\tag h$$
$$A_{\tb+e}=\ha_{\tb+1}\ha_{\tb+2}\do\ha_{\tb+e}\ha_{\tb}\do\ha_2\ha_1a_0f\i\qua
(e\in[1,b]).\tag i$$
From (e),(f),(g) we see that $A_e\in W_\cl$ for $e\in[0,b+\tb]$. From the definitions we 
have
$$A_{e-1}\i A_e=fS_ef\i\qua(e\in[1,\tb]),$$
$$A_{\tb+e-1}A_{\tb+e}\i=S_e\qua(e\in[1,b]).$$
We have
$$\align
A_{\tb+b}&=\ha_{\tb+1}\ha_{\tb+2}\do\ha_{\tb+b}\ha_{\tb}\do\ha_2\ha_1a_0f\i=\o a_\r\tio\i
f\i\\&=(\o c)(c\i a_\r c)(\tio c)\i f\i=(\o c)a_0(\tio c)\i f\i=(\o c)A_0f(\tio c)\i f\i.
\endalign$$
In particular, $(\o c)wf(\tio c)\i f\i=w'$ for some $w\in W_\cl$. We have 
$(\o c)w=w_1(\o c)$ for some $w_1\in W_\cl$ hence $(\o c)f(\tio c)\i f\i=w_1\i w'$ belongs
to $W_\cl\cap W'{}^0_\cl=\{1\}$. Thus $(\o c)f(\tio c)\i f\i=1$ that is $f\in\fF$. We also
see that 
$$A_{\tb+b}=(\o c)A_0f(\tio c)\i f\i=(\o c)A_0(\o c)\i.$$
Thus to each $(a_0,a_1,\do,a_\r)\in\ca(W,c,\cl,\ss,\tss)$ we have associated $f\in\fF$ and
an element $(A_0,A_1,\do,A_{b+\tb})$ in $\ca(W_\cl,\o c,\bbq,\SS,{}^f\tSS)$. This defines 
the map $\Psi$ in the proposition.

Conversely, assume that we are given $f\in\fF$ and an element $(A_0,A_1,\do,A_{b+\tb})$ in
$\ca(W_\cl,\o c,\bbq,\SS,{}^f\tSS)$. We will construct a sequence $(a_0,a_1,\do,a_\r)$ in
$W$ as follows.

We set $a_0=A_0f$. We define $\ha_1,\ha_2,\do,\ha_{\tb}$ inductively so that (h) holds. We
define $\ha_{\tb+1},\ha_{\tb+2},\do,\ha_{\tb+b}$ inductively so that (i) holds. We define 
the elements $a_0,a_1,a_2,\do,a_{\y}$ by induction as follows. Note that $a_0$ is already 
defined. Assume that $a_0,a_1,\do,a_{u-1}$ are already defined for some $u\in[1,\y]$. If 
$u=j_e$ for some $e\in[1,\tb]$ we set $a_u=\ha_ea_{u-1}$. If
$u\notin\cits$ we set $a_u=a_{u-1}\ts_u$. This completes the definition of
$a_0,a_1,a_2,\do,a_{\y}$.
We define the elements $a_{\y},a_{\y+1},\do,a_{\y+r}$ by induction as follows. Note that 
$a_{\y}$ is already defined. Assume that $a_{\y},a_{\y+1},\do,a_{\y+u-1}$ are already 
defined for some $u\in[1,r]$. If $u=i_e$ for some $e\in[1,b]$ we set
$a_{\y+u}=a_{\y+u-1}a_{\y}\i\ha_{\tb+e}a_{\y}$. If $u\notin\ciss$ we set
$a_{\y+u}=s_ua_{\y+u-1}$. This completes the definition of $a_{\y},a_{\y+1},\do,a_{\y+r}$.

From the definitions we see that $(a_0,a_1,a_2,\do,a_\r)\in\ca(W,c,\cl,\ss,\tss)$. We have
thus constructed a map 
$$\sqc_{f\in\fF}\ca(W_\cl,\o c,\bbq,\SS,{}^f\tSS)@>>>\ca(W,c,\cl,\ss,\tss).$$
From the definitions we see that this is an inverse to $\Psi$. The proposition is proved.

\head 6. A basis for a space of intertwining operators\endhead
\subhead 6.1\endsubhead
Let $\cl,\tcl\in\cs(T)$. Let $\ss=(s_1,s_2,\do,s_r)$, $\tss=(\ts_1,\ts_2,\do,\ts_{\y})$ be
two sequences in $\II$ such that $([\ss]F)^*\cl\cong\cl,([\tss]F)^*\tcl\cong\tcl$. Define 
$\ciss\sub[1,r]$ in terms of $\ss,\cl$ as in 2.4 or 5.1; define $\cits\sub[1,\y]$ similarly
in terms of $\tss,\tcl$. 

We set $\r=r+\y$. Until the end of 6.13 we fix a sequence $\aa=(a_0,a_1,\do,a_\r)$ in $W$
such that 

$a_{j-1}\i a_j\in\{1,\ts_j\}$ for $j\in[1,\y]$;

$a_{\y+i}a_{\y+i-1}\i\in\{1,s_i\}$ for $i\in[1,r]$;

if $j\in[1,\y]-\cits$, $a_{j-1}\ts_j>a_{j-1}$ then $a_{j-1}\i a_j=\ts_j$;

if $i\in[1,r]-\ciss$, $s_ia_{\y+i-1}>a_{\y+i-1}$ then $a_{\y+i}a_{\y+i-1}\i=s_i$;

$a_\r=\boc(a_0)$.
\nl
We define representatives $\dda_i$ for $a_i$ in $N(T)$ ($i\in[0,\r]$) as follows. We set 
$\dda_0=\doa_0$. For $j\in[1,\y]$ we define $\dda_j$ inductively by

$\dda_j=\dda_{j-1}$ if $a_j=a_{j-1}$,

$\dda_j=\dda_{j-1}\dts_j$ if $a_j=a_{j-1}\ts_j$.
\nl
For $i\in[1,r]$ we define $\dda_{\y+i}$ inductively by

$\dda_{\y+i}=\dda_{\y+i-1}$ if $a_{\y+i-1}=a_{\y+i}$,

$\dda_{\y+i}=\ds_i\i\dda_{\y+i-1}$ if $a_{\y+i}=s_ia_{\y+i-1}$.
\nl
Consider the commutative diagrams
$$\CD     
Z\T\tZ@<d_0\T\td_0<<Z_0\T\tZ_0@<d_1\T\td_1<<Z_1\T\tZ_1@>d_2\T\td_2>>Z_2\T\tZ_2@>=>>
Z_2\T\tZ_2\\
@AeAA @Ae_0AA @Ae_1AA @Ae_2AA @Ae_3AA \\
X@<b_0<<X_0@<b_1<<X_1@>b_2>>X_2@>b_3>>X_3\endCD$$
$$\CD     
Z_2\T\tZ_2@>=>>Z_2\T\tZ_2@>=>>Z_2\T\tZ_2@<=<<Z_2\T\tZ_2@. X_8\\
@Ae_3AA @Ae_4AA  @Ae_5AA @Ae_6AA  @Ab_8AA\\
X_3@<b_4<<X_4@>b_5>>X_5@<b_6<<X_6@>b_7>>X_7\endCD$$
where the following notation is used.

$Z,Z_i,d_i$ are as in 2.6, $\tZ,\tZ_i,\td_i$ are the analogous objects defined in terms of 
$\tss,\tcl$ instead of $\ss,\cl$. 

$X$ is the set of all 
$((B_0,B_1,\do,B_r),(\tB_0,\tB_1,\do,\tB_{\y}))\in\cb^{[0,r]}\T\cb^{[0,\y]}$
such that $(A0)-(A5)$ below hold:

($A0$) $\po(B_0,\tB_j)=a_j(j\in[0,\y])$, $\po(B_i,\tB_{\y})=a_{\y+i}(i\in[0,r])$,

($A1$) $\po(B_{i-1},B_i)=s_i$ if $i\in[1,r]-\ciss$,

($A2$) $\po(B_{i-1},B_i)\in\{1,s_i\}$ if $i\in\ciss$,

($A3$) $\po(\tB_{j-1},\tB_j)=\ts_j$ if $j\in[1,\y]-\cits$,

($A4$) $\po(\tB_{j-1},\tB_j)\in\{1,\ts_j\}$ if $j\in\cits$,

($A5$) $F(B_0)=B_r, F(\tB_0)=\tB_{\y}$.
\nl
$X_0$ is the set of all 
$$((g_0U,g_1U,\do,g_rU),(\tg_0U,\tg_1U,\do,\tg_{\y}U))\in(G/U)^{[0,r]}\T(G/U)^{[0,\y]}$$ 
such that $(B0)-(B5)$ below hold:

($B0$) $k(g_0\i\tg_j)=\dda_j(j\in[0,\y]),k(g_i\i\tg_{\y})=\dda_{\y+i}(i\in[0,r])$,

($B1$) $g_{i-1}\i g_i\in P_{s_i}-B$ if $i\in[1,r]-\ciss$,

($B2$) $g_{i-1}\i g_i\in P_{s_i}$ if $i\in\ciss$,

($B3$) $\tg_{j-1}\i\tg_j\in P_{\ts_j}-B$ if $j\in[1,\y]-\cits$,

($B4$) $\tg_{j-1}\i\tg_j\in P_{\ts_j}$ if $j\in\cits$,

($B5$) $g_r\i F(g_0)\in U$, $\tg_{\y}\i F(\tg_0)\in U$.
\nl
$b_0$ is 
$$\align&((g_0U,g_1U,\do,g_rU),(\tg_0U,\tg_1U,\do,\tg_{\y}U))\m \\&
((g_0Bg_0\i,g_1Bg_1\i,\do,g_rBg_r\i),
(\tg_0B\tg_0\i,\tg_1B\tg_1\i,\do,\tg_{\y}B\tg_{\y}\i)).\endalign$$
$X_1$ is the set of all 
$((g_0,g_1,\do,g_r),(\tg_0,\tg_1,\do,\tg_{\y}))\in G^{[0,r]}\T G^{[0,\y]}$
such that $(B0)-(B5)$ hold. 
\nl
$b_1$ is  
$((g_0,g_1,\do,g_r),(\tg_0,\tg_1,\do,\tg_{\y}))\m 
((g_0U,g_1U,\do,g_rU)),(\tg_0U,\tg_1U,\do,\tg_{\y}U))$.
\nl
$X_2$ is the set of all $(x,x',u,u',y_0,y_1,\do,y_\r)\in G\T G\T U\T U\T G^{[0,\r]}$ such 
that $(C0)-(C'5)$ below hold:

$(C0$) $k(y_z)=\dda_z (z\in[0,\r])$    

$(C1)$ $y_{\y+i-1}y_{\y+i}\i\in P_{s_i}-B$ if $i\in[1,r]-\ciss$,

$(C2)$ $y_{\y+i-1}y_{\y+i}\i\in P_{s_i}$ if $i\in\ciss$,

$(C3)$ $y_{j-1}\i y_j\in P_{\ts_j}-B$ if $j\in[1,\y]-\cits$,

$(C4)$ $y_{j-1}\i y_j\in P_{\ts_j}$ if $j\in\cits$,

$(C5)$ $uF(y_0)=y_\r u'$,

$(C'5)$ $y_0x'=xF(y_0)$, $y_\r y_{\y}\i x=u$. 
\nl
$b_2$ is $((g_0,g_1,\do,g_r),(\tg_0,\tg_1,\do,\tg_{\y}))\m(x,x',u,u',y_0,y_1,\do,y_\r)$
\nl
where
$$\align&x=g_0\i F(g_0),x'=\tg_0\i F(\tg_0),u=g_r\i F(g_0),u'=\tg_{\y}\i F(\tg_0),\\&
y_j=g_0\i\tg_j(j\in[0,\y]),y_{\y+i}=g_i\i\tg_{\y}(i\in[0,r]).\endalign$$
$X_3$ is the set of all $(u,u',y_0,y_1,\do,y_\r)\in U^2\T G^{[0,\r]}$ such that 
$(C0)-(C5)$ hold. 

$b_3$ is $(x,x',u,u',y_0,y_1,\do,y_\r)\m(u,u',y_0,y_1,\do,y_\r)$.

$X_4$ is the set of all $(u,u',v,v',\tv,\tv',y_0,y_1,\do,y_\r)\in U^6\T G^{[0,\r]}$ such 
that $(C0)-(C4)$ hold and 
$$y_0=v\dda_0v',y_\r=\tv F(\dda_0)\tv',uF(v)F(\dda_0)F(v')=\tv F(\dda_0)\tv'u'.$$
$b_4$ is $(u,u',v,v',\tv,\tv',y_0,y_1,\do,y_\r)\m(u,u',y_0,y_1,\do,y_\r)$.

$X_5$ is the set of all $(v,v',y_0,y_1,\do,y_\r)\in U^2\T G^{[0,\r]}$ such that $(C0)-(C4)$
hold and $y_0=v\dda_0v'$.

$b_5$ is $(u,u',v,v',\tv,\tv',y_0,y_1,\do,y_\r)\m(v,v',y_0,y_1,\do,y_\r)$.

$X_6$ is the set of all 
$((g_0,g_1,\do,g_r),(\tg_0,\tg_1,\do,\tg_{\y}))\in G^{[0,r]}\T G^{[0,\y]}$
such that $(B0)-(B4)$ hold and $g_0\in U,\tg_0\in\dda_0U$.

$b_6$ is $((g_0,g_1,\do,g_r),(\tg_0,\tg_1,\do,\tg_{\y}))\m(v,v',y_0,y_1,\do,y_\r)$
\nl
where 
$$v=g_0\i,v'=\dda_0\i\tg_0,y_j=g_0\i\tg_j(j\in[0,\y]),y_{\y+i}=g_i\i\tg_{\y}(i\in[0,r]).$$
$X_7$ is the set of all 
$$((g_0U,g_1U,\do,g_rU),(\tg_0U,\tg_1U,\do,\tg_{\y}U))\in(G/U)^{[0,r]}\T(G/U)^{[0,\y]}$$
such that $(B0)-(B4)$ hold and $g_0\in U,\tg_0\in\dda_0U$.

$b_7$ is given by the same formula as $b_1$.

$X_8$ is the set of all 
$((B_0,B_1,\do,B_r),(\tB_0,\tB_1,\do,\tB_{\y}))\in\cb^{[0,r]}\T\cb^{[0,\y]}$	
such that $(A0)-(A4)$ hold and $B_0=B,\tB_0=\dda_0B\dda_0\i$.

$b_8$ is given by the same formula as $b_0$.

$e,e_0$ are the obvious imbeddings. 

$e_1$ is
$$((g_0,g_1,\do,g_r),(\tg_0,\tg_1,\tg_{\y}))\m
((g_0,g_0g_1\i,\do,g_{r-1}\i g_r),(\tg_0,\tg_0\tg_1\i,\do,\tg_{\y-1}\i\tg_{\y})).$$
$e_2$ is $(x,x',u,u',y_0,y_1,\do,y_\r)\m\l(y_0,y_1,\do,y_\r)$ where 
$$\l(y_0,y_1,\do,y_\r)=((y_{\y}y_{\y+1}\i,y_{\y+1}y_{\y+2}\i,\do,y_{\y+r-1}y_{\y+r}\i),
(y_0\i y_1,y_1\i y_2,\do,y_{\y-1}\i y_{\y})).$$
$e_3$ is $(u,u',y_0,y_1,\do,y_\r)\m\l(y_0,y_1,\do,y_\r)$.

$e_4$ is $(u,u',v,v',\tv,\tv',y_0,y_1,\do,y_\r)\m\l(y_0,y_1,\do,y_\r)$.

$e_5$ is $(v,v',y_0,y_1,\do,y_\r)\m\l(y_0,y_1,\do,y_\r)$.

$e_6$ is
$$\align&((g_0,g_1,\do,g_r),(\tg_0,\tg_1,\do,\tg_{\y}))\m\\&
((g_0\i g_1,g_1\i g_2,\do,g_{r-1}\i g_r),(\tg_0\i\tg_1,\tg_1\i\tg_2,\do,\tg_{r-1}\i\tg_r)).
\endalign$$

\subhead 6.2\endsubhead
Let $\uB,\bcl,\bcl_i$ be as in 2.4, 2.11. Define $\utB,\tbcl,\tbcl_i$ similarly, in terms 
of $\tss,\tcl$ instead of $\ss,\cl$. Note that $\G,\uB$ act on $Z$ and $Z_i(i\in[0,2])$ as
in 2.7. Similarly, $\G,\utB$ act on $\tZ$ and $\tZ_i(i\in[0,2])$. These actions give rise
to (commuting) actions of $\G\T\G,\uB\T\utB$ on $Z\T\tZ$ and $Z_i\T\tZ_i(i\in[0,2])$ hence
to actions of $\G\T\G\T\uB\T\utB$. By 2.11, $\bcl$ and $\bcl_i$ are $\uB$-equivariant local
systems with natural $\G$-equivariant structures. Similarly, $\tbcl$ and $\tbcl_i$ are 
$\utB$-equivariant local systems with natural $\G$-equivariant structures. Hence 
$\bcl\bxt\tbcl$ and $\bcl_i\bxt\tbcl_i$ are $\uB\T\utB$-equivariant local systems on 
$Z\T\tZ$ and $Z_i\T\tZ_i$ with natural $\G\T\G$-equivariant structures. These structures 
give rise to $\G\T\G\T\uB\T\utB$-equivariant structures on $\bcl\bxt\tbcl$ and 
$\bcl_i\bxt\tbcl_i$. From the definitions we have 
$$\align&(d_0\T\td_0)^*(\bcl\bxt\tbcl)=\bcl_0\bxt\tbcl_0,
(d_1\T\td_1)^*(\bcl_0\bxt\tbcl_0)=\bcl_1\bxt\tbcl_1,\\&
(d_2\T\td_2)^*(\bcl_2\bxt\tbcl_2)=\bcl_1\bxt\tbcl_1,\endalign$$
compatibly with the $\G\T\G\T\uB\T\utB$-equivariant structures.

\subhead 6.3\endsubhead
Let 
$$\align\ct&=\{((t_0,t_1,\do,t_r),(\tit_0,\tit_1,\do,\ttir))\in T^{[0,r]}\T T^{[0,\y]};\\&
t_r=F(t_0),\ttir=F(\tit_0),\tit_j=a_j\i(t_0)\text{ for }j\in[0,\y], 
t_i=a_{\y+i}(\ttir)\text{ for }i\in[0,r].\endalign$$
This is a subgroup of $T^{[0,r]}\T T^{[0,\y]}$ isomorphic to the finite subgroup
$$\ct_0=\{t_0\in T;a_{\y}\i(t_0)=F(a_0\i(t_0))\}$$ 
of $T$ under $((t_i),(\tit_j))\m t_0$. Hence $\ct$ is finite.

$\G\T\ct$ acts on $X$ by
$$\align&(g,(t_i),(\tit_j)):((B_0,B_1,\do,B_r),(\tB_0,\tB_1,\do,\tB_{\y}))\m\\&
((gB_0g\i,gB_1g\i,\do,gB_rg\i),(g\tB_0g\i,g\tB_1g\i,\do,g\tB_{\y}g\i));\endalign$$
on $X_0$ by   
$$\align&(g,(t_i),(\tit_j)):((g_0U,g_1U,\do,g_rU),(\tg_0U,\tg_1U,\do,\tg_{\y}U))\m\\&
((gg_0t_0\i U,gg_1t_1\i U,\do,gg_rt_r\i U),(g\tg_0\tit_0\i U,g\tg_1\tit_1\i U,\do,
g\tg_{\y}\ttir\i U));\endalign$$
on $X_1$ by
$$\align&(g,(t_i),(\tit_j)):((g_0,g_1,\do,g_r),(\tg_0,\tg_1,\do,\tg_{\y}))\m\\&
((gg_0t_0\i,gg_1t_1\i,\do,gg_rt_r\i),(g\tg_0\tit_0\i,g\tg_1\tit_1\i,\do,g\tg_{\y}\ttir\i));
\endalign$$
on $X_2$ by
$$\align&(g,(t_i),(\tit_j)):(x,x',u,u',y_0,y_1,\do,y_\r)\m
(t_0xF(t_0)\i,\tit_0x'F(\tit_0)\i,t_rut_r\i,\ttir u'\ttir\i, \\&
t_0y_0\tit_0\i,t_0y_1\tit_1\i,\do,t_0y_{\y}\ttir\i,t_1y_{\y+1}\ttir\i,\do,t_ry_\r\ttir\i);
\endalign$$
on $X_3$ by
$$\align&(g,(t_i),(\tit_j)):(u,u',y_0,y_1,\do,y_\r)\m\\&
(t_rut_r\i,\ttir u'\ttir\i,t_0y_0\tit_0\i,t_0y_1\tit_1\i,\do,t_0y_{\y}\ttir\i,
t_1y_{\y+1}\ttir\i,\do,t_ry_\r\ttir\i);\endalign$$
on $X_4$ by
$$\align&(g,(t_i),(\tit_j)):(u,u',v,v',\tv,\tv',y_0,y_1,\do,y_\r)\m\\&
(t_rut_r\i,\ttir u'\ttir\i,t_0vt_0\i,\tit_0v'\tit_0\i,t_0y_0\tit_0\i,t_0y_1\tit_1\i,\do,
t_0y_{\y}\ttir\i,t_1y_{\y+1}\ttir\i,\do,t_ry_\r\ttir\i);\endalign$$
on $X_5$ by
$$\align&(g,(t_i),(\tit_j)):(v,v',\tv,\tv',y_0,y_1,\do,y_\r)\m\\&
(t_0vt_0\i,\tit_0v'\tit_0\i,t_0y_0\tit_0\i,t_0y_1\tit_1\i,\do,t_0y_{\y}\ttir\i,
t_1y_{\y+1}\ttir\i,\do,t_ry_\r\ttir\i);\endalign$$
on $X_6$ by 
$$\align&(g,(t_i),(\tit_j)):((g_0,g_1,\do,g_r),(\tg_0,\tg_1,\do,\tg_{\y}))\m\\&
((t_0g_0t_0\i,t_0g_1t_1\i,\do,t_0g_rt_r\i),
(t_0\tg_0\tit_0\i,t_0\tg_1\tit_1\i,\do,t_0\tg_{\y}\ttir\i));\endalign$$
on $X_7$ by 
$$\align&(g,(t_i),(\tit_j)):((g_0U,g_1U,\do,g_rU),(\tg_0U,\tg_1U,\do,\tg_{\y}U))\m\\&
((t_0g_0t_0\i U,t_0g_1t_1\i U,\do,t_0g_rt_r\i U),
(t_0\tg_0\tit_0\i U,t_0\tg_1\tit_1\i U,\do,t_0\tg_{\y}\ttir\i U));\endalign$$
on $X_8$ by
$$\align&(g,(t_i),(\tit_j)):((B_0,B_1,\do,B_r),(\tB_0,\tB_1,\do,\tB_{\y}))\m\\&
((t_0B_0t_0\i,t_0B_1t_0\i,\do,t_0B_rt_0\i),
(t_0\tB_0t_0\i,t_0\tB_1t_0\i,\do,t_0\tB_{\y}t_0\i)).\endalign$$
The maps $b_i$ are compatible with the $\G\T\ct$-actions. Let 
$$\align\cg&=\{((b_0,b_1,\do,b_r),(\tb_0,\tb_1,\do,\tb_{\y}))\in\uB\T\utB;\\&
((k(b_0),k(b_1),\do,k(b_r)),(k(\tb_0),k(\tb_1),\do,k(\tb_{\y})))\in\ct\},\endalign$$
a subgroup of $\uB\T\utB$. The $\G\T\ct$-action on $X_6$ extends to a $\G\T\cg$-action on
$X_6$:
$$\align&(g,(b_i),(\tb_j)):((g_0,g_1,\do,g_r),(\tg_0,\tg_1,\do,\tg_{\y}))\m\\&
((k(b_0)g_0b_0\i,k(b_0)g_1b_1\i,\do,k(b_0)g_rb_r\i),
(k(b_0)\tg_0\tb_0\i,k(b_0)\tg_1\tb_1\i,\do,k(b_0)\tg_{\y}\tb_{\y}\i)).\endalign$$

\subhead 6.4\endsubhead
Now 

(a) $b_0$ is a principal $\ct$-bundle
\nl
($\ct$ acts on $X_0$ by restriction of the $\G\T\ct$-action) and induces an isomorphism 
$\ct\bsl X_0@>\si>>X$. (See \cite{\HB, 3.4}.)

(b) $b_1$ is a principal $U^{[0,r]}\T U^{[0,\y]}$-bundle.

(c) $b_2$ is a principal $\G$-bundle
\nl
($\G$ acts on $X_1$ by restriction of the $\G\T\ct$-action) and induces an isomorphism
$\G\bsl X_1@>\si>>X_2$. (See \cite{\HB, 3.8}.)

(d) $b_3$ is an isomorphism.

(e) $b_4$ is a quasi-vector bundle (see{\cite{\HB, 3.2}) with fibres of dimension 
$2(\dim U-l(a_0))$.

(f) $b_5$ is a quasi-vector bundle with fibres of dimension $2(\dim U-l(a_0))$.

(g) $b_6$ is an isomorphism.

(h) $b_7$ is a principal $U^{[0,r]}\T U^{[0,\y]}$-bundle.

(i) $b_8$ is an isomorphism. (See \cite{\HB, 3.24}.)

\subhead 6.5\endsubhead
Now $\ct$ is naturally a subgroup of $\uB\T\utB$ and $\G$ is a subgroup of $\G\T\G$ (the 
diagonal) hence $\G\T\ct$ is a subgroup of $\G\T\G\T\uB\T\utB$. Hence the actions in 6.2 
give by restriction actions of $\G\T\ct$ on $Z\T\tZ$ and $Z_i\T\tZ_i(i\in[0,2])$ and the 
equivariant structures on $\bcl\bxt\tbcl$ and $\bcl_i\bxt\tbcl_i$ in 6.2 restrict to 
$\G\T\ct$-equivariant structures on these local systems. Since $e$ and $e_i(i\in[0,6])$ are
compatible with the $\G\T\ct$-actions we see that the local systems 
$\ce=e^*(\bcl\bxt\tbcl)$, $\ce_i=e_i^*(\bcl_i\bxt\tbcl_i)$, $(i\in[0,2])$ and
$\ce_i=e_i^*(\bcl_2\bxt\tbcl_2)$, $(i\in[3,6])$ have natural $\G\T\ct$-equivariant 
structures. Moreover, the $\G\T\ct$-equivariant structure on $\ce_6$ extends to a
$\G\T\cg$-equivariant structure since $e_6$ is compatible with the $\G\T\cg$ actions (see
6.3). Since the restriction of the $\G\T\cg$ action on $X_6$ to the subgroup
$U^{[0,r]}\T U^{[0,\y]}$ is the free action which makes $X_6$ a principal bundle over $X_7$
(see 6.4(h)) it follows that there is a well defined local system $\ce_7$ on $X_7$ with a 
natural $\G\T\ct$-equivariant structure such that $b_7^*\ce_7=\ce_6$. Since $b_8$ is an
isomorphism, there is a well defined local system $\ce_8$ on $X_7$ with a natural 
$\G\T\ct$-equivariant structure such that $b_8^*\ce_8=\ce_7$. We have
$$\align& b_0^*\ce=\ce_0;\qua b_1^*\ce_1=\ce_0;\qua b_2^*\ce_2=\ce_1;
b_3^*\ce_3=\ce_2;\qua b_4^*\ce_3=\ce_4;\qua b_5^*\ce_5=\ce_4;\\&
\qua b_6^*\ce_5=\ce_6;b_7^*\ce_7=\ce_6;\qua b_8^*\ce_8=\ce_7;\tag a\endalign$$
compatibly with the $\G\T\ct$-equivariant structures. We show:

(b) {\it $\ct$ acts trivially on any stalk of $\ce$.}

(c) {\it $\G$ acts trivially on any stalk of $\ce_i$ (if $i\in[2,8]$).}
\nl
(b) follows from the fact that $\uB\T\utB$ acts trivially on $Z\T\tZ$ and, being connected,
it acts trivially on any stalk of $\bcl\bxt\tbcl$. Now (c) follows from the fact that 
$\G\T\G$ acts trivially on $Z_2\T\tZ_2$ and on any stalk of $\bcl_2\bxt\tbcl_2$.

\subhead 6.6\endsubhead
For any $h\in[0,\y]$ let $\fY_h$ be the set of all $(\tB_0,\tB_1,\do,\tB_h)\in\cb^{[0,h]}$
such that 

$\po(B,\tB_j)=a_j(j\in[0,h])$,  

$\po(\tB_{j-1},\tB_j)=\ts_j$ if $j\in[1,h],j\n\cits$,

$\po(\tB_{j-1},\tB_j)\in\{1,\ts_j\}$ if $j\in[1,h]\cap\cits$,

$\tB_0=\dda_0B\dda_0\i$.
\nl
Note that $\fY_0$ is a point. Moreover, if $h\in[1,\y]$ then we have an obvious map 
$\fY_h@>>>\fY_{h-1}$ which is 

-an isomorphism if $a_h=a_{h-1}\ts_h<a_{h-1}$;

-an isomorphism if $h\in\cits$, $a_{h-1}\ts_h>a_{h-1}=a_h$;

-a line bundle minus the zero section if $h\n\cits$, $a_{h-1}\ts_h<a_{h-1}=a_h$;

-a line bundle if $a_h=a_{h-1}\ts_h>a_{h-1}$;

-a line bundle if $h\in\cits$, $a_{h-1}\ts_h<a_{h-1}=a_h$.
\nl
For any $h\in[0,r]$ let $\fY_{\y+h}$ be the set of all 
$$((B_0,B_1,\do,B_h),(\tB_0,\tB_1,\do,\tB_{\y}))\in\cb^{[0,h]}\T\cb^{[0,\y]}$$ 
such that 

$\po(B_0,\tB_j)=a_j(j\in[0,\y])$, $\po(B_i,\tB_{\y})=a_{\y+i}(i\in[0,h])$,

$\po(B_{i-1},B_i)=s_i$ if $i\in[1,h],j\n\ciss$,

$\po(B_{i-1},B_i)\in\{1,s_i\}$ if $i\in [1,h]\cap\ciss$,

$\po(\tB_{j-1},\tB_j)=\ts_j$ if $j\in[1,\y],j\n\cits$,

$\po(\tB_{j-1},\tB_j)\in\{1,\ts_j\}$ if $j\in\cits$,

$B_0=B,\tB_0=\dda_0B\dda_0\i$.
\nl
Note that $\fY_{\y}$ in the last definition may be identified with $\fY_{\y}$ in the 
earlier one. Moreover, if $h\in[1,r]$ then we have an obvious map 
$\fY_{\y+h}@>>>\fY_{\y+h-1}$ which is 

-an isomorphism if $a_{\y+h}=s_ha_{\y+h-1}<a_{\y+h-1}$;

-an isomorphism if $h\in\ciss$, $s_ha_{\y+h-1}>a_{\y+h-1}=a_{\y+h}$;

-a line bundle minus the zero section if $h\n\ciss$, $s_ha_{\y+h-1}<a_{\y+h-1}=a_{\y+h}$;

-a line bundle if $a_{\y+h}=s_ha_{\y+h-1}>a_{\y+h-1}$;

-a line bundle if $h\in\ciss$, $s_ha_{\y+h-1}<a_{\y+h-1}=a_{\y+h}$.

\subhead 6.7\endsubhead
Assume that for some $h\in[1,r]-\ciss$ we have $a:=a_{\y+h}=a_{\y+h-1}$, $s:=s_h$, $sa<a$.
We show that 

(a) $H^*_c(X_8,\ce_8)=0$.
\nl
We have an obvious map $\ph:X_8@>>>\fY_{\y+h-1}$. (See 6.6.) It is enough to show that for
any $p=((B_0,\do,B_{h-1}),(\tB_0,\tB_1,\do,tB_{\y}))\in\fY_{\y+h-1}$ we have 
$H^*_c(\fY',\ce_8)=0$ where $\fY'=\ph\i(p)$. We may identify $\fY'$ with the set of all 
$(B_h,B_{h+1},\do,B_r)\in\cb^{[h,r]}$ such that

$\po(B_i,\tB_{\y})=a_{\y+i}(i\in[h,r])$,

$\po(B_{i-1},B_i)\in\{s_i,1\}$ if $i\in[h+1,r]\cap\ciss$.

$\po(B_{i-1},B_i)=s_i$ if $i\in[h,r],i\n\ciss$.
\nl
Since $sa<a$, there is a unique $D\in\cb$ such that $\po(B_{h-1},D)=s$,
$\po(D,\tB_{\y})=sa$. Pick $C\in\cb$ such that $\po(\tB_{\y},C)=a\i w_\II$. Since 
$\po(D,C)=sw_\II$, $V:=U_D\cap U_C$ is a one dimensional connected unipotent group. Now 
$\tB_{\y}$ is the unique Borel such that $\po(D,\tB_{\y})=sa$, $\po(\tB_{\y},C)=a\i w_\II$.
Hence it is normalized by $V$ and $V\sub U_{\tB_{\y}}$. Since $\po(B_{h-1},C)=w_\II$, we 
have $U_{B_{h-1}}\cap U_C=\{1\}$ hence $U_{B_{h-1}}\cap V=\{1\}$. Let 
$$\align&\Xi=\{E\in\cb;\po(B_{h-1},E)=s,\po(E,\tB_{\y})=a\}\\&
=\{E\in\cb;\po(B_{h-1},E)=s,\po(E,D)=s\}.\endalign$$ 
Since $V\sub U_D$, $V\cap U_{B_{h-1}}=\{1\}$, $V$ acts simply transitively (by conjugation)
on the affine line $\{E'\in\cb;\po(E',D)=s\}=\Xi\sqc\{B_{h-1}\}$. Pick $B_h\in\xi$. Define 
$v_0\in V-\{1\}$ by $v_0B_hv_0\i=B_{h-1}$. Let 
$$\fY''=\{(B_{h+1},\do,B_r)\in\cb^{[h+1,r]};(B_h,B_{h+1},\do,B_r)\in\fY'\}.$$ 
The map $\z:(V-\{v_0\})\T\fY''@>>>\fY'$, 
$$(v,(B_{h+1},\do,B_r))\m(vB_hv\i,vB_{h+1}v\i,\do,vB_rv\i)$$
is an isomorphism. Hence it is enough to show that 
$$H^*_c((V-\{v_0\})\T\fY'',\z^*\ce_8)=0.$$
Let $\p'':(V-\{v_0\})\T\fY''@>>>\fY''$ be the projection. It is enough to show that for any
$p'=(B_{h+1},\do,B_r)\in\fY''$ we have $H^*_c(\p''{}\i(p'),\z^*\ce_8)=0$ or equivalently
that $H^*_c(V-\{v_0\},\z'{}^*\ce_8)=0$ where $\z':V-\{v_0\}@>>>X_8$ is  
$$v\m((B_0,B_1,\do,B_{h-1},vB_hv\i,vB_{h+1}v\i,\do,vB_rv\i),(\tB_0,\tB_1,\do,\tB_{\y})).$$
Since $\z'{}^*\ce_8$ is a local system of rank $1$ on $V-\{v_0\}\cong\kk^*$ with monodromy
of finite order invertible in $\kk$, it is enough to show that $\z'{}^*\ce_8\not\cong\bbq$.
We can find $\e=((g_0,g_1,\do,g_r),(\tg_0,\tg_1,\do,\tg_{\y}))\in X_6$ such that
$$b_8(b_7(\e))=((B_0,B_1,\do,B_{h-1},B_h,B_{h+1},\do,B_r),(\tB_0,\tB_1,\do,\tB_{\y})),$$
$g_{h-1}\i\tg_{\y}=\dda$ and $V=g_{h-1}y_{s_h}(\kk)g_{h-1}\i$. Define $\l_0\in\kk$ by 
$y_{s_h}(\l_0)g_{h-1}\i g_h\in B$. Define $\tiz:\kk-\{\l_0\}@>>>X_6$ by 
$$\l\m((g_0,g_1,\do,g_{h-1},g_{h-1}y_{s_h}(\l)g_{h-1}\i g_h,\do,
g_{h-1}y_{s_h}(\l)g_{h-1}\i g_r),(\tg_0,\tg_1,\do,\tg_{\y})).$$
We have $\z'{}^*\ce_8=\tiz^*\ce_6$. It is enough to prove that $\tiz^*\ce_6\not\cong\bbq$ 
or that $\tiz^*e_6^*(\bcl_2\bxt\tbcl_2)\not\cong\bbq$. Note that 
$e_6\tiz:\kk-\{\l_0\}@>>>Z_2\T\tZ_2$ is 
$$\l\m((y_1,y_2,y_{h-1},y_{s_h}(\l)y_h,y_{h+1},\do,y_r),(\ty_1,\ty_2,\do,\ty_{\y}))$$
where 
$$((y_1,y_2,y_{h-1},y_h,y_{h+1},\do,y_r),(\ty_1,\ty_2,\do,\ty_{\y}))=e_6(\e).$$
From this and the definition of $\bcl_2=\uucl$ (see 2.4) we see that 
$(e_6\tiz)^*(\bcl_2\bxt\tbcl_2)$ is isomorphic to the inverse image of $\cl$ under a map 
$\kk-\{\l_0\}@>>>T$ of the form 
$$\l\m t'\ds_1\do\ds_{h-1}k(y_{s_h}(\l-\l_0))\ds_h\i\do\ds_1\i$$
where $t'$ is a fixed element of $T$. This is also of the form $\l\m\b(\l-\l_0)t'$ 
where $\b:\kk^*@>>>T$ is one of the two coroots with associated reflection 
$s_1s_2\do s_h\do s_2s_1$. The desired result follows from the fact that 
$\b^*(\cl)\not\cong\bbq$. (Recall that $h\n\ciss$).

\subhead 6.8\endsubhead
Assume that for some $h\in[1,\y]-\cits$ we have $a:=a_h=a_{h-1}$, $s:=\ts_h$, $as<a$. We
have 

(a) $H^*_c(X_8,\ce_8)=0$.
\nl
The proof is entirely similar to that of 6.7(a).

\subhead 6.9\endsubhead
Let 
$$\align&N_\aa=\\&|\{h\in[1,\y];a_{h-1}\le a_h\ge a_{h-1}\ts_h\}|
+|\{h\in[1,r];a_{\y+h-1}\le a_{\y+h}\ge s_ha_{\y+h-1}\}|.\endalign$$
Consider the following condition on $\aa$:

(a) {\it for any $h\in[1,r]-\ciss$ we have $a_{\y+h}=s_ha_{\y+h-1}$; for any 
$h\in[1,\y]-\cits$ we have $a_h=a_{h-1}\ts_h$.}
\nl
If $\aa$ satisfies (a), we have the following results:

(b) {\it $X_8$ is isomorphic to an affine space of dimension $N_\aa$}.

(c) {\it $X$ has pure dimension $N_\aa$}.
\nl
Indeed, from the results in 6.6 we see by induction on $h\in[0,\r]$ that $\fY_h$ is an 
affine space. (In this case each of the maps $\fY_h@>>>\fY_{h-1}$, ($h\in[1,\r]$) in 6.6 is
either an isomorphism or an affine line bundle.) The same argument yields the dimension of
each $\fY_h$. This yields (b). Now (c) follows from (b) and the results in 6.4.

\subhead 6.10\endsubhead
We assume that $\aa$ satisfies 6.9(a). Define $x=((g_i),\tg_j))\in G^{[0,r]}\T G^{[0,\y]}$
by $g_i=\dda_{\y}\dda_{\y+i}\i$, $\tg_j=\dda_j$. Our assumption on $\aa$ implies that 
$x\in X_6$. From the definitions we see that $x$ is a fixed point of the $\ct$-action 6.3
on $X_6$. Hence, in the $\G\T\ct$-structure of $\ce_6$, $\ct$ acts on the stalk 
$\ce_{6,x}$ of $\ce_6$ at $x$. We show:

(a) {\it the $\ct$-action on $\ce_{6,x}$ is trivial if and only if 
$\cl\ot(a_0\i)^*\tcl\cong\bbq$.}
\nl
Let 
$$\cj=\{i\in[1,r];a_{\y+i-1}=a_{\y+i}\}, \tcj=\{j\in[1,\y];a_{j-1}=a_j\}.$$
From our assumption on $\aa$ we see that $\cj\sub\ciss$, $\tcj\sub\cits$. Let 
$f^\cj:\cy^\cj@>>>T$ be as in 2.4; let $\tf^{\tcj}:\tcy^{\tcj}@>>>T$ be the analogous map
defined in terms of $\tss,\tcl,\tcj$ instead of $\ss,\cl,\cj$. Define sequences 
$\sscj=(s'_1,s'_2,\do,s'_r)$, $\tsstcj=(\ts'_1,\ts'_2,\do,\ts'_{\y})$ by

$s'_i=1$ if $i\in\cj$, $s'_i=s_i$ if $i\in[1,r]-\cj$, $\ts'_j=1$ if $j\in\tcj$, 
$\ts'_j=\ts_j$ if $j\in[1,\y]-\tcj$.
\nl
We have $e_6(x)=((y_1,\do,y_r),(\ty_1,\do,\ty_{\y}))$ where 
$y_i=\dda_{\y+i-1}\dda_{\y+i}\i$ for $i\in[1,r]$ and $\ty_j=\dda_{j-1}\i\dda_j$ for 
$j\in[1,\y]$. Equivalently, $y_i=\ds'_i$ for $i\in[1,r]$ and $\ty_j=\dot{\ts'}_j$ for 
$j\in[1,\y]$. Thus we have $e_6(x)\in\cy^\cj\T\tcy^{\tcj}$. Moreover, we have
$(f^\cj\T\tf^{\tcj})e_6(x)=(1,1)\in T\T T$. It is enough to show that the $\ct$-action on 
the stalk of $\bcl_2\bxt\tbcl_2$ at $e_6(x)$ is trivial if and only if 
$\cl\ot(a_0\i)^*\tcl\cong\bbq$. (The $\ct$-equivariant structure on $\bcl_2\bxt\tbcl_2$ is
obtained by restricting the $\uB\T\utB$-equivariant structure in 6.2.) By 2.4(a) we have 
canonically 
$$(\bcl_2\bxt\tbcl_2)|_{\cy^\cj\T\tcy^{\tcj}}=(f^\cj\T\tf^{\tcj})^*(\cl\bxt\tcl).$$
Now $\cy^\cj\T\tcy^{\tcj}$ is stable under the $\uB\T\utB$-action on $Z_2\T\tZ_2$ and
the previous equality shows that $(f^\cj\T\tf^{\tcj})^*(\cl\bxt\tcl)$ has a 
$\uB\T\utB$-equivariant structure; this structure is unique since $\uB\T\utB$ is connected.
By restriction to the subgroup $\ct$ of $\uB\T\utB$ we obtain a $\ct$-equivariant structure
on $(f^\cj\T\tf^{\tcj})^*(\cl\bxt\tcl)$ and it is enough to show that the $\ct$-action on 
the stalk of $(f^\cj\T\tf^{\tcj})^*(\cl\bxt\tcl)$ at $e_6(x)$ is trivial if and only if 
$\cl\ot(a_0\i)^*\tcl\cong\bbq$. We define a $\uB\T\utB$-action on $T\T T$ by
$$((b_i),(\tb_j)):(t,\tit)\m(k(b_0)t([\sscj]Fk(b_0\i)),k(\tb_0)\tit([\tsstcj]Fk(\tb_0\i))).
$$
From the definitions we see that $f^\cj\T\tf^{\tcj}:Z_2\T\tZ_2@>>>T\T T$ is compatible with
the $\uB\T\utB$-actions. Moreover, $\cl\bxt\tcl$ is equivariant for the $\uB\T\utB$-action
on $T\T T$ as above. (In fact it is equivariant for the action of the bigger group 
$B^{[0,r]}\T B^{[0,\y]}$ given by the same formula; this follows from 1.8 using that 
$\cl\in\cs(T)^{[\sscj]F}$ and $\cl\in\cs(T)^{[\tsstcj]F}$.) By restriction to the subgroup
$\ct$ of $\uB\T\utB$ we obtain a $\ct$-equivariant structure on $\cl\bxt\tcl$. We may 
identify the stalk of $(f^\cj\T\tf^{\tcj})^*(\cl\bxt\tcl)$ at $e_6(x)$ with the stalk of 
$\cl\bxt\tcl$ at $(1,1)$ as $\ct$-modules. Let $\c:\ct@>>>\bbq^*$ be the character by which
$\ct$ acts on the stalk of $\cl\bxt\tcl$ at $(1,1)$. It is enough to show that $\c=1$ if 
and only if $\cl\ot(a_0\i)^*\tcl\cong\bbq$. 

Now let $T'$ be the subgroup of $B^{[0,r]}\T B^{[0,\y]}$ consising of all elements \lb
$((t_0,t_1,\do,t_r),(\tit_0,\tit_1,\do,\ttir))$ with coordinates in $T$ such that 
$\tit_j=a_j\i(t_0)\text{ for }j\in[0,\y]$, $t_i=a_{\y+i}(\ttir)\text{ for }i\in[0,r]$. We 
have $\ct\sub T'$ and $((t_i),(\tit_j))\m t_0$ is an isomorphism $T'@>\sim>>T$. Moreover
$\cl\bxt\tcl$ is $T'$-equivariant where $T'$ acts on $T\T T$ by restriction of the 
$B^{[0,r]}\T B^{[0,\y]}$ action. Now $T'$ acts on $T$ by $t_0:t\m t_0tF'(t_0\i)$ where 
$F':T@>>>T$, $F'(t_0)=(a_{\y}a_\r\i F)(t_0)$ is the Frobenius map for an $\FF_q$-rational 
structure on $T$. The map $m:T@>>>T\T T$, $t\m(t,a_0\i(t))$ is compatible with the 
$T'$-actions. Hence $m^*(\cl\bxt\tcl)=\cl\ot(a_0\i)^*\tcl$ is a $T'$-equivariant local
system on $T$ and the natural action of $\ct$ at any stalk of this local system is via 
$\c$. If we identify $T'$ with $T$ as above, $\ct$ becomes the subgroup $T^{F'}$ of $T$. It
remains to use the bijection (i)-(iv) in 1.9, with $\cl$ replaced by $\cl\ot(a_0\i)^*\tcl$.

\subhead 6.11\endsubhead
We assume that $\aa$ satisfies 6.9(a). We show:

(a) $\ce_6\cong\bbq$.
\nl
We write $\bcl_2=\uucl=\ot_{i\in[1,r]}\cf_i$ as in 2.4. Similarly
$\tbcl_2=\ot_{j\in[1,\y]}\tcf_j$ where $\tcf_j$ are local systems on $\tZ_2$. From the 
definitions we have 
$$\ce_6=\ot_{i\in[1,r]}e_6^*(\cf_i\bxt\bbq)\ot\ot_{j\in[1,\y]}e_6^*(\bbq\bxt\tcf_j).$$
It is enough to show:

(b) $e_6^*(\cf_i\bxt\bbq)\cong\bbq$ for any $i\in[1,r]$;

(c) $e_6^*(\bbq\bxt\tcf_j)\cong\bbq$ for any $j\in[1,\y]$.
\nl
We prove (b). The general case can be reduced to the case where $G$ has simply connected
derived group, which we now assume. Let $\psi:Z_2\T\tZ_2@>>>Z_2$ be the projection. Let 
$p_i,f_{s_i}$ be as in 2.4. Let $\tf_i$ be the obvious map from $P_{s_i}$ to the quotient 
of $P_{s_i}/U_{P_{s_i}}$ by its derived subgroup. Then $e_6^*(\cf_i\ot\bbq)$ is the inverse
image of a local system of rank $1$ under 

$\tf_ip_i\psi e_6$ (if $i\in\ciss$), or $f_{s_i}p_i\psi e_6$ (if $i\in[1,r]-\ciss$).
\nl
It is enough to show that the image of $\tf_ip_i\psi e_6$ (if $i\in\ciss$) or of
$f_{s_i}p_i\psi e_6$ (if $i\in[1,r]-\ciss$) is a point.

Let $\xi=((g_0,g_1,\do,g_r),(\tg_0,\tg_1,\do,\tg_{\y}))\in X_6$. We have 
$p_i\psi e_6(\xi)=g_{i-1}\i g_i$. We have $g_i\i\tg_{\y}\in U\dda_{\y+i}U$, 
$g_{i-1}\i\tg_{\y}\in U\dda_{\y+i-1}U$, hence
$g_{i-1}\i g_i\in u_1\dda_{\y+i-1}u_2\dda_{\y+i}\i u_3$ with $u_1,u_2,u_3$ in $U$. Moreover
$g_{i-1}\i g_i\in P_{s_i}$.

{\it Case 1.} Assume that $i\in\ciss$ and $a_{\y+i-1}=a_{\y+i}$. Then 
$\dda_{\y+i-1}=\dda_{\y+i}\i$ and
$\dda_{\y+i-1}u_2\dda_{\y+i}\i=\dda_{\y+i-1}u_2\dda_{\y+i-1}\i\in P_{s_i}$ is unipotent.
Hence $g_{i-1}\i g_i$ is a product of three unipotent elements of $P_{s_i}$ so that 
$\tf_i(g_{i-1}\i g_i)=1$.

{\it Case 2.} Assume that $i\in\ciss$ and $a_{\y+i-1}=s_ia_{\y+i}$. Then
$\dda_{\y+i}=\ds_i\i\dda_{\y+i-1}\i$ and
$$\dda_{\y+i-1}u_2\dda_{\y+i}\i=\dda_{\y+i-1}u_2\dda_{\y+i-1}\i\ds_i=u'_2\ds_i$$
where $u'_2\in G$ is unipotent. Thus $g_{i-1}\i g_i=u_1u'_2\ds_iu_3$. Since 
$u_1u'_2\ds_iu_3,u_1,\ds_i,u_3$ belong to $P_{s_i}$ we see that $u'_2$ is a unipotent 
element of $P_{s_i}$. Note also that $\ds_i$ belongs to the derived subgroup of $P_{s_i}$.
We see that $\tf_i(g_{i-1}\i g_i)=1$.

{\it Case 3.} Assume that $i\n\ciss$. By our assumption we have $a_{\y+i-1}=s_ia_{\y+i}$. 
Then $\dda_{\y+i}=\ds_i\i\dda_{\y+i-1}\i$ and as before we have 
$g_{i-1}\i g_i=u_1u'_2\ds_iu_3$ where $u'_2$ belongs to the unipotent group 
$P_{s_i}\cap\dda_{\y+i-1}U\dda_{\y+i-1}\i$. This unipotent group is normalized by $T$. By
\cite{\BO, 14.4} we have $P_{s_i}\cap\dda_{\y+i-1}U\dda_{\y+i-1}\i=U_1U_2\do U_n$ where
$U_1,U_2,\do,U_n$ are the connected one dimensional unipotent subgroups of
$P_{s_i}\cap\dda_{\y+i-1}U\dda_{\y+i-1}\i$ normalized by $T$, in any order. Thus any of 
$U_1,U_2,\do,U_n$ is either $y_{s_i}(\kk)$ or else is contained in $U$; moreover we can 
assume that $U_1,U_2,\do,U_{n-1}$ are contained in $U$ and $U_n=y_{s_i}(\kk)$. Thus 
$u'_2\in Uy_{s_i}(\kk)$ and
$$g_{i-1}\i g_i\in Uy_{s_i}(\kk)\ds_i U=U\ds_ix_{s_i}(\kk)U=U\ds_iU.$$
Note that $f_{s_i}(g_{i-1}\i g_i)=1$.

We see that in Case 1 and 2 we have $\tf_ip_i\psi e_6(\x)=1$ for any $\x\in X_6$. In Case 3
we have $f_{s_i}p_i\psi e_6(\x)=1$ for any $\x\in X_6$. This completes the proof of (b).

The proof of (c) is entirely similar.

\subhead 6.12\endsubhead
We set $X_{-1}=X$, $\ce_{-1}=\ce$. Let $\d_i=\dim X_i (i\in[-1,8])$. For $i\in[-1,8]$ we 
show:

(a) {\it Assume that $\aa$ satisfies 6.9(a); if $n\ne2\d_i$ or if $\tcl\not\cong a_0^*\ccl$
then\lb  $H^n_c(X_i,\ce_i)^{\G\T\ct}=0$; if $\tcl=w^*\ccl$ where $w\in W$ and 
$w^*\cl\cong a_0^*\cl$ then \lb $H^{2\d_i}_c(X_i,\ce_i)^{\G\T\ct}(\d_i)=\bbq$ canonically. 
Assume that $\aa$ does not satisfy 6.9(a); then $H^n_c(X_i,\ce_i)^{\G\T\ct}=0$ for all $n$.
}
\nl
Here $(\d_i)$ is a Tate twist, The upper index denotes $\G\T\ct$-invariants (the action of
$\G\T\ct$ comes from the $\G\T\ct$-equivariant structure of $\ce_i$).

Let $P_i$ be the statement of (a). From 6.4 we see that the statements $P_0,P_1$ are 
equivalent and that the statements $P_2,P_3,\do,P_8$ are equivalent. From 6.4 we see also 
that $H^n_c(X_2,\ce_2)=H^n_c(X_1,\ce_1)^\G$ and that
$H^n_c(X_{-1},\ce_{-1})=H^n_c(X_0,\ce_0)^\ct$ for any $n$. Since $\G$ acts trivially on any
stalk of $\ce_2$ (see 6.5(c)) and $\ct$ acts trivially on any stalk of $\ce_{-1}$ (see
6.5(b)) it follows that $H^n_c(X_2,\ce_2)^\G=H^n_c(X_1,\ce_1)^\G$ and that
$H^n_c(X_{-1},\ce_{-1})^{\ct}=H^n_c(X_0,\ce_0)^\ct$ for any $n$. This shows that the
statements $P_1,P_2$ are equivalent and the statements $P_{-1},P_0$ are equivalent. We see
that the statements $P_{-1},P_0,\do,P_8$ are all equivalent. Thus it is enough to show that
$P_8$ holds.

Assume first that $\aa$ does not satisfy 6.9(a). Then the result follows from 6.7(a) or
6.8(a).

Next we assume that $\aa$ satisfies 6.9(a). Let 
$\bx=b_8(b_7(x))\in X_8$ where $x\in X_6$ is as in 6.10. We can write 6.11(a) in the form 
$b_7^*b_8^*\ce_8\cong\bbq$. Since $b_8$ is an isomorphism and $b_7$ is a principal
$U^{[0,r]}\T U^{[0,\y]}$-bundle, it follows that $\ce_8\cong\bbq$. Using this and 6.9(b) we
see that 
$$H^n_c(X_8,\ce_8)=0\text{ for }i\ne2\d_8, H^{2\d_8}_c(X_8,\ce_8)(\d_8)=\ce_{8,\bx},$$
where $\ce_{8,\bx}$ is the stalk of $\ce_8$ at $\bx$. Moreover the last equality is 
compatible with the natural $\ct$-actions (coming from the $\ct$-equivariant structure of 
$\ce_8$). Now $\ce_{8,\bx}$ may be canonically identified with the stalk $\ce_{6,x}$ of 
$\ce_6$ at $x$. Thus we have $H^{2\d_8}_c(X_8,\ce_8)(\d_8)=\ce_{6,x}$, compatibly with the
$\ct$-actions. Using 6.10(a) we see that the $\ct$-action on $H^{2\d_8}_c(X_8,\ce_8)$ is 
trivial if and only if $\cl\ot(a_0\i)^*\tcl\cong\bbq$. Taking $\ct$-invariants we see that

$H^{2\d_8}_c(X_8,\ce_8)^\ct(\d_8)=\ce_{6,x}$ if $\cl\ot(a_0\i)^*\tcl\cong\bbq$;

$H^{2\d_8}_c(X_8,\ce_8)^\ct=0$ if $\cl\ot(a_0\i)^*\tcl\not\cong\bbq$.
\nl
By 6.5(c), $\G$ acts trivially on $H^{2\d_8}_c(X_8,\ce_8)$ hence

$H^{2\d_8}_c(X_8,\ce_8)^{\G\T\ct}=H^{2\d_8}_c(X_8,\ce_8)^\ct$.

We see that the first assertion in (a) for $i=8$ holds. It remains to prove that, if 
$\tcl=w^*\ccl$ where $w\in W$ and $w^*\cl\cong a_0^*\cl$, then $\ce_{6,x}=\bbq$
canonically. By the proof of 6.10(a), $\ce_{6,x}$ is canonically isomorphic to the stalk of
$\cl\ot(a_0\i)^*\tcl$ at $1$ that is to the stalk of $\cl\ot(wa_0\i)^*\ccl$ at $1$. The
stalk of $(wa_0\i)^*\ccl$ at $1$ is the same as the stalk of $\ccl$ at $wa_0(1)=1$. Thus
$\ce_{6,x}$ is canonically isomorphic to the stalk of $\cl\ot\ccl=\bbq$ at $1$. Thus $P_8$
holds. We see that (a) holds for $i\in[-1,8]$.

\subhead 6.13\endsubhead
We now write $X_\aa,\ce_\aa,e_\aa$ instead of $X,\ce,e$ in 6.1. We identify $X_\aa$ with 
its image under the imbedding $e_\aa$. Let $\bX_\aa$ be the closure of $X_\aa$ in $Z\T\tZ$.
Recall that $\ce_\aa=e_\aa^*(\bcl\bxt\tbcl)$. Let $\be_\aa:\bX_\aa@>>>Z\T\tZ$ by the
inclusion. By 6.9(c), $X_\aa$ (hence also $\bX_\aa$) has pure dimension $N_\aa$ and 
$\dim(\bX_\aa-X_\aa)<N_\aa$. Hence the natural map 
$$H^{2N_\aa}_c(X_\aa,e^*_\aa(\bcl\bxt\tbcl))@>>>
H^{2N_\aa}_c(\bX_\aa,\be^*_\aa(\bcl\bxt\tbcl))\tag a$$
(induced by the open imbedding $X_\aa\sub\bX_\aa$) is an isomorphism. Let 
$$\x'_\aa:H^{2N_\aa}_c(Z\T\tZ,(\bcl\bxt\tbcl))@>>>
H^{2N_\aa}_c(X_\aa,e^*_\aa(\bcl\bxt\tbcl))$$
be the linear map obtained by composing the linear map
$$H^{2N_\aa}_c(Z\T\tZ,(\bcl\bxt\tbcl))@>>>H^{2N_\aa}_c(\bX_\aa,\be^*_\aa(\bcl\bxt\tbcl))$$
with the inverse of (a). By taking $\G$-invariants and applying a Tate twist we obtain 
from $\x'_\aa$ a linear map
$$\x_\aa:H^{2N_\aa}_c(Z\T\tZ,(\bcl\bxt\tbcl))^\G(N_\aa)@>>>
H^{2N_\aa}_c(X_\aa,e^*_\aa(\bcl\bxt\tbcl))^\G(N_\aa).\tag b$$

\subhead 6.14\endsubhead
Let $\ca$ be the set of all $\aa$ as in 6.1. Note that the subvarieties $X_\aa$
($\aa\in\ca$) form a partition of $Z\T\tZ$. 

(a) {\it For any $n\in\ZZ$, the linear map
$$\x_n:H^{2n}_c(Z\T\tZ,\bcl\bxt\tbcl)^\G(n)@>>>
\op_{\aa\in\ca;N_\aa=n}H^{2n}_c(X_\aa,e^*_\aa(\bcl\bxt\tbcl))^\G(n)$$
whose components are the maps $\x_\aa$ $(\aa\in\ca,N_\aa=n)$, is an isomorphism. Moreover,
$H^{n'}_c(Z\T\tZ,\bcl\bxt\tbcl)^\G=0$ for any odd $n'$.}
\nl
The proof is almost identical to that of \cite{\HB, 2.7}. (We use the fact that \lb
$H^m_c(X_\aa,e^*_\aa(\bcl\bxt\tbcl))^\G=0$ for any $\aa\in\ca$ and any $m\ne2N_\aa$; this 
follows from 6.12(a) for $i=-1$; note that, in view of 6.5(b), 6.12(a) for $i=-1$ remains 
valid if $\G\T\ct$-invariants are replaced by $\G$-invariants.)

From (a) and 6.12(a) for $i=-1$ we see that the following holds.

(b) {\it If for any $w\in W$ we have $\tcl\not\cong w^*\ccl$, then
$H^m_c(Z\T\tZ,\bcl\bxt\tbcl)^\G=0$ for any $m\in\ZZ$.}
\nl
Now assume that $\tcl=w^*\ccl$ for some $w\in W$. Let 

$\ca_w=\{\aa=(a_0,a_1,\do,a_\r)\in\ca; \aa \text{ satisfies 6.9(a)}, 
a_0^*\cl\cong w^*\cl.\}$.
\nl
Note that for $w=1$ we have $\ca_1=\ca(W,\boc,\cl,\ss,\tss)$, see 5.9.

By 6.12(a) for $i=-1$, the summand in the target of $\x_n$ corresponding to $\aa$ is 
canonically $\bbq$ if $\aa\in\ca_w$ and is $0$ if $\aa\n\ca_w$. Hence for any $\aa\in\ca_w$
with $N_\aa=n$ there is a unique element $b_\aa^w\in H^{2n}_c(Z\T\tZ,\bcl\bxt\tbcl)^\G(n)$
such that $\x_n(b_\aa^w)$ is contained in the summand corresponding to $\aa$ and, as an 
element of that summand, it corresponds to $1\in\bbq$. Moreover,

(c) {\it $\{b_\aa^w;\aa\in\ca_w,N_\aa=n\}$ is a $\bbq$-basis of
$H^{2n}_c(Z\T\tZ,\bcl\bxt\tbcl)^\G(n)$.}
\nl
Taking $w=1$ and taking direct sum over $n$ we obtain

(d) {\it $\{b_\aa^1;\aa\in\ca(W,\boc,\cl,\ss,\tss)\}$ is a $\bbq$-basis of
$\op_nH^{2n}_c(Z\T\tZ,\bcl\bxt\tbcl)^\G(n)$.}

\subhead 6.15\endsubhead
More generally, if $J$ is as in 3.1, we set 
$\ca_J=\{\aa=(a_0,a_1,\do,\a_\r)\in\ca; a_0\in W_J\}$. Let 
$$(Z\T\tZ)_J=\{((B_0,B_1,\do,B_r),(\tB_0,\tB_1,\do,\tB_{\y}))\in Z\T\tZ;
P_{B_0,J}=P_{\tB_0,J}\}.$$
Note that the subvarieties $X_\aa$ ($\aa\in\ca_J$) form a partition of $(Z\T\tZ)_J$. As in
6.14 we see that for any $n\in\ZZ$ we have an isomorphism 
$$H^{2n}_c((Z\T\tZ)_J,\bcl\bxt\tbcl)^\G(n)@>\si>>
\op_{\aa\in\ca_J;N_\aa=n}H^{2n}_c(X_\aa,e^*_\aa(\bcl\bxt\tbcl))^\G(n).\tag a$$
Moreover,
$$H^{n'}_c((Z\T\tZ)_J,\bcl\bxt\tbcl)^\G=0\text{ for any odd }n'.\tag b$$
As in 6.14 we see that the following holds:

(c) {\it If for any $w\in W_J$ we have $\tcl\not\cong w^*\ccl$, then
$H^m_c((Z\T\tZ)_J,\bcl\bxt\tbcl)^\G=0$ for any $m\in\ZZ$.}
\nl
Define $h:(Z\T\tZ)@>>>\cp_J$ by 
$((B_0,B_1,\do,B_r),(\tB_0,\tB_1,\do,\tB_{\y}))\m P_{B_0,J}=P_{\tB_0,J}$. We have 
$$H^m_c((Z\T\tZ)_J,\bcl\bxt\tbcl)=H^m_c(\cp_J,h_!(\bcl\bxt\tbcl))=
H^m_c(\cp_J,(\Upss_!\bcl)\ot(\Uptss_!\tbcl)).$$
(Notation of 3.2.) From 2.13 we see that 
$$\Upss_!\bcl=\bUp^\ss_!\bcl^\sh\tag d$$ 
(notation of 3.2). By the decomposition theorem \cite{\BBD}, 
$\bUp^\ss_!\bcl^\sh\cong\op_j{}^pH^j(\bUp^\ss_!\bcl^\sh)[-j]$. Hence
$\Upss_!\bcl\cong\op_j{}^pH^j(\Upss_!\bcl)[-j]$. Similarly,
$\Uptss_!\tbcl\cong\op_{\tj}{}^pH^{\tj}(\Uptss_!\tbcl)[-\tj]$. We see that 
$$H^m_c((Z\T\tZ)_J,\bcl\bxt\tbcl)=
\op_{j,j'}H^{m-j-\tj}_c(\cp_J,{}^pH^j(\Upss_!\bcl)\ot{}^pH^{\tj}(\Uptss_!\tbcl))\tag e$$
(compatibly with the $\G$-actions). 

\proclaim{Proposition 6.16} Let $J\sub\II$. Let $A$ be a simple object of $\cm_\G(\cp_J)$.
Assume that $A\dsv_\G{}^pH^j(\Upss_!\bcl)$ and $A\dsv_\G{}^pH^{j'}(\Uptss_!\tbcl)$ where 
$j,j'\in\ZZ$. Then $\tcl\cong w^*\ccl$ for some $w\in W_J$ and $j=j'\mod2$.
\endproclaim
Using 6.15(d) and the fact that $\bUp^\ss$ is proper we see that
$$\fD(\Upss_!\bcl)=\bUp^\ss_!\fD(IC(\bZ^\ss,\bcl))=
\bUp^\ss_!(IC(\bZ^\ss,\bar{\ccl}))[2r]=\Upss_!\bar{\ccl}[2r].$$
Hence $\fD({}^pH^j(\Upss_!\bcl))={}^pH^{-j+2r}(\Upss_!\bar{\ccl})$. We see that \lb
$\fD(A)\dsv_\G{}^pH^{-j+2r}(\Upss_!\bar{\ccl})$. By 1.6(a) we have 
$\dim H^0_c(\cp_J,\fD(A)\ot A)^\G=1$. Since $\fD(A)$ (resp. $A$) is a direct summand of
${}^pH^{-j+2r}(\Upss_!\bar{\ccl})$ (resp. ${}^pH^{j'}(\Uptss_!\tbcl)$) we see that
$$H^0_c(\cp_J,{}^pH^{-j+2r}(\Upss_!\bar{\ccl})\ot{}^pH^{j'}(\Uptss_!\tbcl))^\G\ne0.$$
Using 6.15(e) we see that $H_c^{-j+2r+j'}((Z\T\tZ)_J,\bcl\bxt\tbcl)^\G\ne0$. Using 
6.15(b),(c) we see that the proposition holds.

\subhead 6.17\endsubhead
In the case where $J=\II$, the first assertion of 6.16 reduces to the disjointness theorem
\cite{\DL, 6.2, 6.3} (we use also the equivalence of 3.2(i), 3.2(vi)).

\subhead 6.18\endsubhead
In this subsection we assume that $\tcl=\ccl$ and 

(a) any $s\in\II$ is in the $\boc$-orbit of some $s_i(i\in[1,r])$.
\nl
Since $\dim(Z\T\tZ)=\r$ we have $H^i_c(Z\T\tZ,\bcl\bxt\tbcl)=0$ for $i<2\r$. We show:

(b) {\it $\dim H^{2\r}_c(Z\T\tZ,\bcl\bxt\tbcl)^\G=n_\cl$ 
where $n_\cl$ is $1$ if $R=R_\cl$ and $0$ if $R\ne R_\cl$.}
\nl
By 6.14(c) we have $\dim H^{2\r}_c(Z\T\tZ,\bcl\bxt\tbcl)^\G=|\fA|$ where
$\fA=\{\aa\in\ca_1;N_\aa=\r\}$. Note that $\fA$ is the set of all sequences 
$\aa=(a_0,a_1,\do,a_\r)$ in $W$ such that

$a_{j-1}\i a_j\in\{1,\ts_j\}$ for $j\in\cits$;

$a_{j-1}\i a_j=\ts_j$ for $j\in[1,\y]-\cits$;

$a_{\y+i}a_{\y+i-1}\i\in\{1,s_i\}$ for $i\in\ciss$;

$a_{\y+i}a_{\y+i-1}\i=s_i$ for $i\in[1,r]-\ciss$;

if $h\in[1,\y]$ then $a_{h-1}\le a_h\ge a_{h-1}\ts_h$;

if $h\in[1,r]$ then $a_{\y+h-1}\le a_{\y+h}\ge s_ha_{\y+h-1}$;

$a_\r=\boc(a_0)$;

$a_0^*\cl\cong\cl$.
\nl
If $\aa\in\fA$ then $l(a_0)\le l(a_1)\le\do\le l(a_\r)=l(a_0)$ hence 
$l(a_0)=l(a_1)=\do=l(a_\r)$ and since $a_0\le a_1\le\do\le a_\r$, we have
$a_0=a_1=\do=a_\r$. This forces $\cits=[1,\y]$. Thus, 
$s_1,s_1s_2s_1,\do,s_1s_2\do s_r\do s_2s_1$ are in $W_\cl$ hence $s_i\in W_\cl$ for 
$i\in[1,r]$. Using $([\ss]F)^*\cl\cong\cl$ and $[\ss]\in W_\cl$ we deduce $F^*\cl\cong\cl$.
Hence if $i\in[1,r]$ then $\boc^n(s_i)\in W_\cl$ for any $n\ge1$. Using (a) we deduce that
$\II\sub W_\cl$ hence $W=W_\cl$ and $R=R_\cl$. 

If $R=R_\cl$ then $\fA$ is in bijection with the set 

$\fA'=\{a_0\in W;a_0=\boc(a_0), a_0\ts_h<a_0\text{ for }h\in[1,\y], s_ha_0<a_0\text{ for }
h\in[1,r]\}$.
\nl
We set $I'=\{s\in\II;s=s_h\text{ for some }h\in[1,r]\}$,
$I''=\{s\in\II;s=\ts_h\text{ for some }h\in[1,\y]\}$. We see that

$\fA'=\{a_0\in W;I'\sub L_{a_0},I''\sub R_{a_0},a_0=\boc(a_0)\}$
\nl
where for $w\in W$ we set $L_w=\{s\in\II;sw<w\},R_w=\{s\in\II;ws<w\}$.

For $a_0$ such that $a_0=\boc(a_0)$, the set $L_{a_0}$ is $\boc$-stable; hence the 
condition $I'\sub L_{a_0}$ is equivalent to
$I'\cup\boc(I')\cup\boc^2(I')\cup\do\sub L_{a_0}$ that is (by (a)) to $\II=L_{a_0}$. We see
that $\fA'$ has exactly one element: $w_\II$. This proves (b).

\head 7. The variety $X$\endhead
\subhead 7.1\endsubhead
In this section we study the variety
$$X=\{(B',g)\in\cb\T G;g\i F(g)\in U_{B'}\}.$$
Note that $X$ is an \'etale covering of the smooth connected variety 
$X'=\{(B',u)\in\cb\T G;u\in U_{B'}\}$ of dimension via the map $\r_1:X@>>>X'$, 
$(B',g)\m(B',g\i F(g))$. Since $\dim X'=2d$ where $d=\dim\cb$ we see that $X$ is smooth of
pure dimension $2d$. Here is one of the main results of this section.

\proclaim{Proposition 7.2} If $G$ is simply connected then $X$ is connected.
\endproclaim
The proof is given in 7.16. Note that $X$ is not necessarily connected without the 
assumption that $G$ is simply connected.

\subhead 7.3\endsubhead
Let $w\in W$. Let $i_w:\cb_w@>>>\cb$ be the inclusion. Let $A'\in\Bbb S(\cb_w)$. Let $A$ be
a simple object of $\cm_\G(\cb)$ such that $A\dsv_\G{}^pH^{\cdot}(i_{w!}A')$. We show:
$$A\in\Bbb S(\cb).\tag a$$
We argue by induction on $l(w)$. If $l(w)=0$ then $i_{w!}A'\in\Bbb S(\cb)$ and the result 
is clear. We now assume that $l(w)>0$. We have $A=\cl_w[l(w)]$ where $\cl\in\cs(T)^{wF}$. 
Let $K=\cl_w^\sh[l(w)]$. We have $K\in\Bbb S(\cb)$. Assume first that 
$A\not\dsv_\G{}^pH^{\cdot}(i_{w'!}i_{w'}^*K)$ for any $w'\in W$ such that $w'<w$. Then
$A\not\dsv_\G{}^pH^{\cdot}(u_!u^*K)$ where $u:\cup_{w';w'<w}\cb_{w'}@>>>\cb$ is the 
inclusion. Let $u':\cup_{w';w'\le w}\cb_{w'}@>>>\cb$ be the inclusion. Since
$i_{w!}i_w^*K=A'$ we see that $A\dsv_\G{}^pH^{\cdot}(u'_!u'{}^*K)$. Thus $A\dsv_\G K$ hence
$A\cong K$ and $A\in\Bbb S(\cb)$.

Next we assume that $A\dsv_\G{}^pH^{\cdot}(i_{w'!}i_{w'}^*K)$ for some $w'\in W$ such that
$w'<w$. It follows that there exists $j'$ such that 
$A\dsv_\G{}^pH^{\cdot}(i_{w'!}({}^pH^{j'}(i_{w'}^*K)))$. Hence there exists a simple object
$A''$ of $\cm_\G(\cb_{w'})$ such that $A''\dsv{}^pH^{j'}(i_{w'}^*K)$ and 
$A\dsv_\G{}^pH^{\cdot}(i_{w'!}A'')$. From 4.14 we see that $A''\in\Bbb S(\cb_{w'})$. From 
the induction hypothesis we see that $A\in\Bbb S(\cb)$. This proves (a).

\subhead 7.4\endsubhead
Let $w\in W$. Let $f_w:\ti\cb_w@>>>\cb_w$ be as in 4.7, a finite principal covering. Let 
$U_w=\{u\in U;\dw F(u)\dw\i\in U\}$. Let $\hat\cb_w=\{z\in G;z\i F(z)\in U\dw\}$. Define 
$a':\hat\cb_w@>>>\ti\cb_w$ by $z\m zU$, a principal $U_w$-bundle. Let 

$X_w=\{(B',g)\in\cb\T G;g\i F(g)\in U_{B'},\po(B',F(B'))=w\}$,

$\tX_w=\{(zU_w,g)\in G/U_w\T G;z\i F(z)\in U\dw,g\i F(g)\in zUz\i\}$,

$\tX'_w=\{(z,y)\in G\T G;z\i F(z)\in U\dw,y\i F(y)\in U\dw\}$,

$\tX''_w=U_w\bsl\tX'_w$,
\nl
where $U_w$ acts (freely) by $u:(z,y)\m(zu\i,yu\i)$. Define $\p_w:X_w@>>>\cb_w$ by
$(B',g)\m B'$. Define $\p':\tX_w@>>>\ti\cb_w$ by $(zU_w,g)\m zU$. Under the identification 
$\tX_w=\tX''_w$, $(zU_w,g)\lra(z,gz)U_w$, $\p'$ becomes $\p'':\tX''_w@>>>\ti\cb_w$, 
$(z,y)U_w\m zU$. Define $a:\tX'_w@>>>\tX''_w$ by $(z,y)\m(z,y)U_w$. Define 
$h:\tX'_w@>>>\hat\cb_w$ by $(z,y)\m z$. Define $\g:\tX_w@>>>X_w$ by $(zU_w,g)\m(zBz\i,g)$.
Now $\G$ acts on $\tX''_w$ by $g_0:(z,y)\m(g_0z,y)$. We show:

(i) if $A\in\cm_\G(\ti\cb_w)$ is simple and 
$A\dsv_\G{}^pH^{\cdot}(h_!\bbq)$ then $A\cong\bbq[\dim U]$;

(ii) if $A\in\cm_\G(\ti\cb_w)$ is simple and 
$A\dsv_\G{}^pH^{\cdot}(a'_!h_!\bbq)$ then $A\cong\bbq[l(w)]$;

(iii) if $A\in\cm_\G(\ti\cb_w)$ is simple and  
$A\dsv_\G{}^pH^{\cdot}(\p''_!a_!\bbq)$ then $A\cong\bbq[l(w)]$;

(iv) if $A\in\cm_\G(\ti\cb_w)$ is simple and 
$A\dsv_\G{}^pH^{\cdot}(\p'_!\bbq)$ then $A\cong\bbq[l(w)]$;

(v) if $A\in\cm_\G(\cb_w)$ is simple and 
$A\dsv_\G{}^pH^{\cdot}(\p_{w!}\bbq))$ then $A\in\Bbb S(\cb_w)$;

(vi) if $A\in\cm_\G(\cb)$ is simple and $A\dsv_\G{}^pH^{\cdot}(i_{w!}\p_{w!}\bbq)$ then 
$A\in\Bbb S(\cb)$.
\nl
Now (i) is obvious; (ii) follows from (i) using $a'_!\bbq=\bbq[-2\d](-\d)$ where 
$\d=\dim U_w$; (iii) follows from (ii) using $a'h=\p''a$; (iv) follows from (iii) using 
$a_!\bbq=\bbq[-2\d](-\d)$ and the identification $\p'=\p''$. Since the diagram formed by 
$\p',f_w,\g,\p_w$ is cartesian we have $\p'_!\bbq=\p'_!\g^*\bbq=f_w^*\p_{w!}\bbq$. 

We prove (v). Let $A$ be as in (v). Let $A'$ be a simple object of $\cm_\G(\ti\cb_w)$ such
that $A'\dsv_\G f_w^*A$. We have
$$A'\dsv_\G f_w^*({}^pH^{\cdot}(\p_{w!}\bbq))={}^pH^{\cdot}(f_w^*\p_{w!}\bbq)
={}^pH^{\cdot}(\p'_!\bbq)$$
and using (iv) we see that $A'\cong\bbq[l(w)]$. Then there exists a non-zero morphism 
$f_w^*A@>>>\bbq[l(w)]$ in $\cm_\G(\ti\cb_w)$. Hence there exists a non-zero morphism 
$A@>>>f_{w!}\bbq[l(w)]$ in $\cm_\G(\cb_w)$. Now $f_{w!}\bbq[l(w)]=\op_{\cl}\cl_w[l(w)]$ 
where $\cl$ runs over the local systems in $\cs(T)^{wF}$ (up to 
isomorphism). We must have $A\cong\cl_w[l(w)]$ for some $\cl$ as above. This proves (v).

We prove (vi). We have $A\dsv_\G{}^pH^{\cdot}(i_{w!}({}^pH^{j'}(\p_{w!}\bbq)))$ for some 
$j'$. Hence there exists a simple object $A_1$ of $\cm_\G(\cb_w)$ such that 
$A_1\dsv_\G{}^pH^{j'}(\p_{w!}\bbq)$ and $A\dsv_\G{}^pH^{\cdot}(i_{w!}A_1)$. From (v) we see
that $A_1\in\Bbb S(\cb_w)$. From 7.3 we see that $A\in\Bbb S(\cb)$. 

\subhead 7.5\endsubhead
Let $J\sub\II$. Define $\p_J:X@>>>\cp_J$ by $(B',g)\m P_{B',J}$. This is compatible with 
the $\G$-actions where $\G$ acts on $X$ by $g_0:(B',g)\m(g_0B'g_0\i,gg_0\i)$. Hence 
${}^pH^j(\p_{J!}\bbq)$ has a $\G$-equivariant structure. We show:

\proclaim{Theorem 7.6}(a) If $A$ is a simple object of $\cm_\G(\cp_J)$ such that 
$A\dsv_\G{}^pH^{\cdot}(\p_{J!}\bbq)$ then $A\in\Bbb S(\cp_J)$.
\endproclaim
Note that $\p_J$ is a composition $\p'\p''$ where $\p'':X@>>>\cb$ is $(B',g)\m B'$ and
$\p'':\cb@>>>\cp_J$ is $B'\m P_{B',J}$. We have a spectral sequence in $\cm_\G(\cp_J)$ with
$E_2={}^pH^{\cdot}(\p'_!({}^pH^{\cdot}(\p''_!\bbq)))$ and $E_\iy$ is an associated graded
of ${}^pH^{\cdot}(\p'_!\p''_!\bbq)$. We have $A\dsv_\G E_\iy$ hence $A\dsv_\G E_2$. Hence 
we can find a simple object $A_1\in\cm_\G(\cb)$ such that
$A_1\dsv_\G{}^pH^{\cdot}(\p''_!\bbq)$ and $A\dsv_\G{}^pH^{\cdot}(\p'_!A_1)$. Using the 
partition $\cb=\sqc_w\cb_w$ we see that $A_1\dsv_\G{}^pH^{\cdot}(i_{w!}i_w^*\p''_!\bbq)$ 
for some $w\in W$ (with $i_w$ as in 7.3). Since $i_w^*\p''_!\bbq=\p_{w!}\bbq$ (with $\p_w$
as in 7.4) we see that $A_1\dsv_\G{}^pH^{\cdot}(i_{w!}\p_{w!}\bbq)$. Using 7.4(vi) we see
that $A_1\in\Bbb S(\cb_w)$. By 4.13 (for $\cb$ instead of $\cp_J$) we see that
$A_1\dsv_\G{}^pH^{\cdot}({}^\em\bUp^\ss_!\bcl^\sh)$ for some sequence $\ss$ in $\II$ and 
some $\cl\in\cs(T)^{[\ss]F}$ (here ${}^\em\bUp^\ss$ is $\bUp^\ss$ of 3.2 with $J$ replaced
by $\em$). We have a spectral sequence in $\cm_\G(\cp_J)$ with
$E_2={}^pH^{\cdot}(\p'_!({}^pH^{\cdot}({}^\em\bUp^\ss_!\bcl^\sh)))$ and $E_\iy$ is an
associated graded of 

${}^pH^{\cdot}(\p'_!{}^\em\bUp^\ss_!\bcl^\sh)={}^pH^{\cdot}(\bUp^\ss_!\bcl^\sh)$
\nl
(with $\bUp^\ss$ as in 3.2). Now $A_1$ is a direct summand of 
${}^pH^{\cdot}({}^\em\bUp^\ss_!\bcl^\sh)$. Hence ${}^pH^{\cdot}(\p'_!A_1)$ is a direct 
summand of ${}^pH^{\cdot}(\p'_!({}^pH^{\cdot}({}^\em\bUp^\ss_!\bcl^\sh)))$. Hence 
$A\dsv_\G E_2$. Our spectral sequence is degenerate (by an argument as in 4.8). Hence 
$A\dsv_\G E_\iy$ and $A\dsv_\G{}^pH^{\cdot}(\bUp^\ss_!\bcl^\sh)$. Using 4.13 we deduce that
$A\in\Bbb S(\cp_J)$. The theorem is proved.

\subhead 7.7\endsubhead
We have a $\G$-action $g_1:(B',g)\m(B',g_1g)$ on $X$ (this is different from the 
$\G$-action on $X$ in 7.5). This induces a $\G$-module structure on $H^{2d}_c(X,\bbq)$ ($d$
as in 7.1).

\proclaim{Proposition 7.8} The $\G$-module $H^{2d}_c(X,\bbq)$ contains a copy of
$\text{Reg}$, the left regular representation of $\G$.
\endproclaim
Let $\r_1:X@>>>X'$ be as in 7.1. Let $X''=\{u\in G;u\text{ unipotent}\}$. Define 
$\r_2:X'@>>>X''$ by $(B',u)\m u$. Let $\r=\r_2\r_1:X@>>>X''$. It is well known that $\r_2$
is a semismall morphism. Recall that $\r_1$ is a finite \'etale covering. Hence $\r$ is a 
semismall morphism. Using this, we see that $\r_!\bbq[2d]$ is a perverse sheaf on $X'$ 
(recall that $X$ is smooth of pure dimension $2d$ and $\r$ is proper). By the decomposition
theorem \cite{\BBD}, ${}^pH^{\cdot}(\r_!\bbq)$ is a semisimple perverse sheaf. Hence 
$\r_!\bbq[2d]$ is a semisimple perverse sheaf. Now $\G$ acts on $X''$ trivially and $\r$ is
compatible with the $\G$-actions. Since $\bbq$ is naturally a $\G$-equivariant local system
on $X$ we see that we have naturally $\r_!\bbq[2d]\in\cm_\G(X'')$. Let $j:\{1\}@>>>X''$ be
the inclusion of the unit element into $X''$. If $E$ is an irreducible $\G$-module over 
$\bbq$ we can regard $E$ as a $\G$-equivariant local system on $\{1\}$. Then $j_!E$ is a 
simple object of $\cm_\G(X'')$. Let $n_E$ be the number of times $j_!E$ appears in a direct
sum decomposition of $\r_!\bbq[2d]$ (a semisimple object of $\cm_\G(X'')$) into a direct 
sum of simple objects. Since $\r_!\bbq[2d]$ is selfdual, we have using 1.6(a):
$$\align&n_E=\dim H^0_c(X'',j_!E\ot\r_!\bbq[2d])^\G
=\dim H^{2d}_c(\{1\},E\ot j^*\r_!\bbq)^\G\\&=
\dim(H^{2d}_c(\cb\T\G,\bbq)\ot E)^\G=\dim(\text{Reg}(-d)\ot E)^\G=\dim E.\endalign$$
We see that $\r_!\bbq[2d]$ contains as a direct summand the perverse sheaf $j_!\text{Reg}$ 
where $Reg$ is regarded as an object of $\cm_\G(\{1\})$. It follows that the $\G$-module 
$H^0_c(\r_!\bbq[2d])$ contains as a direct summand the $\G$-module 
$H^0_c(j_!\text{Reg})=\text{Reg}$. Equivalently the $\G$-module $H^{2d}(X,\bbq)$ contains 
$\text{Reg}$ as a direct summand. The proposition is proved.

\proclaim{Corollary 7.9}Let $E$ be an irreducible $\G$-module over $\bbq$. There exists 
$w\in W$ such that $E$ appears in the $\G$-module $\op_iH^i_c(\ti\cb_w,\bbq)$. ($\G$ acts 
on $\ti\cb_w$ by $g_1:zU@>>>g_1zU$.)
\endproclaim
From 7.8 we see that $E$ appears in the $\G$-module $H^{2d}_c(X,\bbq)$. Using this and the
partition $X=\sqc_{w\in W}X_w$ (with $X_w$ as in 7.4, $\G$-stable) we see that there exists
$w\in W$ such that $E$ appears in the $\G$-module $H^{2d}_c(X_w,\bbq)$. 

Now $\G$ acts:

on $\tX_w$ by $g_1:(zU_w,g)\m(zU_w,g_1g)$, 

on $\tX'_w$ by $g_1:(z,y)\m(z,g_1y)$, 

on $\tX''_w$ by $g_1:(z,y)U_w\m(z,g_1y)U_w$, 

on $\hat\cb_w$ by $g_1:z\m g_1z$. 
\nl
(Notation of 7.4.) Moreover the maps $\g,a,a'$ in 7.4 are compatible with
the $\G$-actions. Since $\g$ is a finite principal covering we see that $E$ appears in the
$\G$-module $H^{2d}_c(\tX_w,\bbq)$. We identify $\tX_w=\tX''_w$ as in 7.4. We see that $E$
appears in the $\G$-module $H^{2d}_c(\tX''_w,\bbq)$. Since $a$ is an affine space bundle we
see that $E$ appears in the $\G$-module $H^{4d-2l(w)}_c(\tX'_w,\bbq)$. We have 
$\tX'_w=\hat\cb_w\T\hat\cb_w$ and
$\op_iH^i_c(\tX'_w,\bbq)=\op_{i,i'}H^i_c(\hat\cb_w,\bbq)\ot H^{i'}_c(\hat\cb_w,\bbq)$ with
$\G$ acting only on the seond factor. It follows that $E$ appears in the $\G$-module 
$\op_iH^i_c(\hat\cb_w,\bbq)$. Since $a'$ is an affine space bundle we see that the 
corollary holds.

\subhead 7.10\endsubhead
For $w\in W$ let $\fT_w=\{t\in T;\dw F(t)\dw\i=t\}$. Now $\G\T\G\T\fT_w$ acts on $X_w$ by 
$(g_0,g_1,t):(B',g)\m(g_0B'g_0\i,g_1gg_0\i)$, on $\tX_w$ by 
$(g_0,g_1,t):(zU_w,g)\m(g_0zt\i U_w,g_1gg_0\i)$, on $\tX'_w$ by 
$(g_0,g_1,t):(z,y)\m(g_0zt\i,g_1yt\i)$, on $\tX''_w$ by 
$(g_0,g_1,t):(z,y)U_w\m(g_0zt\i,g_1yt\i)U_w$.
Moreover the maps $\g,a$ in 7.4 are compatible with the $\G\T\G\T\fT_w$-actions. Also 
$\G\T\fT_w$ acts on $\hat\cb_w$ by $(g_1,t):z\m g_1zt\i$, on $\ti\cb_w$ by 
$(g_1,t):zU\m g_1zt\i U$; the map $a'$ in 7.4 is compatible with the $\G\T\fT_w$-actions. 

For any $\th\in\hat\fT_w:=\Hom(\fT_w,\bbq^*)$ let $H^i_c(\ti\cb_w)_\th$ be the subspace of
$H^i_c(\ti\cb_w)$ on which $\fT_w$ acts via $\th$; this is naturaly a $\G\T\fT_w$-module.
Now $\G\T\G$ acts on $X$ by $(g_0,g_1):(B',g)\m(g_0B'g_0\i,g_1gg_0\i)$.

Let $\cg(\G)$ (resp. $\cg(\G\T\G)$ or $\cg(\G\T\G\T\fT_w)$) be the Grothendieck group of 
$\G$-modules (resp $\G\T\G$-modules or $\G\T\G\T\fT_w$-modules) of finite dimension over 
$\bbq$. Let $\Pi:\cg(\G\T\G\T\fT_w)@>>>\cg(\G\T\G\T\fT_w)$ be the homomorphism which takes
an irreducible $\G\T\G\T\fT_w$-module to the space of $\fT_w$-invariants (an irreducible 
$\G\T\G\T\fT_w$-module or $0$). In the setup of 7.9 we show:

\proclaim{Proposition 7.11}
$$\sum_i(-1)^iH^i_c(X,\bbq)=\sum\Sb w\in W\\ \th\in\hat\fT_w\\i,i'\eSb
(-1)^{i+i'}H^i_c(\ti\cb_w,\bbq)_{\th\i}\ot H^{i'}_c(\ti\cb_w,\bbq)_{\th}.\tag a$$
equality in $\cg(\G\T\G)$.
\endproclaim
By the arguments in the proof of 7.9 we have
$$\sum_i(-1)^iH^i_c(X,\bbq)=\sum_{w\in W}\sum_i(-1)^iH^i_c(X_w,\bbq),$$
equality in $\cg(\G\T\G)$ and
$$\align&\sum_i(-1)^iH^i_c(X_w,\bbq)=\Pi(\sum_i(-1)^iH^i_c(\tX_w,\bbq))=
\Pi(\sum_i(-1)^iH^i_c(\tX''_w,\bbq))\\&=\Pi(\sum_i(-1)^iH^i_c(\tX'_w,\bbq))=
\Pi(\sum_{i,i'}(-1)^{i+i'}H^i_c(\hat\cb_w,\bbq)\ot H^{i'}_c(\hat\cb_w,\bbq))\\&
=\sum_{\th\in\hat\fT_w}\sum_{i,i'}
(-1)^{i+i'}H^i_c(\hat\cb_w,\bbq)_{\th\i}\ot H^{i'}_c(\hat\cb_w,\bbq)_{\th}\\&=
\sum_{\th\in\hat\fT_w}\sum_{i,i'}
(-1)^{i+i'}H^i_c(\ti\cb_w,\bbq)_{\th\i}\ot H^{i'}_c(\ti\cb_w,\bbq)_{\th},
\endalign$$
equalities in $\cg(\G\T\G\T\fT_w)$. This proves (a).

\subhead 7.12\endsubhead
Restricting the $\G\T\G$-action on the modules in 7.11(a) to $\G$ by $g_1\m(1,g_1)$ gives
$$\sum_i(-1)^iH^i_c(X,\bbq)=\sum\Sb w\in W\\ \th\in\hat\fT_w\\ i,i'\eSb
(-1)^{i+i'}\dim(H^i_c(\ti\cb_w,\bbq)_{\th\i})H^{i'}_c(\ti\cb_w,\bbq)_{\th}.$$
equality in $\cg(\G)$. From \cite{\DL, 7.1, (7.6.3)} we see that the right hand side is 
equal to $|W|\text{Reg}$. Thus we have

\proclaim{Proposition 7.13}
$$\sum_i(-1)^iH^i_c(X,\bbq)=|W|\text{Reg}  \text{ in } \cg(\G).\tag a$$
\endproclaim

\proclaim{Lemma 7.14} Assume that $G$ is simply connected. Let $\ss=(s_1,\do,s_r),\cl$ be 
as in 2.6 and that 6.18(a) holds. Let $Z^\ss$ be as in 2.5. Let $\dZ$ be as in 2.5 (with 
$\ww=\ss$). Let $Z=\cz^\ss$ be as in 2.6. Then:

(a) $\dim H^{2r}_c(Z,\bcl)=n_\cl$;

(b) $\dim H^{2r}_c(Z^\ss,\cl_\ss)=n_\cl$;

(c) $\dim H^{2r}_c(\dZ,\bbq)=1$.
\endproclaim
We prove (a). We shall use 6.18(b) assuming that $\tss=\ss$. Then $Z=Z'$ and $\dim Z=r$ 
hence $H^i_c(Z,\bcl)=0$ for $i<2r$ and $H^i_c(Z,\tbcl)=0$ for $i<2r$. Hence
$H^{4r}_(Z\T Z,\bcl\bxt\tbcl)=H^{2r}_c(Z,\bcl)\ot H^{2r}_c(Z,\tbcl)$. Since in this case 
$\r=2r$ we see from 6.18(b) that $\dim(H^{2r}_c(Z,\bcl)\ot H^{2r}_c(Z,\tbcl))^\G=n_\cl$. 
Using Poincar\'e duality on the smooth variety $Z$ of pure dimension $r$ (see 2.9) we 
deduce that $\dim(H^0(Z,\tbcl)\ot H^0(Z,\bcl))^\G=n_\cl$. Note that the $\G$-modules 
$H^0(Z,\tbcl)$, $H^0(Z,\bcl)$ are dual to each other. Hence the previous equality can be 
written as $\dim\End_\G(H^0(Z,\bcl))=n_\cl$ where $\End_\G()$ is the space of $\G$-module 
endomorphisms. If $H^0(Z,\bcl)=0$, it follows that $n_\l=0$; in this case we have also
$H^0(Z,\tbcl)=0$ and by Poicar\'e duality, $H^{2r}_c(Z,\bcl)=0$ and (a) holds. Thus we may
assume that $H^0(Z,\bcl)\ne0$. Then $\dim\End_\G(H^0(Z,\bcl))\ge1$ hence 
$\dim\End_\G(H^0(Z,\bcl))=n_\cl=1$. Then $R=R_\cl$. Since $G$ is simply connected, it 
follows that $\cl\cong\bbq$. Then $H^0(Z,\bcl)=H^0(Z,\bbq)$ may be identified with the
permutation representation $V$ of $\G$ on the set $C$ of connected components of $Z$. Let 
$\G'$ be the isotropy group in $\G$ of some connected component of $Z$. Since $\G$ acts 
transitively on $C$ (see 2.9) we see that $\dim V=|\G/\G'|$, 
$1=\dim\End_G(V)=|\G'\bsl\G/\G'|$. It follows that $\G=\G'$ hence $\dim V=1$. Then we have
$\dim H^0(Z,\tbcl)=1$ and, by Poincar\'e duality, $\dim H^{2r}_c(Z,\bcl)=1$. This proves
(a).

We prove (b). Note that $Z^\ss$ is an open dense subset of $Z$ (using the commutative 
diagram in 2.6 with $\cj=\em$ and 2.8, this statement is reduced to the statement that 
$Z_2^\em$ in 2.6 is open dense in $Z_2$ which is clear). Note also that 
$\bcl|_{Z^\ss}=\cl_\ss$ (see 2.12) and that $Z$ has pure dimension $r$ (see 2.9). We see 
that (b) follows from (a).

We prove (c). Let $f_!^\c\bbq$ be as in 2.5 with $\ww=\ss$. From the definitions we have 
$H^{2r}_c(\dZ,\bbq)=\op_\c H^{2r}_c(Z^\ss,f^\c_!\bbq)$ where $\c$ runs over 
$\Hom(T^{F'},\bbq^*)$ (as in 2.5). Using (b) and 2.5(a) we see that 
$\dim H^{2r}_c(Z^\ss,f^\c_!\bbq)$ is $1$ if $\c=1$ and is $0$ if $\c\ne1$. The result 
follows.

\subhead 7.15\endsubhead
Assume that $G$ is simply connected. Let $d$ be as in 7.1. Let $w=w_\II$. We show:

(a) {\it $\tcb_w$ is connected;}

(b) {\it $X_w$ is connected (notation of 7.4);}
\nl
Assume that $\ss$ is a reduced expression for $w$. The associated variety $\dZ$ (see 2.5)
is connected by 7.14(c) (note that $\dZ$ has pure dimension $d=l(w)$). But $\dZ$ may be 
identified with $\tcb_w$. Hence (a) holds. Since $a':\hat{\cb}_w@>>>\tcb_w$ is a principal
$U_w$-bundle (as in 7.4) we see that $\hat{\cb}_w$ is connected. Since 
$\tX'_w\cong\hat{\cb}_w\T\hat{\cb}_w$ (as in 7.4) we see that $\tX'_w$ is connected. Since
$a:\tX'_w@>>>\tX''_w$ (as in 7.4) is a principal $U_w$-bundle we see that $\tX''_w$ is
connected. Since $\tX_w=\tX''_w$ (as in 7.4) we see that $\tX_w$ is connected. Since 
$\g:\tX_w@>>>X_w$ (as in 7.4) is surjective, we see that $X_w$ is
connected. Hence (b) holds. 

\mpb

Note that (a),(b) above do not necessarily hold without the assumption that $G$ is
simply connected.

\subhead 7.16\endsubhead
{\it Proof of Proposition 7.2.} Note that $X=\cup_{w'\in W}X_{w'}$, that $X$ is of pure 
dimension $2d$, that $X_w$ is an open subset of $X$ and that for any $w'\in W-\{w\}$, 
$X_{w'}$ has pure dimension equal to $d+l(w')<2d$. It follows that $X$ is connected if and
only if $X_w$ is connected. Hence 7.2 is a consequence of 7.15(b).

\head 8. A conjecture\endhead
\subhead 8.1\endsubhead
In this section we assume that $G$ has connected centre. Let $\cl\in\cs(T)$ be such that
$(wF)^*\cl\cong\cl$ for some $w\in W$. Our assumption on $G$ guarantees that 

(a) $w\in W'_\cl\imp w\in W_\cl$.
\nl
(Notation of 5.7.) Let $\cx_\cl$ be the set of all sequences $\ss=(s_1,\do,s_r)$ in $\II$ 
such that $\cl\in\cs(T)^{[\ss]F}$. Note that $\cx_\cl\ne\em$.

To any $\ss\in\cx_\cl$ we associate an element $\o\in W$ as in 5.5. We show that $\o$ is 
independent of the choice of $\ss$. We must show that if $\tss$ is another element of 
$\cx_\cl$ and $\tio\in W$ is associated to $\tss$ in the same way as $\o$ is associated to
$\ss$ then $\tio=\o$. Using 5.6 we see that $F^*\o^*\cl\cong\cl$, $F^*\tio^*\cl\cong\cl$.
Hence $\o^*\cl\cong\tio^*\cl$ so that $\tio\o\i\in W'_\cl$ and, using (a), 
$\tio\o\i\in W_\cl$. By the proof of 5.8 (with $c=\boc$) we have $\o\boc(R^+_\cl)=R^+_\cl$ and
similarly $\tio\boc(R^+_\cl)=R^+_\cl$. Thus $\o\i(R^+_\cl)=\tio\i(R^+_\cl$ and
$\tio\o\i(R^+_\cl)=R^+_\cl$. This together with $\tio\o\i\in W_\cl$ yields $\tio\i\o=1$, as
desired.

Using the previous paragraph and (a) we see that for any $\ss,\tss\in\cx_\cl$, the set 
$\fF$ defined as in 5.9 (with $c=\boc$) is equal to $\{1\}$.

\subhead 8.2\endsubhead
For $\ss,\tss\in\cx_\cl$ let $Z,\tZ,r,\y$ be as in 6.1. Let $\bcl$ be as in 6.2 and let 
$\tbcl$ be the analogous local system on $\tZ$ defined in terms of $\cl$. Then 
$\check{\tbcl}$ (dual of $\tbcl$) is the analogous local system on $\tZ$ defined in terms
of $\ccl$. Let 
$$V_{\tss,\ss}=\Hom_\G(\op_iH^i_c(\tZ,\tbcl)(i/2),\op_{i'}H^{i'}_c(Z,\bcl)(i'/2)).$$
where $\Hom_\G$ is the space of homomorphisms of $\G$-modules. Using 2.13(a) we see that
Poincar\'e duality holds on $\tZ$ in the form
$$\Hom(H^i_c(\tZ,\tbcl)(i/2),\bbq)=H^{2\y-i}_c(\tZ,\check{\tbcl})(\y-i/2).$$
Hence
$$\Hom(\op_iH^i_c(\tZ,\tbcl)(i/2),\bbq)=\op_iH^i_c(\tZ,\check{\tbcl})(i/2),$$
$$V_{\ss,\tss}
=((\op_{i'}H^{i'}_c(Z,\bcl)(i'/2))\ot(\op_iH^i_c(\tZ,\check{\tbcl})(i/2)))^\G,$$
$$V_{\ss,\tss}=(\op_nH^n_c(Z\T\tZ,\bcl\bxt\check{\tbcl})(n/2))^\G.$$
By 6.14(d) and 6.15(b), the last vector space has a distinguished basis 
$\{b_\aa^1;\aa\in\ca_1\}$, with $\ca_1$ as in 6.14.

Let $C_\cl$ be the category whose objects are the elements of $\cx_\cl$ and in which the 
set of morphisms from $\tss$ to $\ss$ is the vector space $V_{\tss,\ss}$. The composition
of morphisms is given by composing linear maps.

\subhead 8.3\endsubhead
We will view $T$ as a maximally $\FF_q$-split torus of a second connected reductive 
algebraic group $G'$ over $\FF_q$ in such a way that $R_\cl$ is the 
set of roots of $G'$ with respect to $T'$ and $R^+_\cl$ is the set of positive roots of 
$G'$ with respect to $T'$ and a Borel subgroup $B'$ of $G'$ which is defined over $\FF_q$ 
and contains $T'$.

Replacing $G,T,B,\cl$ by $G',T,B',\bbq$ in the definition of the set $\cx_\cl$ and of the 
category $C_\cl$ in 8.2 we obtain a set $\cx'_{\bbq}$ and a category $C'_{\bbq}$. Note that
the objects of $C'_{\bbq}$ are the elements of $\cx'_{\bbq}$ that is the sequences
$\SS=(S_1,S_2,\do,S_b)$ in $W_\cl$. For $\tSS,\SS$ in $\cx'_{\bbq}$ we denote by 
$V'_{\tSS,\SS}$ the vector space of morphisms from $\tSS$ to $\SS$ in $C'_{\bbq}$.

The following is conjecturally a functor $\Ph:C_\cl@>>>C'_{\bbq}$. To an object $\ss$ of
$C_\cl$, $\Ph$ associates the object $\SS$ of $C'_{\bbq}$ defined as in 5.5. Given two 
objects $\tss,\ss$ of $C_\cl$ we set $\tSS=\Ph(\tss),\SS=\Ph(\ss)$ and we define a linear 
map $\Ph:V_{\tss,\ss}@>>>V'_{\tSS,\SS}$ to be the isomorphism which maps the distinguished
basis of $V_{\tss,\ss}$ onto the analogous distinguished basis of $V'_{\tSS,\SS}$ according
to the bijection $\ca(W,\boc,\cl,\ss,\tss)@>\sim>>\ca(W_\cl,\o\boc,\bbq,\SS,\tSS)$ described in
5.10. (As pointed out in 8.1, the set $\fF$ which appears in 5.10 is in our case equal to 
$\{1\}$.) We expect that $\Ph$ is a functor and that moreover it is an equivalence of
categories.

\head Index of Notation\endhead

1.1. $\kk,G,\cb,B,T,N(T),W,\po,l(w),\II$

1.2. $\le,{}^JW,W^{J'},{}^JW^{J'},w_\II$

1.3. $\cp,U_P,U,\cp_J,L_J,P_{B',J},P^Q,k(g)$

1.4, $R,\chR,U_\a,R^+,R^-,\a_s$

1.5. $x_s(),y_s(),\dw,[\ww],[\ww]^\bul$

1.6. $\cd(X),\fD(K),\cm(X),{}^pH^{\cdot}(K),f^\bst(A),\cm_\G(Y),\cf^\sharp,K\dsv_\G K',E_X$

1.7. $\cs(T)$

1.10. $R_\cl,R^+_\cl,W_\cl,\II_\cl,\chR_\cl$

2.4. $\ciss,\uucl,F:G>>G,\G,\cs(T)^{wF},F_0:T@>>>T,\boc$

2.5. $Z^\ww,\fT,\cb_w,\cl_\ww,\cl_w$

2.6. $\bZ^\ss,\cz^\ss$

2.11. $\bcl$

3.2. $\Pi^\ww,\Up^\ww,\bar\Up^\ss$

3.7. $\Bbb S'(\cp_J)$

4.1 $\cp_J^\tt$

4.2. $\ct'(J,\boc),\vt$

4.3. $\tcp_J^\tt$

4.4. $\uM,\Bbb S(\cp_J^\tt),\Bbb S(\cp_J)$

5.6. $W'_\cl$

5.9. $\ca(W,c,\cl,\ss,\tss)$

6.3 $\ct$

6.9. $N_\aa$

7.1. $X$


\widestnumber\key{ABC}
\Refs
\ref\key{\BE}\by R.B\'edard\paper On the Brauer liftings for modular representations\jour 
J.Algebra\vol93\yr1985\pages332-353\endref
\ref\key{\BBD}\by A.Beilinson, J.Bernstein and P.Deligne\paper Faisceaux pervers\jour 
Ast\'erisque\vol100\yr1982\pages5-171\endref
\ref\key{\BO}\by A.Borel\book Linear algebraic groups\publ Benjamin,\yr1969\publaddr
New York, Amsterdam\endref
\ref\key{\DL}\by P.Deligne and G.Lusztig\paper Representations of reductive groups over 
finite fields\jour Ann. Math.\vol103\yr1976\pages103-161\endref
\ref\key{\ORA}\by G.Lusztig\book Characters of reductive groups over a finite field, 
Ann.Math.Studies 107\publ Princeton U.Press\yr1984\endref
\ref\key{\CS}\by G.Lusztig\paper Character sheaves I\jour Adv.Math.\vol56\yr1985\pages
193-237\moreref II\jour Adv.Math.\vol57\yr1985\pages226-265\endref
\ref\key{\GR}\by G.Lusztig\paper Green functions and character sheaves\jour Ann.Math.\vol
131\yr1990\pages355-408\endref
\ref\key{\HB}\by G.Lusztig\paper Homology bases arising from reductive groups over a finite
field, in "Algebraic groups and their representations", ed. R.W.Carter et al.\publ Kluwer 
Acad. Publ.\yr1998\pages53-72\endref
\ref\key{\PAR}\by G.Lusztig\paper Parabolic character sheaves I, II\jour Moscow Math.J.\vol
4\yr2004\pages153-179, 869-896\endref
\ref\key{\CSD}\by G.Lusztig\paper Character sheaves on disconnected groups, VI\jour 
Represent.Th.(electronic)\vol8\yr2004\pages377-413\endref 
\endRefs
\enddocument